\documentclass[11pt]{amsart}
\usepackage[margin=1.15in]{geometry}
\usepackage{amsmath,amssymb,amsthm,mathrsfs}
\usepackage{enumitem}
\usepackage[hidelinks]{hyperref}
\numberwithin{equation}{section}
\theoremstyle{plain}
\newtheorem{theorem}{Theorem}[section]
\newtheorem{proposition}[theorem]{Proposition}
\newtheorem{corollary}[theorem]{Corollary}
\newtheorem{lemma}[theorem]{Lemma}
\theoremstyle{definition}

\newtheorem{example}[theorem]{Example}
\theoremstyle{remark}
\newtheorem{remark}[theorem]{Remark}
\newcommand{\C}{\mathbb C}
\newcommand{\R}{\mathbb R}
\newcommand{\Z}{\mathbb Z}
\newcommand{\cO}{\mathcal O}

\newcommand{\tr}{\operatorname{tr}}
\newcommand{\Sym}{\operatorname{Sym}}
\newcommand{\im}{\operatorname{im}}

\title[Algebraicity and geometry of period images]
{Algebraicity, Hodge Gradings, and the Geometry of Period Images}

\author{Kefeng Liu}
\address{Mathematical Sciences Research Center, Chongqing University of Technology, Chongqing 400054, China;\\
Department of Mathematics, University of California at Los Angeles, Los Angeles, CA 90095-1555, USA}
\email{liu@math.ucla.edu}
\author{Yang Shen}
\address{Mathematical Sciences Research Center, Chongqing University of Technology, Chongqing 400054, China}
\email{syliuguang2007@163.com}
\date{}

\begin{document}
\begin{abstract}
We prove the algebraicity of period images in arithmetic quotients
without using o-minimality. The proof uses infinitesimal neighborhoods
and formal algebraization to recover the images and the maps from
their algebraic sources.
Hodge positivity and an algebraic ampleness criterion yield the
Griffiths polarization and quasi-projectivity. For lifted images,
we construct Stein neighborhoods and holomorphic matrix realizations
of the Hodge filtration. On closed horizontal subspaces, we bound
canonical Bergman density below by the Hodge volume form. We compare
the Bergman and Griffiths--Schmid metrics near one-parameter boundaries
with nontrivial unipotent monodromy and nondegenerate graded period maps,
and in all directions near general one- and two-node quintic threefolds.
A period surface arising from abelian threefolds shows that the upper
metric comparison fails in general.
\end{abstract}
\maketitle
\setcounter{tocdepth}{1}
\tableofcontents

\section{Introduction}

Let $D$ be a connected period domain for polarized integral
Hodge structures of weight $w$, and let $\Gamma$ be an
arithmetic subgroup of the polarization group preserving $D$.
A period map to $\Gamma\backslash D$ is a
locally liftable holomorphic map whose differential is horizontal
on the regular locus.

If $a$ and $b$ are the smallest and
largest indices with nonzero Hodge numbers, the Griffiths line bundle is
$$
 \Lambda=\bigotimes_{p=a+1}^{b}\det F^p.
$$
On a neat quotient this is a line bundle; in general it is a
$\mathbb Q$-line bundle.

We study the algebraic structure of period images, including at
nonnormal points, and the geometry of their lifts to $D$.
Holomorphic matrix realizations of the Hodge filtration provide
potentials for the Hodge metric and estimates for the canonical
Bergman kernel.

\subsection*{Algebraicity and the Griffiths polarization}

\begin{theorem}[Algebraicity and polarization]\label{intro:algebraicity}
Let $B$ be a reduced separated algebraic space of finite type
over $\C$, and let $\Phi:B^{\rm an}\to\Gamma\backslash D$ be a period map
defined on all of $B^{\rm an}$. There are a reduced
quasi-projective scheme $Y$, a dominant algebraic map $f$,
and a closed analytic immersion $i$, such that
\begin{equation}\label{intro:algebraic-factor}
 \Phi=i\circ f^{\rm an},\qquad
 B\xrightarrow{f}Y,\qquad
 i:Y^{\rm an}\hookrightarrow\Gamma\backslash D,\qquad
 i(Y^{\rm an})=\overline{\Phi(B^{\rm an})}_{\rm red}.
\end{equation}
Here the closure is taken in $\Gamma\backslash D$, with its
reduced analytic structure. The pullback $i^*\Lambda$ is the
analytification of an ample
algebraic $\mathbb Q$-line bundle $L$ on $Y$.
If $\Phi$ is proper, then $f$ is proper and surjective.
For neat $\Gamma$, the Hodge bundle, its filtration and its
polarization have algebraic models on $Y$, compatible with
pullback to every smooth algebraic source.
\end{theorem}

We first prove the neat case: Sections~\ref{aq:section}--\ref{aq:formal-section}
establish algebraicity, and Section~\ref{aq:griffiths-section}
algebraizes the Hodge bundles and proves ampleness.
For general $\Gamma$, choose a normal neat subgroup
$\Gamma'$ of finite index and a proper surjection $S\to B$
from a finite disjoint union of smooth quasi-projective varieties.
Normalizing $(S^{\rm an}\times_{\Gamma\backslash D}(\Gamma'\backslash D))_{\rm red}$
gives a finite \'etale cover of $S^{\rm an}$.
Riemann existence algebraizes this cover, so the neat theorem applies upstairs;
the finite quotient gives the image, the source map, and the
Griffiths $\mathbb Q$-line bundle.
This descent is proved in Theorem~\ref{aq:finite-level}.

This theorem was proved by Bakker--Brunebarbe--Tsimerman
\cite[Theorem~1.1]{BBT}. Our proof constructs the infinitesimal
images of period maps using their algebraic flat connections,
and recovers the analytic image by formal algebraization.
Once the image and the Hodge bundles have been algebraized, their
canonical extensions satisfy the general ampleness criterion,
Theorem~\ref{aq:general-ampleness}. In particular,
\begin{equation}\label{intro:ample-closure}
 j:Y\hookrightarrow\mathbf P(V^\vee),\qquad
 j^*\cO(1)\simeq L^{m_0},\qquad
 \cO_{\overline{j(Y)}}(1)|_Y\simeq L^{m_0}
\end{equation}
for a sufficiently divisible $m_0>0$ and a finite system
$V\subset H^0(Y,L^{m_0})$.
In general, $f(B)$ is a dense constructible subset of $Y$;
when $\Phi$ is proper, it equals $Y$. The assumption that $B$ is reduced
is essential: arbitrary analytic nilpotent deformations of a period map
need not algebraize
\cite[Example~6.2]{BBT}.

Here is an outline of the proof. Griffiths' extension theorem reduces
the case of a smooth source to a proper period map. Cattani--Deligne--Kaplan
give the algebraicity of the reduced period relation. A Hilbert-scheme
construction gives a smooth quasi-projective modification
$$
 p:M^{\rm an}\to Y_{\rm an},\qquad
 W=(p^{-1}Z)_{\rm red}\subset M,
$$
with algebraic exceptional set and $\dim Z<\dim Y_{\rm an}$.
After algebraizing $Z$ by induction, put
$$
 W_n=V(\mathcal I_W^{n+1}),\qquad
 T_n=\operatorname{Im}_{\rm sch}(W_n^{\rm an}\to Y_{\rm an}).
$$
An algebraic integrable connection gives finite parallel transport
on infinitesimal neighborhoods. A transverse slice to the Hodge
section in a relative flag variety algebraizes the maps
$W_n\to T_n$, including their nonreduced structures and
transition maps.
At singular centers, a coherent submodule is recovered from its
restriction off the center and its restrictions to all finite
neighborhoods (Lemma~\ref{aq:coherent-recovery}). Applied also to
graphs of functions, this recovers the required subrings.
Artin's theorem then contracts the formal completion of $M$ along $W$.
Two applications of Krull's intersection theorem in the proper direct
image identify the resulting analytic local rings with those
of $Y_{\rm an}$, including at nonnormal points. The same
parallel transport and coherent sheaf arguments algebraize the
Hodge bundles with their given analytic identifications. Their canonical
extensions are compatible with pullback. Hodge positivity makes these
extensions nef and big on resolutions of closed subspaces, as required
by the ampleness criterion.

\subsection*{A Euclidean realization of the lifted flags}

Let $G$ be the connected complex polarization group,
$\mathfrak g$ its Lie algebra and $\check D$ the compact dual.
At a reference Hodge decomposition define
$$
 S_o|_{H_o^{p,w-p}}=(2p-w)I,
 \qquad \mathcal A=G\cdot S_o\subset\mathfrak g.
$$
The ordered spectral filtration defines a map
$\varpi:\mathcal A\to\check D$ from this closed affine orbit.

\begin{theorem}[Hodge-grading realization]\label{intro:realization}
For a polarized integral variation of pure Hodge structure over a
smooth connected
quasi-projective base, let $I\subset D$ be its lifted period
image and $Z=\overline I^{\,D}_{\rm red}$.
Then $Z$ is a closed analytic subspace of $D$. It is Stein, and its
regular locus is horizontal.
Inside every prescribed open neighborhood of $Z$ there is
an open Stein neighborhood $\Omega\subset D$ and a
holomorphic embedding
$$
 S:\Omega\hookrightarrow\mathcal A\subset\mathfrak g,
 \qquad \varpi\circ S=\operatorname{id}_{\Omega}.
$$
The spectral projections are
\begin{equation}\label{intro:spectral-recovery}
 \Pi_p(S)=\prod_{\substack{q\ne p\\H^{q,w-q}\ne0}}
       \frac{S-(2q-w)I}{2(p-q)},
 \qquad
 F^p=\bigoplus_{q\ge p}\operatorname{im}\Pi_q(S).
\end{equation}
Thus the restriction of $S$ to $I$ preserves its flags
and subspace topology.
\end{theorem}

No properness or Torelli assumption on the original period map is
required. The spectral projections recover every step of the Hodge
filtration.
The compatible gradings form a torsor under $T^*\check D$; Cartan's
Theorem B supplies a section on $\Omega$. On the algebraic quotient,
the descended Hodge bundle has algebraic gradings compatible with
its filtration on affine charts. These are sections of its endomorphism
bundle satisfying the same spectral identities.

The Stein neighborhood is obtained from Griffiths' extension
theorem, Remmert's proper mapping theorem and Siu's neighborhood
theorem. Boundedness of the matrix coordinates, as required in
\cite[Conjecture~(10.1)]{Griffiths1970}, remains open.

\subsection*{Bergman geometry near the boundary}

For a closed reduced pure $n$-dimensional horizontal space
$Z\subset D$, consider the Hilbert space of all square-integrable
holomorphic top forms
$$
 \mathcal H_Z=H^0_{(2)}(Z_{\rm reg},K_{Z_{\rm reg}}),\qquad
 \|\sigma\|^2=\int_{Z_{\rm reg}}i^{n^2}
                         \sigma\wedge\overline\sigma.
$$
Write an orthonormal basis of $\mathcal H_Z$ locally as
$\sigma_\nu=f_\nu\,dz_1\wedge\cdots\wedge dz_n$, and put
$$
 b_Z=\sum_\nu|f_\nu|^2,\qquad
 g_{B,Z}=\partial\bar\partial\log b_Z,
 \qquad
 dV_{B,Z}=i^{n^2}\sum_\nu\sigma_\nu\wedge\overline{\sigma_\nu}.
$$
The canonical density $dV_{B,Z}$ need not equal the Riemannian
volume of $g_{B,Z}$.
Write $g_{\rm GS}$ for the Griffiths--Schmid metric and
$dV_{\rm GS}$ for its induced volume density on $Z_{\rm reg}$.

\begin{theorem}[Bergman density and boundary comparison]
\label{intro:bergman-boundary}
For fixed $D$ and $n$,
$$
 dV_{B,Z}\ge c(D,n)\,dV_{\rm GS}\quad\text{on }Z_{\rm reg}.
$$
Consider a full-dimensional admissible one-parameter period branch in
$Z$, with nonzero logarithm of unipotent monodromy, fixed
limiting mixed Hodge numbers, and jointly injective differentials
of the period maps on the weight-graded pieces along a compact
set of transverse parameters. At regular points of $Z$ on this
branch, uniformly near the boundary,
\begin{equation}\label{intro:bergman-cusp}
 g_{B,Z}\asymp g_{\rm GS}\asymp
 \frac{|dq|^2}{|q|^2(\log|q|)^2}+g_S.
\end{equation}
Here $g_S$ is a fixed smooth Hermitian metric on the transverse
parameter space. On the effective deformation spaces of general
quintic threefolds with one or two nodes, the comparison holds in
all $101$ directions; at two nodes
the model has two cusp summands, with constants independent of
the relative degeneration rates. Every fixed covariant derivative
of the Bergman curvature is bounded on these regions.
\end{theorem}

The boundary comparison holds for arbitrary weight and nilpotency
index, without assuming that the limiting mixed Hodge structure is
$\R$-split.

Let $H_\C$ be the underlying complex vector space, $h$ the
Hodge metric and $h_o$ a fixed reference Hodge metric. The proof
starts from the trace estimate
$$
 E=\operatorname{tr}(h_o^{-1}h),\qquad
 |\partial E|^2\le c_rE\,\partial\bar\partial E,
 \qquad c_r=\cos\frac{\pi}{r+1}<1,
$$
where $r=b-a+1$ is the length of the Hodge filtration, counting
indices between the two extremes. The exterior-power sum
$$
 \sum_{k=1}^{\dim H_\C-1}
      \log\operatorname{tr}_{\Lambda^kh_o}(\Lambda^kh)
$$
has Levi form uniformly comparable to the Griffiths--Schmid metric.
Its one-parameter asymptotics give embedded period charts of dimension
$n$ with uniform estimates. Bounded trace weights
give uniform $L^2$ estimates for interpolation of jets of top
forms. The $\bar\partial$-equation is solved on a smooth Stein
open subset, and $L^2$ removability extends the resulting forms
to elements of $\mathcal H_Z$.

The nondegeneracy hypothesis on the graded period maps is essential.
The product family of Legendre elliptic curves
$$
 E_t\times E_{1/2+ts}\times E_{1/2+s^2}
$$
has a closed smooth lifted period surface with complete Griffiths--Schmid
metric, but along $s=0$,
$$
 g_{\rm GS}(\partial_s,\overline{\partial_s})\asymp|t|^2,
 \qquad g_{B,Z}(\partial_s,\overline{\partial_s})\ge c>0
 \quad(t\to0).
$$
Thus $g_{B,Z}\le Cg_{\rm GS}$ fails in general.
The boundary comparison also holds on the monodromy cover of the
source, with its own Hilbert space of square-integrable top forms.

\subsection*{Relation to previous work}

Griffiths showed that period maps extend across boundary points with
finite local monodromy and that the resulting extension is proper
\cite[Theorem~9.5 and the proof of Theorem~9.6]{Griffiths1970III}.
In his 1970 survey of open problems, he conjectured that period images
are algebraic and lifted period images embed into bounded domains
in complex Euclidean space
\cite[Conjectures~(10.5) and~(10.1)]{Griffiths1970}.
Sommese treated smooth images and developed the analytic
adjoint-section approach \cite{Sommese1975,Sommese1978}.
Bakker--Brunebarbe--Tsimerman proved the general result using
definable coherent sheaves and Artin algebraization \cite{BBT}.

Our modification follows Sommese's Chow-space construction
\cite[Propositions~III and IV]{Sommese1978} and its Hilbert-scheme
form in \cite[Section~4.2 and Proposition~4.6]{BBT}; we also retain
formal contraction. The infinitesimal-image argument outlined above
replaces definable coherent GAGA.

The Stein property and plurisubharmonic potentials for the Hodge metric
were studied by Griffiths--Schmid and Eyssidieux
\cite{GriffithsSchmid1969,Eyssidieux1997}.
The boundary asymptotics use the nilpotent-orbit and $\mathrm{SL}_2$
theorems \cite{Schmid1973,CKS1986}.
Our boundary comparison adapts Zimmer's methods for domains with bounded
intrinsic geometry \cite{Zimmer2020}, using the full space of
square-integrable top forms on the regular locus of a period image.
It applies near general quintic threefolds with two ordinary double points.

The ample extension in \eqref{intro:ample-closure} concerns a
chosen projective closure; the canonical Baily--Borel extension
\cite[Theorem~1.2]{BFMT2025} and the conjectural projectivity of a
suitable generalized toroidal completion
\cite[Conjecture~2.12]{DengTsimerman2025} are stronger statements.

\subsection*{Historical remark}

This paper continues our program of recovering the algebraicity of
period images through standard algebraic models. Our withdrawn preprints
\href{https://arxiv.org/abs/1507.01860}{arXiv:1507.01860}
(2015) and
\href{https://arxiv.org/abs/1608.06029}{arXiv:1608.06029} (2016)
sought to use the bounded realization of the lifted period image
envisaged in \cite[Conjecture~(10.1)]{Griffiths1970}, followed by
monodromy descent and Riemann extension. Those arguments contained serious gaps,
including an incorrect claim to have proved that conjecture,
which is not established here.

We subsequently turned to the algebraicity of the reduced period relation
given by Cattani--Deligne--Kaplan, an approach also
discussed in \cite[p.~166 and Example~4.13]{BBT}. The infinitesimal-image
argument outlined above completes this approach.

AI tools helped us develop and refine this argument, check references,
and clarify its historical context. Responsibility for the mathematics
and its attribution remains with the authors.

\subsection*{Organization}

Sections~\ref{aq:section}--\ref{aq:formal-section} establish
algebraization from infinitesimal neighborhoods of period images.
Section~\ref{aq:griffiths-section} recovers the Hodge bundles and
proves ampleness. Sections~\ref{sec:grading}--\ref{stein:section}
construct the grading realization and its Stein neighborhood.
Section~\ref{sec:boundary-hodge} develops the exterior trace and
normalized boundary charts. The final sections prove the
canonical estimates and the comparison counterexample.

\section{Algebraic period relations and modifications}
\label{aq:section}

Fix a polarized integral lattice, its period domain $D$, and a neat
arithmetic group $\Gamma$ preserving a connected component of $D$.
Point stabilizers are compact and discrete, hence finite. Neatness
therefore makes $\Gamma\backslash D$ a manifold and
$D\to\Gamma\backslash D$ a covering.
In the next four sections, period maps take values in this fixed
quotient. We prove algebraicity in the neat case of
Theorem~\ref{intro:algebraicity}; after establishing ampleness in
Section~\ref{aq:griffiths-section}, we remove the neat hypothesis by
descent along a finite cover in Theorem~\ref{aq:finite-level}.

\begin{theorem}\label{aq:algebraization}
Let $B$ be a reduced separated algebraic space of finite type over
$\C$, and let $\Phi:B^{\mathrm{an}}\to\Gamma\backslash D$ be holomorphic
on all of $B^{\mathrm{an}}$ and satisfy Griffiths transversality on its
regular locus.
There are a reduced separated algebraic space $Y$ of finite type,
a dominant algebraic morphism $f:B\to Y$, and a closed analytic
immersion $\iota:Y^{\mathrm{an}}\hookrightarrow\Gamma\backslash D$ such that
$$
 \Phi=\iota f^{\mathrm{an}},\qquad
 \iota(Y^{\mathrm{an}})=
 \overline{\Phi(B^{\mathrm{an}})}^{\,\Gamma\backslash D}_{\mathrm{red}}.
$$
If $\Phi$ is proper, then $f$ is proper and surjective, and
$Y^{\mathrm{an}}$ is the reduced analytic image of $\Phi$.
Moreover, $Y$ admits a proper algebraic surjection $q:M\to Y$
from a finite disjoint union of smooth quasi-projective schemes,
and $\iota q^{\mathrm{an}}$ is a period map of polarized integral
variations on these schemes.
\end{theorem}

The image is algebraized with its possibly nonnormal local rings.
This is the algebraicity assertion of \cite[Theorem~1.1]{BBT}.
The proof uses Hilbert schemes and formal contraction
\cite{Sommese1978,Artin1970} with parallel transport on infinitesimal
neighborhoods and algebraization of specified analytic subrings.
It proceeds in three steps:
$$
 \begin{array}{c}
 M^{\mathrm{an}}\longrightarrow Y,\quad
 W^{\mathrm{an}}=(p^{-1}Z)_{\mathrm{red}}
 \quad\text{(Proposition~\ref{aq:modification})},\\[2pt]
 W_n\longrightarrow T_n=
 \im_{\mathrm{sch}}(W_n^{\mathrm{an}}\to Y)
 \quad\text{(Theorem~\ref{aq:finite-images})},\\[2pt]
 \{W_n\to T_n\}_n\longrightarrow
 \bigl(M\to Y_{\mathrm{alg}},\ Y_{\mathrm{alg}}^{\mathrm{an}}=Y\bigr)
 \quad\text{(Proposition~\ref{aq:formal-recovery})}.
 \end{array}
$$
An algebraic model includes the prescribed analytic identification
and the morphism from the source. We use the comparison theorem for
coherent direct images under proper morphisms and faithful
analytification \cite[Expos\'e XII, 1.1, 1.3.1, 4.2]{SGA1}
on affine \'etale charts.

\subsection{A smooth modification with algebraic exceptional locus}

\begin{lemma}[Descent through a proper source]\label{aq:graph-descent}
Let $a:E\to T$ be proper and schematically dominant, where $E$ is
a separated finite-type scheme and $T,S$ are separated finite-type
algebraic spaces. If $b:E\to S$ is algebraic and an analytic map
$v:T^{\mathrm{an}}\to S^{\mathrm{an}}$ satisfies
$va^{\mathrm{an}}=b^{\mathrm{an}}$, then $v$ is algebraic.
\end{lemma}

\begin{proof}
The map $(a,b):E\to T\times S$ is proper. Its scheme-theoretic image $G$
satisfies
$$
 \mathcal I_G=
 \ker(\cO_{T\times S}\to(a,b)_*\cO_E),\qquad
 G^{\mathrm{an}}=\operatorname{Graph}(v),
$$
by proper comparison; schematic dominance is preserved as well.
The projection $\pi:G\to T$ is proper \cite[Tag 0AGD]{Stacks}.
Its positive-dimensional fiber locus is empty, since it is closed
and has no complex closed points. Thus $G\to T$ is finite.
On an affine \'etale chart, faithfulness of analytification shows
that the kernel and cokernel of $\cO_T\to\pi_*\cO_G$ vanish.
Hence $\pi$ is an isomorphism, and its inverse followed by the second
projection gives $v$. This also proves uniqueness.
\end{proof}

\begin{lemma}[The reduced period relation]\label{aq:relation}
For smooth quasi-projective varieties $B_1,B_2$ with period maps,
the reduced relation
$$
 R_{12}^{\mathrm{an}}
 =(B_1^{\mathrm{an}}\times_{\Gamma\backslash D}B_2^{\mathrm{an}})_{\mathrm{red}}
$$
is algebraic. If the second period map is proper, then $R_{12}\to B_1$
is proper.
\end{lemma}

\begin{proof}
On $B_1\times B_2$ take the weight-zero integral variation
$\mathbb W=\operatorname{Hom}_{\Z}(\mathbb V_1,\mathbb V_2)$.
Its polarization is an integral multiple of
$q_H(A,B)=\tr(A^{\dagger_Q}B)$: with Weil operators $C_i$,
$$
 C_H(A)=C_2AC_1^{-1},\qquad
 q_H(C_HA,\overline A)=\tr(A^{*h}A)>0.
$$
Every polarization-preserving Hodge isometry $u$ therefore has the
same norm $Q_{\mathbb W}(u,u)=c\,\operatorname{rank}\mathbb V_i$.
The space of integral Hodge classes of bounded norm is algebraic
and finite over the base \cite[Theorem~1.1]{CDK}. Retain the isometry
itself in this space. The condition $u\in\Gamma$ is open and closed in the integral
local system and invariant under $u\mapsto\gamma_2u\gamma_1^{-1}$.
Its intersection with this space is an algebraic union of connected
components. Its reduced finite image is $R_{12}$.
The final assertion follows from analytic properness and
\cite[Expos\'e XII, 3.2(v)]{SGA1}.
\end{proof}

\begin{proposition}\label{aq:modification}
Suppose $B$ is a finite disjoint union of smooth quasi-projective
varieties with period maps, $\Phi$ is proper, and
$Y=\Phi(B^{\mathrm{an}})_{\mathrm{red}}$ has dimension $d$.
There exist a finite disjoint union of smooth quasi-projective
varieties $M$ with period maps, a proper map $p:M^{\mathrm{an}}\to Y$,
and a reduced closed algebraic subspace $W\subset M$ such that, for
$Z=p(W)_{\mathrm{red}}$,
\begin{equation}\label{aq:center}
 W^{\mathrm{an}}=(p^{-1}Z)_{\mathrm{red}},\qquad
 M^{\mathrm{an}}\setminus W^{\mathrm{an}}\simeq Y\setminus Z,
 \qquad \dim Z<d
\end{equation}
whenever $Z$ is nonempty.
\end{proposition}

\begin{proof}
First let $B,Y$ be irreducible. The algebraic Hodge filtration and
connection \cite[(4.13)(a)]{Schmid1973} give the algebraic Higgs map
$$
 \vartheta_B:T_B\longrightarrow
 \bigoplus_p\operatorname{Hom}(\operatorname{Gr}_F^p\mathcal V,
                              \operatorname{Gr}_F^{p-1}\mathcal V).
$$
Let $C$ be the locus where its rank is less than $d$, and let
$N_s,N_t$ be the nonsmooth loci of the proper projections
$s,t:R_{BB}\rightrightarrows B$.
The closed algebraic saturated set
$$
 E=s\bigl(t^{-1}(C\cup s(N_s)\cup t(N_t))\bigr)_{\mathrm{red}}
$$
is a proper subset of $B$: over a nonempty open of
$Y_{\mathrm{reg}}$, the proper map $\Phi$ is a submersion by Sard's
theorem, and its complete
fibers avoid the displayed bad loci. Set $U=B\setminus E$.
The smooth proper equivalence relation $R_U\rightrightarrows U$
has an algebraic-space quotient $u:U\to Q$
\cite[Tag 04S6]{Stacks}. Composition is algebraic because its
source is smooth over $U$, hence reduced, and the map given by the
two outer coordinates factors through $R_U$ on closed points. The
quotient is smooth, separated and of finite type; $u$ is smooth and
proper.

Local analytic sections of $u$ give an injective immersion
$Q^{\mathrm{an}}\to\Gamma\backslash D$ with image $V=Y\setminus\Phi(E^{\mathrm{an}})$.
The source has dimension $d$ and the differential has rank $d$.
For compact $K\subset V$, its inverse image is
$u^{\mathrm{an}}(\Phi^{-1}K)$, hence compact. Consequently
$Q^{\mathrm{an}}\simeq V$.

Choose a projective compactification $\bar B$. The closed immersion
$(j,u):U\hookrightarrow\bar B\times Q$ is a flat family of complete
fibers, and defines
$\beta:Q\to H_U\subset H=\operatorname{Hilb}_{P_0}(\bar B)$,
where $P_0$ is their Hilbert polynomial and $H_U$ is the open locus
of subschemes supported in $U$.
The map $\beta$ is a monomorphism: if two embedded families agree,
their maps to $Q$ agree after the fppf cover given by that family.
To prove properness, consider a flat DVR family
$C_A\subset U\times\operatorname{Spec}A$ extending a generic fiber.
The two maps $C_A\times_A C_A\rightrightarrows Q$ agree on the
schematically dense generic fiber, hence everywhere. Fppf descent
gives $\operatorname{Spec}A\to Q$. Flatness identifies both families
with the closure of their generic fiber. Thus $\beta$ is a closed
immersion, and $Q$ is quasi-projective.

Let $Q_B=\bar Q\cap H_B$, where $\bar Q$ is the reduced closure in
$H$ and $H_B$ is the open locus of subschemes supported in $B$.
Normalize and resolve $Q_B$, without changing its smooth open $Q$,
to obtain $M$.
The universal family and its fiber square are flat over $Q_B$.
Their inverse images of $Q$ are dense; hence the two evaluation maps
to $\Gamma\backslash D$, equal there, agree everywhere. The support
of each Hilbert fiber lies in one $\Phi$-fiber. The proper projection
of the universal family therefore descends $\Phi$ to a continuous
map $Q_B^{\mathrm{an}}\to Y$.
On its normalization the map is holomorphic over $Q$ and locally
bounded; the Riemann extension theorem makes it holomorphic
everywhere.

For closed $C_0\subset\bar B^{\mathrm{an}}$, the support condition
$S_{C_0}=\{h:\operatorname{Supp}\mathcal F_h\subset C_0\}$ is closed
in $H^{\mathrm{an}}$: its complement is the open image of
$e^{-1}(\bar B^{\mathrm{an}}\setminus C_0)$ under the flat universal
projection. For compact $K\subset Y$, properness of $\Phi$ gives
$$
 h^{-1}(K)=\bar Q^{\mathrm{an}}\cap S_{\Phi^{-1}(K)}
 \subset Q_B^{\mathrm{an}},
$$
a compact set. Thus the induced $p:M^{\mathrm{an}}\to Y$ is proper,
surjective and an isomorphism over $V$. Pullback from
$\Gamma\backslash D$ gives a polarized integral variation on $M$:
Griffiths transversality holds on the dense open $Q$ and hence
everywhere by holomorphicity.

For reducible $Y$, perform this construction on a source component
dominating each maximal $Y_i$. In $M=\coprod M_i$, take the union of
the complements of the open sets where the maps are isomorphisms
and the images of all mixed relations $R_{M_i,M_j}$, $i\ne j$, and
saturate once by
$R_{MM}$. Lemma~\ref{aq:relation} and properness make the resulting
$W$ algebraic and closed. Its intersection with each $M_i$ is a
proper subset, since general points of $Y_i$ avoid the other maximal
components. It is empty on
zero-dimensional components. This gives \eqref{aq:center}.
\end{proof}

\section{Infinitesimal neighborhoods and parallel transport}
\label{aq:finite-section}

\begin{lemma}[Finite horizontal parallel transport]\label{aq:taylor}
Let $S$ be smooth algebraic, $G\subset S$ reduced closed with ideal
$I$, and $(\mathcal E,\nabla)$ an algebraic integrable connection.
If an algebraic $t_0\in H^0(G,\mathcal E|_G)$ is the restriction of
local horizontal sections, their restrictions modulo every ordinary
power $I^n$ are algebraic and compatible.
\end{lemma}

\begin{proof}
On $S=\operatorname{Spec}A$, let $J$ be the diagonal ideal and set
$$
 B_N=(A\otimes A)/(I\otimes A+J^{N+1}),\quad
 Q_N=p_{2*}(p_2^*\mathcal E|_{B_N}),\quad
 \eta_N:\mathcal E\to Q_N.
$$
The connection identifies fibers over infinitesimally neighboring
points by parallel transport;
see \cite[Chapter~I, \S\S2.2--2.3, 2.21]{Deligne} for the
infinitesimal interpretation. In \'etale coordinates, its finite
form is the usual truncated Taylor formula
$$
 \epsilon_N(p_1^*e)=
 \sum_{|\alpha|\le N}\frac{(p_1^*x-p_2^*x)^\alpha}{\alpha!}
             p_2^*(\nabla_\partial^\alpha e),
 \qquad \epsilon_N:p_1^*\mathcal E\xrightarrow{\sim}p_2^*\mathcal E.
$$
Leibniz's rule and integrability give coordinate independence, the
cocycle identity and compatibility in $N$. For the principal-parts map,
$$
 \ker\eta_N=I^{\langle N+1\rangle}\mathcal E,\qquad
 I^{\langle m\rangle}=
 \{b:D(b)\in I\text{ for }\operatorname{ord}D<m\}.
$$
These differential powers are cofinal with ordinary powers.
Indeed, write $I^n=\bigcap Q_j$, with $Q_j$ primary to $\mathfrak p_j$,
and choose $m$ such that $\mathfrak p_j^m\subset Q_j$ for every $j$.
The characteristic-zero Zariski--Nagata theorem
\cite[Theorem~3.6]{ZN} gives
$$
 I^{\langle m\rangle}\subset
 \mathfrak p_j^{\langle m\rangle}
 =\mathfrak p_j^mA_{\mathfrak p_j}\cap A
 \subset Q_jA_{\mathfrak p_j}\cap A=Q_j.
$$
The algebraic section $q_N=\epsilon_N(p_1^*t_0)$ has zero class in
$\operatorname{coker}\eta_N$ after analytification, because it is the
restriction of the given horizontal section under this transport.
Faithful flatness shows that this class is zero. For $N+1\ge m$,
its image in $\mathcal E/I^n\mathcal E$ is the desired section.
Faithfulness gives compatibility, so the local sections glue.
\end{proof}

Fix a smooth quasi-projective period source $X$, an algebraic
map $h:\Delta\to X$ from a separated finite-type scheme, and put
$V=\Delta_{\mathrm{red}}$, $\mathcal N=\ker(\cO_\Delta\to\cO_V)$,
$\mathcal N^a=0$. Suppose a reduced separated finite-type algebraic
space $C$ has a specified closed analytic immersion $i:C^{\mathrm{an}}
\hookrightarrow\Gamma\backslash D$, and a proper algebraic
surjection $f:V\to C$ with
\begin{equation}\label{aq:finite-data}
 \Phi_Xh|_{V^{\mathrm{an}}}=if^{\mathrm{an}}.
\end{equation}
The same argument applies to a locally closed $C$ after restricting
the ambient space to an open neighborhood. Over each target chart,
we take the full inverse image in the source.

\begin{proposition}[A smooth reduced center]\label{aq:smooth-center}
If $C$ is smooth, the scheme-theoretic image of
$\Phi_Xh:\Delta^{\mathrm{an}}\to\Gamma\backslash D$, together with
the induced map from $\Delta$, is algebraic.
The construction preserves the given Hodge filtration and its tensor
operations, and is compatible with restriction to closed subschemes
of $\Delta$ with the same reduction.
\end{proposition}

\begin{proof}
The analytic map is proper. Its scheme-theoretic image $T$ is
a closed analytic subspace \cite{GrauertRemmert}, and
$I_{C/(\Gamma\backslash D)}^{a}\subset I_{T/(\Gamma\backslash D)}$.
Work on an affine \'etale chart $U\to C$.
Its unique lift to $T$ pulls back to the unique algebraic \'etale
lift $\Delta_U\to\Delta$ of $V\times_CU\to V$; retain this entire
proper source and suppress the subscript.

The map $i$ is horizontal: resolve a component of $V$ dominating
each component of $C$ and use the generically surjective differential
of $f$; Griffiths transversality then holds everywhere on the smooth
space $C$. Griffiths' algebraicity theorem gives algebraic Hodge
bundles and filtration, with a regular singular connection.
Deligne's full faithfulness
\cite[II, Theorem~5.9]{Deligne} algebraizes the prescribed comparisons on
overlaps.

Let $G$ be the reduced image of the proper map $(h|_V,f):V\to X\times U$.
Equality of period points gives a specified flat isomorphism
$$
 u_0:p_X^*\mathcal H_X|_{G^{\mathrm{an}}}
 \simeq p_U^*\mathcal H_U|_{G^{\mathrm{an}}}.
$$
On a proper resolution $r:\widetilde G\to G$, this comparison is
algebraic by full faithfulness for connections with regular
singularities. For
the Hom bundle $\mathcal E$,
$\mathcal E\hookrightarrow r_*r^*\mathcal E$ has coherent quotient;
the class of the algebraic comparison vanishes analytically, hence
algebraically. Thus $u_0$, its inverse and its compatibility with
the Hodge filtrations descend algebraically to $G$. In local flat frames $u_0$ belongs to
$\Gamma$. Neatness makes this isomorphism unique at each period point,
so all branches through a point define one ambient horizontal germ
along $G$. Lemma~\ref{aq:taylor} gives its restrictions to the
infinitesimal neighborhoods defined by ordinary powers of the ideal
of $G$.

Let $\mathcal Q_U\to U$ be the relative flag variety with fiber
$\check D$ and Hodge section $\sigma$. Any local algebraic lift
$r:\Delta\to U$ of $f$ has $(h,r)^*I_G\subset\mathcal N$, and the
resulting horizontal isomorphism transports the filtration to an
algebraic map $F_r:\Delta\to\mathcal Q_U$.
On the open subset parametrizing polarized Hodge structures, there
is a holomorphic evaluation map
$$
 \begin{gathered}
 q:(\widetilde U\times D)/\Gamma\longrightarrow\Gamma\backslash D,
 \qquad[\widetilde u,F]\mapsto[F],\\
 \widetilde U=U^{\mathrm{an}}\times_{\Gamma\backslash D}D.
 \end{gathered}
$$
The horizontal isomorphism gives
\begin{equation}\label{aq:evaluation}
 qF_r^{\mathrm{an}}=\Phi_Xh
 \quad\text{on all of }\Delta^{\mathrm{an}}.
\end{equation}

Write $\mathcal T=T_U$, $\mathcal V=\sigma^*T_{\mathcal Q_U/U}$.
The algebraic period differential $\kappa:\mathcal T\hookrightarrow
\mathcal V$ is a subbundle inclusion. After affine refinement choose
$\mathcal V=\kappa\mathcal T\oplus E$ and a projection $\ell$ with
$\ell\kappa=1$. Relative algebraic coordinates $u$ near the section
can be chosen with $u\sigma=0$ and vertical differential $1$:
lift a conormal basis on an affine neighborhood and restrict to the
locus where the Jacobian is invertible. The slice $H=\{\ell u=0\}$
satisfies
$$
 d(q|_H):\mathcal T\oplus E\xrightarrow{\sim}i^*T_{\Gamma\backslash D},
 \qquad(v,e)\longmapsto\kappa(v)+e.
$$
Thus $H_{a-1}=V(I_{\sigma/H}^{a})$ is the unique \'etale lift of $U$
over the thickening of $C$ defined by $I_C^a$.

We show that the induced lift of $\Delta$ is algebraic. At the layer
$L_k=\mathcal N^k/\mathcal N^{k+1}$, a temporary lift $r$ has error
$b_r=(\ell u)(F_r)\in f^*\mathcal T\otimes L_k$.
A change by a derivation $\delta$ has base component $\delta$;
\eqref{aq:evaluation} gives vertical component $-\kappa\delta$.
Since $L_k^2=\mathcal NL_k=0$,
\begin{equation}\label{aq:slice-correction}
 b_{r+\delta}=b_r-f^*(\ell\kappa)\delta,
 \qquad \delta=f^*(\ell\kappa)^{-1}b_r.
\end{equation}
The corrected local lifts are unique and therefore glue at every
order.
We obtain algebraic maps $r_1:\Delta\to U$ and
$j:\Delta\to H_{a-1}$, with $qj^{\mathrm{an}}=\Phi_Xh$ and
\begin{equation}\label{aq:filtered-transport}
 j^*F_{\mathrm{taut}}^\bullet\simeq h^*F_X^\bullet,
 \qquad j^*L_{\mathrm{taut}}\simeq h^*L_X
\end{equation}
for every prescribed determinant tensor $L$ of the Hodge filtration.
These identifications agree with the prescribed analytic comparisons.

The reduction of $r_1$ is proper, so $r_1$ and $j$ are proper.
Since $H_{a-1}\to U$ is finite, the coherent finite algebra
$$
 \im\bigl((H_{a-1}\to U)_*\cO_{H_{a-1}}
                    \longrightarrow r_{1*}\cO_\Delta\bigr)
$$
defines an algebraic image whose analytification is precisely the
chosen \'etale lift $T_U$, by proper comparison. On overlaps the
source pullbacks are the same \'etale lift over $\Delta$. Lemma
\ref{aq:graph-descent} algebraizes their specified analytic target
identifications. Uniqueness gives the cocycle condition, descent
to an algebraic space, and compatibility with restriction to closed
subschemes with the same reduction.
The specified filtered comparisons on overlaps descend as well:
view each as a section of the algebraic Hom bundle. Its pullback
to the common proper schematically dominant source is the
algebraic comparison \eqref{aq:filtered-transport}. Apply Lemma
\ref{aq:graph-descent} to its map into the total Hom bundle;
faithfulness verifies invertibility and compatibility with the
filtration and tensor operations.
\end{proof}

\pagebreak[3]
\section{Algebraization of infinitesimal images}
\label{aq:images-section}

\subsection{Specified submodules and square-zero subrings}

\begin{lemma}[Recovery of a specified submodule]\label{aq:coherent-recovery}
Let $Y$ be a separated finite-type algebraic space, $S=V(I)\subset Y$, and
$\mathcal E$ coherent algebraic. An analytic coherent submodule
$\mathcal H\subset\mathcal E^{\mathrm{an}}$ is algebraic if it is
algebraic as a specified submodule off $S$ and every
$\mathcal H+I^n\mathcal E^{\mathrm{an}}$ is algebraic as a specified
submodule.
\end{lemma}

The lemma also applies when $Y$ is nonreduced.

\begin{proof}
Extend the submodule on $Y\setminus S$ algebraically to
$\mathcal H_0\subset\mathcal E$ and saturate:
$\overline{\mathcal H}=\bigcup_r(\mathcal H_0:I^r)$.
It is coherent and $\mathcal E/\overline{\mathcal H}$ has no
$I$-power torsion. Flat analytification preserves this property,
since each annihilator is the kernel of a map of finite modules.
The image of $\mathcal H$ in this quotient is supported on $S^{\mathrm{an}}$,
hence locally $I$-power torsion, and is therefore zero. Put
$\mathcal F=\overline{\mathcal H}^{\mathrm{an}}/\mathcal H$.
Artin--Rees gives a uniform $c$ on a finite affine cover with
$\overline{\mathcal H}\cap I^{n+c}\mathcal E\subset I^n\overline{\mathcal H}$.
Consequently the following specified submodules are algebraic:
\begin{align*}
 \mathcal G_n
 &=\bigl(\overline{\mathcal H}^{\mathrm{an}}\cap
       (\mathcal H+I^{n+c}\mathcal E^{\mathrm{an}})\bigr)
       +I^n\overline{\mathcal H}^{\mathrm{an}}\\
 &=\mathcal H+I^n\overline{\mathcal H}^{\mathrm{an}},\qquad
 \mathcal G_n/\mathcal G_{n+1}=\operatorname{Gr}_I^n\mathcal F.
\end{align*}
The supports of these algebraic graded quotients decrease, because
$(I/I^2)\otimes\operatorname{Gr}_I^n\mathcal F\twoheadrightarrow
\operatorname{Gr}_I^{n+1}\mathcal F$. Their stable support is empty:
any complex closed point in it would contradict the local $I$-power
torsion of $\mathcal F$. Thus $\operatorname{Gr}_I^N\mathcal F=0$
for some $N$; Nakayama gives $I^N\mathcal F=0$, and
$\mathcal H=\mathcal G_N$. Uniqueness inside $\mathcal E$ follows by
faithful analytification and gives \'etale descent.
\end{proof}

\begin{lemma}[Recovery of a square-zero subring]\label{aq:subring}
Let
$$
 0\longrightarrow J\longrightarrow\mathcal R
 \longrightarrow B=\cO_{Z_-}\longrightarrow0,\qquad J^2=0,
$$
be an algebraic square-zero extension with coherent kernel, where
$Z_-$ is separated of finite type. Let $T\subset\mathcal R^{\mathrm{an}}$
be the structure sheaf of an analytic thickening mapping surjectively to
$B^{\mathrm{an}}$ with coherent kernel. For an algebraic closed
$S\subset Z_-$ put $I=\ker(B\to\cO_S)$ and
$\mathfrak a=\ker(\mathcal R\to\cO_S)$.
If $T$ is algebraic as a specified subring off $S$ and all rings
$$
 A_m=\im\bigl(T\to(\mathcal R/\mathfrak a^m)^{\mathrm{an}}\bigr)
$$
are algebraic as specified subrings, then so is $T$.
\end{lemma}

\begin{proof}
Put $K=T\cap J^{\mathrm{an}}$. Surjectivity onto $B^{\mathrm{an}}$
and $J^2=0$ make $K$ an ideal of $\mathcal R^{\mathrm{an}}$.
Consider the algebraic modules
$$
 L_m=\ker(A_m\to B/I^m),\qquad
 Q_m=J/(J\cap\mathfrak a^m).
$$
The first kernel is the image of $K$: if $t\in T$ has image in
$I^m$, lift each factor in an expression as a sum of $m$-fold
products from $I$ to $T\cap\mathfrak a$ and subtract their products
from $t$. The result belongs to $K$ and has the same truncated
image. Both $L_m$ and $Q_m$ are $B/I^m$-modules, since $J^2=0$.
Hence $J\times_{Q_m}L_m$ algebraizes
$J^{\mathrm{an}}\cap(T+\mathfrak a^m\mathcal R^{\mathrm{an}})$.
The section $T/K\simeq B^{\mathrm{an}}$ gives
$$
 \mathcal R^{\mathrm{an}}/K=B^{\mathrm{an}}\oplus M,
 \quad M=J^{\mathrm{an}}/K,\quad
 (\mathfrak a/K)^m=I^m\oplus I^{m-1}M.
$$
Thus
\begin{equation}\label{aq:kernel-truncations}
 J^{\mathrm{an}}\cap(T+\mathfrak a^m\mathcal R^{\mathrm{an}})
 =K+I^{m-1}J^{\mathrm{an}}.
\end{equation}
Lemma~\ref{aq:coherent-recovery} shows that $K$ is algebraic as a
submodule.

Write $\mathcal R'=\mathcal R/K$, $\mathfrak a'=\mathfrak a/K$.
The specified analytic section $s:B^{\mathrm{an}}\to\mathcal R'^{\mathrm{an}}$
has algebraic truncations
\begin{equation}\label{aq:section-truncations}
 s_m:B/I^m\longrightarrow\mathcal R'/(\mathfrak a')^m.
\end{equation}
Indeed, the algebraic image of $A_m$ in this quotient maps
isomorphically to $B/I^m$, as the displayed splitting shows after
analytification. The image is algebraic by finite proper comparison:
the opposite morphism is the identity on the common reduction.

To recover the section itself, let $Y'=(Z_-,\mathcal R')$, take an
affine \'etale chart $A\to Y'$, and put $U=A\times_{Y'}Z_-$.
The analytic retraction $\rho:Y'^{\mathrm{an}}\to Z_-^{\mathrm{an}}$
corresponding to $s$ lifts uniquely to
$\rho_A:A^{\mathrm{an}}\to U^{\mathrm{an}}$.
For algebra generators $z_i$ of $\Gamma(U,\cO_U)$ put
$f_i=\rho_A^*z_i$. The retraction $\rho_m$ given by
\eqref{aq:section-truncations} identifies, by uniqueness of \'etale lifts,
$$
 A\times_{Y'}V((\mathfrak a')^m)
 \simeq V((\mathfrak a')^m)
       \times_{V(I^m),\rho_m}(U\times_{Z_-}V(I^m)).
$$
Therefore every restriction of $f_i$ to a finite infinitesimal
neighborhood of $S_A=V(\mathfrak a'\cO_A)$ is algebraic; so is its restriction off
$S_A$. The graph $\mathcal G_i$ of $f_i$ is algebraic off $S_A$, and
the submodules $\mathcal G_i+(\mathfrak a')^m
(\cO_A\oplus\cO_A)^{\mathrm{an}}$ are algebraic by the preceding
restrictions. Lemma~\ref{aq:coherent-recovery} shows that
$\mathcal G_i$ is algebraic. Its first projection is analytically,
hence algebraically, an isomorphism; the second gives $f_i$
algebraically.
Relations, the retraction identity and agreement on overlaps follow
by faithful analytification. Thus the prescribed section $s$ is
algebraic.

Finally
$$
 T_{\mathrm{alg}}=\mathcal R\times_{\mathcal R',s}B
$$
is a square-zero extension of $B$ by coherent $K$, realized by an
algebraic space \cite[Tag 05ZT]{Stacks}. It is of finite type and
$\operatorname{Spec}\mathcal R\to\operatorname{Spec}T_{\mathrm{alg}}$
is finite, since it is the identity on the common reduction. On this
target, the kernel sequence
$$
 0\to T_{\mathrm{alg}}\to\mathcal R\oplus B
   \xrightarrow{(r,b)\mapsto\bar r-s(b)}\mathcal R'
$$
and proper comparison recover exactly $T\subset\mathcal R^{\mathrm{an}}$.
\end{proof}

\subsection{Filtrations under proper morphisms and singular centers}

\begin{lemma}[Proper comparison of the two filtrations]\label{aq:rees}
Let $g:\Delta\to C_0$ be proper and schematically dominant, with
$C_{0,\mathrm{red}}=C$. Suppose $\psi:C_0^{\mathrm{an}}\to\Gamma\backslash D$
restricts to a specified immersion of $C^{\mathrm{an}}$, and has
image ring $T\subset\cO_{C_0}^{\mathrm{an}}$.
For $S\subset C$ closed algebraic, put
$$
 \mathfrak a=I_{S/C_0},\quad
 \mathfrak b=I_{(g^{-1}S)_{\mathrm{red}}/\Delta},\quad
 \Delta_l=V(\mathfrak b^l).
$$
If the scheme-theoretic images $U_l$ of
$\Delta_l^{\mathrm{an}}\to\Gamma\backslash D$, with their proper
source maps, are algebraic, then all images of $T$ in
$(\cO_{C_0}/\mathfrak a^m)^{\mathrm{an}}$ are algebraic as specified
subrings.
More generally, for every coherent $\mathcal E$ on $\Delta$ there
exist $c,d\ge1$ such that
\begin{equation}\label{aq:rees-coefficients}
 g_*(\mathfrak b^{c(m+d)}\mathcal E)
 \subset\mathfrak a^m g_*\mathcal E\qquad(m\ge0).
\end{equation}
\end{lemma}

\begin{proof}
Let $B_l=\im_{\mathrm{sch}}(\Delta_l\to C_0)$, with ideal
$\mathfrak b_l$. Ambient functions vanishing on $U_l$ vanish on
$\Delta_l$ and hence on $B_l$. Lemma~\ref{aq:graph-descent} therefore
algebraizes the specified $B_l\to U_l$.
On an affine \'etale chart write $A=\cO_{C_0}$ and
$M=g_*\cO_\Delta$, so $A\hookrightarrow M$ is finite.
The graded module $\bigoplus_qg_*(\mathfrak a^q\cO_\Delta)t^q$
is finite over $\mathscr R=\bigoplus_q\mathfrak a^qt^q$:
the relative spectrum of $\bigoplus_q\mathfrak a^q\cO_\Delta t^q$
is closed in $\Delta\times_{C_0}\operatorname{Spec}\mathscr R$,
and proper coherence applies. Consequently, for some $d,e,c$,
$$
 g_*(\mathfrak a^q\cO_\Delta)\subset\mathfrak a^{q-d}M,
 \quad A\cap\mathfrak a^{m+e}M\subset\mathfrak a^mA,
 \quad\mathfrak b^c\subset\mathfrak a\cO_\Delta.
$$
The middle inclusion is Artin--Rees; the last follows from equality
of radicals. They give
\begin{equation}\label{aq:rees-bound}
 l=c(m+d+e)\quad\Longrightarrow\quad
 \mathfrak b_l=A\cap g_*\mathfrak b^l\subset\mathfrak a^m.
\end{equation}
Take the scheme-theoretic image of the finite morphism
$V(\mathfrak a^m)\hookrightarrow B_l\to U_l$; its reduction is the
identity of $S_{\mathrm{red}}$. Its analytic ring is the image of
$T$ in $(A/\mathfrak a^m)^{\mathrm{an}}$: both consist of pullbacks
of the same ambient functions. Thus the truncated source recovers
the required truncated image. With coherent coefficients,
the Rees module
$\bigoplus_q\mathfrak a^q\mathcal E\,t^q$ is finite over the same
source Rees algebra. Its proper direct image is finite over
$\mathscr R$, giving
$g_*(\mathfrak a^q\mathcal E)\subset\mathfrak a^{q-d}g_*\mathcal E$.
Together with $\mathfrak b^c\subset\mathfrak a\cO_\Delta$, this is
\eqref{aq:rees-coefficients}.
\end{proof}

\begin{theorem}[Infinitesimal images]\label{aq:finite-images}
Under \eqref{aq:finite-data}, the scheme-theoretic image
$$
 T_\Delta=\im_{\mathrm{sch}}(\Delta^{\mathrm{an}}
                                \xrightarrow{\Phi_Xh}\Gamma\backslash D)
$$
is algebraic, with its specified proper schematically dominant
source map. The algebraic model is unique with this identification.
Closed immersions of compatible algebraic sources with the same
reduction induce algebraic closed immersions of these images,
compatible with composition.
\end{theorem}

\begin{proof}
Induct on $\dim C$, and within a fixed dimension on the nilpotence
index of $\mathcal N$. A zero-dimensional center is smooth; a
reduced source has the given image $C$.
Suppose $\mathcal N^r=0$ and set
$\Delta_-=V(\mathcal N^{r-1})$, $\mathcal L=\mathcal N^{r-1}$.
The inner induction gives the algebraic image $Z_-$ and a
proper schematically dominant $f_-:\Delta_-\to Z_-$.

We first construct a square-zero envelope. On the \'etale site of
$Z_-$, take the sheaf of rings $f_{-*}\cO_\Delta$: a chart is lifted
uniquely from $\Delta_-$ to $\Delta$ before taking sections. No map
$\Delta\to Z_-$ or module structure over $\cO_{Z_-}$ is assumed yet.
The connecting homomorphism of the additive extension is
$$
 \partial:f_{-*}\cO_{\Delta_-}\longrightarrow R^1f_{-*}\mathcal L.
$$
Products of local lifts make $\partial$ a derivation, and it kills
$\cO_{Z_-}$. To see this, choose affine \'etale $U\to Z_-$ and the
entire proper source $D_U=\Delta_-\times_{Z_-}U$ with its lift
$D_U^+$. Proper comparison, affine and Stein acyclicity give the
injective map
$$
 H^1(D_U,\mathcal L_U)
 =\Gamma(U,R^1f_{U*}\mathcal L_U)
 \hookrightarrow
 \Gamma(U^{\mathrm{an}},(R^1f_{U*}\mathcal L_U)^{\mathrm{an}}).
$$
The scheme-theoretic image of $\Delta^{\mathrm{an}}$ is a square-zero extension
of $Z_-^{\mathrm{an}}$: its kernel pulls into $\mathcal L^{\mathrm{an}}$
and hence squares to zero by schematic dominance.
On its chosen \'etale lift over $U$, each $z\in\Gamma(U,\cO_U)$ lifts
locally on the target, hence over the whole inverse image of those
neighborhoods. Its analytic connecting section is zero. The
displayed injection shows that its algebraic obstruction is zero.
Naturality is seen on a finite affine cover of $D_U^+$: the same
differences of local lifts give the algebraic and analytic
$\mathcal L_U$-valued \v Cech cocycles.

It follows that the sheaf of rings
\begin{equation}\label{aq:envelope}
 \mathcal R=\cO_{Z_-}
       \times_{f_{-*}\cO_{\Delta_-}}f_{-*}\cO_\Delta
 \longrightarrow\cO_{Z_-}
\end{equation}
is surjective with square-zero coherent kernel $J=f_{-*}\mathcal L$.
Its algebraic-space realization \cite[Tag 05ZT]{Stacks} is
$C_0=(Z_-,\mathcal R)$, with a proper schematically dominant
$g:\Delta\to C_0$. This fiber-product construction also appears in
\cite[Section~4.3, especially Lemma~4.10]{BBT}.
Now the terms in \eqref{aq:envelope} are coherent modules over $C_0$,
and proper comparison identifies its analytification with the same
analytic fiber product. Ambient period functions define a specified
$\psi:C_0^{\mathrm{an}}\to\Gamma\backslash D$ whose image ring
$T\subset\mathcal R^{\mathrm{an}}$ surjects onto $\cO_{Z_-}^{\mathrm{an}}$.

Over $C_{\mathrm{reg}}$, Proposition~\ref{aq:smooth-center} and
Lemma~\ref{aq:graph-descent} make this specified subring algebraic.
For $S=(C\setminus C_{\mathrm{reg}})_{\mathrm{red}}$, all sources
$\Delta_l=V(I_{(g^{-1}S)_{\mathrm{red}}/\Delta}^{l})$ still map
algebraically to the smooth period source $X$. Their reduced maps
to $S$ are algebraic proper surjections, and $\dim S<\dim C$.
The outer induction makes their scheme-theoretic images algebraic.
Lemma
\ref{aq:rees} supplies all finite restrictions required by Lemma
\ref{aq:subring}; that lemma recovers $T$ itself inside
$\mathcal R^{\mathrm{an}}$. The composite $\Delta\to C_0\to T_{\mathrm{alg}}$
is the required proper schematically dominant source map.

For descent, lift an affine \'etale chart of $C$ uniquely through
the two nilpotent targets. Their source
pullbacks are canonically the same lift over $\Delta$.
Lemma~\ref{aq:graph-descent} applies to that entire proper scheme
source. It gives descent, uniqueness, and the specified transition
maps. Faithful analytification detects the closed-immersion
property and all composition identities.
\end{proof}

\section{Formal algebraization and the original source}
\label{aq:formal-section}

\subsection{The formal target and its convergent local rings}

\begin{lemma}[Control of proper analytic kernels]\label{aq:analytic-kernels}
Let $p:X\to V$ be proper analytic, $\mathcal E$ coherent, and
$p^{-1}(y)\subset V(I)$. Set
$$
 D=(p_*\mathcal E)_y,\qquad
 H_n=\ker[D\to(p_*(\mathcal E/I^{n+1}\mathcal E))_y].
$$
For every $k$ some $n$ satisfies $H_n\subset\mathfrak m_y^kD$.
\end{lemma}

\begin{proof}
Near the compact fiber, $I\subset\sqrt{\mathfrak m_y\cO_X}$.
A finite cover and properness give $I^d\subset\mathfrak m_y\cO_X$
after shrinking the base. Put
$D_j=(p_*(\mathcal E/\mathfrak m_y^j\mathcal E))_y$ and
$E_j=\ker(D\to D_j)$. The analytic theorem on formal functions
\cite[Chapter VI, Corollary~4.7(i)]{BS} gives
$\widehat D\simeq\varprojlim_jD_j$. The $D_j$ have finite length;
their kernels in $\widehat D$ decrease to zero. Chevalley's lemma
and faithful flat contraction give $E_j\subset\mathfrak m_y^kD$
for large $j$. If $n+1\ge dj$, then $H_n\subset E_j$.
\end{proof}

\begin{proposition}[Algebraization from infinitesimal images]\label{aq:formal-recovery}
Let $p:M^{\mathrm{an}}\to Y$ be a proper surjection, with $M$ smooth
quasi-projective, allowing a finite disjoint union, and $Y$ reduced.
Suppose $p$ is an isomorphism off a reduced center $Z$, the whole
$W=(p^{-1}Z)_{\mathrm{red}}\subset M$ is algebraic, and $Z$ has a
specified separated finite-type algebraic model for which $W\to Z$
is proper algebraic. Write
$$
 I=I_{W/M},\qquad W_n=V(I^{n+1}),\qquad
 T_n=\im_{\mathrm{sch}}(W_n^{\mathrm{an}}\to Y).
$$
Assume all $T_n$ have separated finite-type algebraic models and
the maps $W_n\to T_n$ and inclusions
$Z\hookrightarrow T_n\hookrightarrow T_{n+1}$ are algebraic with
the prescribed compatible analytic identifications. Then there
exist a reduced separated finite-type algebraic space $Y_{\mathrm{alg}}$, a proper
$q:M\to Y_{\mathrm{alg}}$, and an analytic isomorphism
$$
 \theta:Y_{\mathrm{alg}}^{\mathrm{an}}\xrightarrow{\sim}Y,
 \qquad\theta q^{\mathrm{an}}=p.
$$
\end{proposition}

\begin{proof}
If $Z$ is empty, take $Y_{\mathrm{alg}}=M$. Otherwise, on the
Noetherian formal completion $\widehat M_W$ set
$$
 P_n=\ker(\cO_{T_n}\to\cO_Z)\cO_{W_n},\qquad P=\varprojlim P_n.
$$
Surjectivity of the transitions gives, on each formal affine chart,
$P_n=(P+I^{n+1})/I^{n+1}$. The finitely generated ideal $P$ is
coherent and glues by completed localization. If
$Q=I_{Z/Y}\cO_{M^{\mathrm{an}}}$, then
$$
 P_n^{\mathrm{an}}=(Q+(I^{\mathrm{an}})^{n+1})/(I^{\mathrm{an}})^{n+1},
 \qquad\sqrt Q=I^{\mathrm{an}}.
$$
There is a uniform $c$ such that
\begin{equation}\label{aq:source-cofinality}
 I^c\subset P\subset I.
\end{equation}
Indeed, at a complex point $x$, choose $c_x$ with
$(I_x^{\mathrm{an}})^{c_x}\subset Q_x$. On the finite algebraic
level $W_{c_x}$, faithful analytification gives a Zariski neighborhood
with $I^{c_x}\subset P+I^{c_x+1}$. Nakayama on
$(I^{c_x}+P)/P$ gives $I^{c_x}\subset P$. These neighborhoods cover
the finite-type Jacobson space $W$; a finite subcover gives $c$.

\emph{Noetherianity.}
Lift an affine \'etale chart $U=\operatorname{Spec}A_0\to Z$
uniquely to affine $U_n\to T_n$, and set
$$
 A=\varprojlim A_n,\qquad K_n=\ker(A\to A_n),\qquad
 J_n=\ker(A_n\to A_0).
$$
The schemes $F_n=W_n\times_{T_n}U_n$ are separated quasi-finite over
$W_n$; they form the \'etale lift $\mathfrak F$ of the formal
source. We first prove that $A$ is Noetherian.
Choose $N\ge c$ and lift generators of $J_N$ to $b_j\in A$.
Their ideal $Q'$ on $\mathfrak F$ satisfies
$$
 P=Q'+I^{N+1}\subset Q'+IP,
$$
so Nakayama gives $P=Q'$. Present
$A_0=\C[x_1,\ldots,x_s]/(r_1,\ldots,r_e)$ and lift $x_i$ to $a_i\in A$.
The Noetherian complete ring
$$
 B=\C[X_1,\ldots,X_s,T_1,\ldots,T_r]^\wedge_{\mathfrak a},
 \quad\mathfrak a=(r_i(X),T_j),
$$
maps to $A$ by $X_i\mapsto a_i,T_j\mapsto b_j$, since each $J_n$
is nilpotent. Moreover $\mathfrak a\cO_{\mathfrak F}=P$.
The $P^k$-quotients of $\mathfrak F$ form a Cartesian system of
proper schemes over $B/\mathfrak a^k$. They are of finite type because
\eqref{aq:source-cofinality} embeds each into a suitable $F_n$.
Their first layer $\mathfrak F/P$ is a nilpotent thickening of
$F_0$; properness follows from the proper reduction $F_0\to U$.
An ample line bundle on $M$ pulls back to an ample line bundle on
$F_0$: $F_0\to W$ is separated and quasi-finite, hence quasi-affine.
Its compatible extensions remain
ample on nilpotent thickenings. Grothendieck algebraization
\cite[Tag 089A]{Stacks} gives proper $F_B\to\operatorname{Spec}B$,
and the theorem on formal functions gives
$$
 H=\Gamma(F_B,\cO)=\Gamma(\mathfrak F,\cO),
 \qquad H\text{ finite over }B.
$$
The injections $A_n\hookrightarrow\Gamma(F_n,\cO)$ embed $A$ as
a $B$-submodule of $H$. Thus $A$ is finite over $B$ and Noetherian.

Nilpotence gives $\mathfrak a^mA\subset K_n$ for large $m$.
Conversely, the proper Rees argument of Lemma~\ref{aq:rees} gives,
for some fixed $d$ and all $k\ge d$,
$$
 \mathfrak a^kH\subset
 \ker\bigl[H\to H^0(F_B,\cO_{F_B}/\mathfrak a^k\cO_{F_B})\bigr]
 \subset\mathfrak a^{k-d}H.
$$
Thus the $P$-adic kernel topology on $H$ is its
$\mathfrak a$-adic topology. By
\eqref{aq:source-cofinality} and Artin--Rees for $A\subset H$,
the original $K_n$-topology on $A$ is the same adic topology.
For $J=K_0$, choose $K_n\subset\mathfrak a A$; nilpotence of
$J_n$ gives
$$
 J^e\subset K_n\subset\mathfrak a A\subset J.
$$
Thus $A$ is $J$-adic. The unique \'etale lifts and their cocycles
glue these charts into a locally Noetherian adic formal algebraic
space $\mathfrak Y$, with reduction $Z$, and
$\widehat p:\widehat M_W\to\mathfrak Y$ satisfying
$J\cO_{\widehat M_W}=P$. The map is proper: its adic Cartesian
truncations are of finite type with proper reduction. The target is
separated, since locally a scheme testing its diagonal factors
through a common separated $T_n$.

\emph{The completed map.}
At $y\in Z$, put
$S_y=\cO_{Y,y}\hookrightarrow D_y=(p_*\cO_{M^{\mathrm{an}}})_y$.
This is an inclusion in a finite $S_y$-module. For
$L_n=\ker(S_y\to(p_*\cO_{W_n^{\mathrm{an}}})_y)$, Lemma
\ref{aq:analytic-kernels} and Artin--Rees give
$L_n\subset\mathfrak m_y^\ell$ for all sufficiently large $n$,
for each $\ell$. Choose $n$ with
$K_n\subset J^\ell$, $L_n\subset\mathfrak m_y^\ell$, and $n+1\ge\ell$.
The specified $W_n\to T_n$ identify the $\ell$-jets of source and
target and their map. At $x$ over $y$, choose an affine formal
source chart $\operatorname{Spf}C$ and its maximal ideal $\mathfrak n$;
$\mathfrak m$ denotes the chosen maximal ideal of the target chart
$A$. Passing to limits gives
\begin{equation}\label{aq:completed-map}
 \begin{array}{ccc}
 \widehat{A_{\mathfrak m}}&\longrightarrow&\widehat{C_{\mathfrak n}}\\
 \downarrow\wr&&\downarrow\wr\\
 \widehat{\cO_{Y,y}}&\longrightarrow&\widehat{\cO_{M^{\mathrm{an}},x}}.
 \end{array}
\end{equation}
The vertical isomorphisms preserve the center ideals and commute
with the indicated homomorphisms.

We use the following uniform annihilation bound: if $C$ is Noetherian
and complete for $L$, $C/L$ is finite type over $\C$, and a finite
module is killed after completion at every closed point of the
center by some power of $L$, then one power of $L$ kills it. Faithful flatness
gives the local assertion; finite generation spreads it to a
Zariski neighborhood. Since $L\subset\operatorname{Jac}(C)$, every
maximal ideal is a point of the center, and a finite cover gives the
uniform exponent.

\emph{The formal modification.}
For a topologically finitely presented algebra $C=A\{u_1,\ldots,u_r\}/F$,
the completed naive cotangent complex is
$$
 [F/F^2\longrightarrow C^r],\qquad\text{in degrees }-1,0.
$$
Its completion at each maximal ideal agrees in both degrees with
that of the analytic relative cotangent complex of $p$. Indeed,
\eqref{aq:completed-map} identifies two presentations
$\widehat C=\widehat A[[u]]/F_u=\widehat A[[v]]/F_v$ of the same map.
The common presentation has ideal $(F_u,v-\varphi(u))$; its
two-term complex of conormal and continuous differential modules
splits off the contractible summand
$[\widehat C^s\xrightarrow{1}\widehat C^s]$.
Interchange $u,v$ to obtain equality of both $H^{-1}$ and $H^0$.
Flatness of completion gives the same identification for these
finite cohomology modules.
The analytic modules are supported over the center, since $p$ is
an isomorphism elsewhere. The uniform annihilation bound shows
that both formal modules are killed by a power of $P$, and the formal map is
rig-\'etale \cite[Tag 0AJU]{Stacks}.

Take the whole $\mathfrak R=\widehat M_W\times_{\mathfrak Y}\widehat M_W$,
with diagonal ideal $D_\Delta$. At a closed center point $(x_1,x_2)$
its completed ring is, by \eqref{aq:completed-map},
$$
 \widehat{\cO_{M^{\mathrm{an}},x_1}}
 \widehat\otimes_{\widehat{\cO_{Y,y}}}
 \widehat{\cO_{M^{\mathrm{an}},x_2}}.
$$
For $x_1=x_2$, the diagonal ideal is the kernel of multiplication.
For $x_1\ne x_2$, it is the unit ideal and the local fiber product
is supported over the center. In either case, the analytic fiber
product equals the diagonal off the center, so a power of the center
ideal kills the diagonal ideal. Faithful flatness and the uniform
annihilation bound yield
$J_{\mathfrak R}^bD_\Delta=0$ on each chart. The full diagonal is
therefore rig-\'etale and rig-surjective \cite[Tag 0GD6]{Stacks}.

Finally algebraize $\mathfrak F$ over the now Noetherian $A$, using
the same ample line bundle. The proper $F_A\to\operatorname{Spec}A$
satisfies $A\hookrightarrow\Gamma(F_A,\cO)$, hence is surjective.
For an adic complete DVR test $\operatorname{Spf}R\to\operatorname{Spf}A$,
its generic fiber has a point over a finite extension of the
fraction field. Extend and complete the valuation; properness then
gives a lift over the resulting complete DVR. This proves
rig-surjectivity \cite[Tag 0AQP]{Stacks}. All formal-modification
conditions of \cite[Tag 0GDK]{Stacks} are now satisfied.

Artin's theorem in the form that keeps the source fixed
\cite[Theorem~3.1]{Artin1970},
as stated in \cite[Tag 0GIB]{Stacks},
gives a proper $q:M\to Y_{\mathrm{alg}}$ and
\begin{equation}\label{aq:fixed-source}
 \begin{gathered}
 \widehat{Y_{\mathrm{alg}}}_{Z_{\mathrm{alg}}}\simeq\mathfrak Y,
 \qquad\widehat q=\widehat p,\qquad
 q^{-1}(|Z_{\mathrm{alg}}|)=|W|,\\
 M\setminus W\simeq Y_{\mathrm{alg}}\setminus Z_{\mathrm{alg}}.
 \end{gathered}
\end{equation}
The target is separated and locally finite type. The map is
surjective on the center and its complement; quasi-compactness of
$M$ makes $Y_{\mathrm{alg}}$ finite type. The open ideals defining
$T_n$ come from finite algebraic center thickenings, and the fixed
source identifies $q|_{W_n}$ with the specified $W_n\to T_n$.
The kernel of $\cO_{Y_{\mathrm{alg}}}\to q_*\cO_M$ vanishes modulo
every such ideal and off the center. Krull intersection makes it
zero. Hence $Y_{\mathrm{alg}}$ is reduced.

\emph{The analytic local rings.}
The two proper maps $p$ and $q^{\mathrm{an}}$ have the same point
fibers, on the center and its complement. Their quotient topologies
give a homeomorphism $Y_{\mathrm{alg}}^{\mathrm{an}}\simeq Y$.
On this common space put
$$
 \mathcal E=p_*\cO_{M^{\mathrm{an}}}=q^{\mathrm{an}}_*\cO_{M^{\mathrm{an}}},
 \quad\mathcal S=\cO_Y,\quad\mathcal R=\cO_{Y_{\mathrm{alg}}^{\mathrm{an}}},
 \quad\mathcal H_n=\ker(\mathcal E\to p_*\cO_{W_n^{\mathrm{an}}}).
$$
The sheaves $\mathcal S,\mathcal R$ inject into $\mathcal E$, which
is finite over each separately. The specified finite source maps
identify their images in $\mathcal E/\mathcal H_n$ for every $n$.
For $s\in\mathcal S_y$, this approximates $s$ by elements of
$\mathcal R_y$ modulo every $\mathcal H_{n,y}$. Applying Lemma
\ref{aq:analytic-kernels} to $q^{\mathrm{an}}$ gives
$$
 [s]\in\bigcap_{k\ge1}\mathfrak m_{\mathcal R_y}^k
               (\mathcal E_y/\mathcal R_y)=0.
$$
Thus $\mathcal S_y\subset\mathcal R_y$. Applying the same argument
to $p$, with the finite $\mathcal S_y$-module
$\mathcal E_y/\mathcal S_y$, gives the reverse inclusion.
Thus $\mathcal S=\mathcal R$ inside the same direct image. This
identifies the convergent local rings, including the nonnormal ones,
and gives $\theta$.
\end{proof}

\subsection{Induction and the original source}

\begin{proof}[Proof of Theorem~\ref{aq:algebraization}]
First suppose $B$ is a finite disjoint union of smooth
quasi-projective sources and $\Phi$ is proper. Induct on the
dimension $d$ of its image; for $d=0$ group the finitely
many source components by their period value.
Proposition~\ref{aq:modification} supplies $p:M^{\mathrm{an}}\to Y$
and $W,Z$. If $W$ is empty, take $Y_{\mathrm{alg}}=M$, $q=\mathrm{id}_M$
and $\theta=p$, and pass directly to the source descent below.
Otherwise resolve all components of $W$ together:
$r:\widetilde W\to W$ is a proper surjection from a finite disjoint
union of smooth quasi-projective varieties with period maps whose
images cover $Z$. Since $\dim Z<d$, induction algebraizes this whole
center and $\widetilde W\to Z$. Reducedness of $W$ makes $r$
schematically dominant, so Lemma~\ref{aq:graph-descent} algebraizes
the specified proper map $W\to Z$.

Apply Theorem~\ref{aq:finite-images} to the algebraic maps
$W_n=V(I_{W/M}^{n+1})\to M$. Their reduced map $W\to Z$ is now
proper algebraic. It gives every $T_n$ and all specified
maps required by Proposition~\ref{aq:formal-recovery}. That
proposition recovers $Y_{\mathrm{alg}}$, $q:M\to Y_{\mathrm{alg}}$
and $\theta:Y_{\mathrm{alg}}^{\mathrm{an}}\simeq Y$.
The induction in Theorem~\ref{aq:finite-images} involves only
centers of dimension less than $d$.

To recover the original source map, use the full algebraic relation
$$
 E=R_{B,M},\quad a=\operatorname{pr}_B:E\to B,\quad
 b=q\operatorname{pr}_M:E\to Y_{\mathrm{alg}}.
$$
The map $a$ is proper, surjective and schematically dominant, and
$(\theta^{-1}\Phi)a^{\mathrm{an}}=b^{\mathrm{an}}$.
Lemma~\ref{aq:graph-descent} yields the original algebraic $f$;
its analytification is proper and surjective, hence so is $f$.
The construction includes components whose images lie inside other
maximal image components.

For a smooth quasi-projective source with nonproper period map,
take a smooth projective compactification with simple normal crossings
boundary
$\bar B\setminus B=\bigcup D_j$ and set
$$
 B^+=\bar B\setminus\bigcup_{T_j\text{ of infinite order}}D_j.
$$
Griffiths' extension theorem
\cite[Theorem~9.5, Propositions~9.10--9.11 and \S11(a)]{Griffiths1970III}
gives a polarized integral variation on $B^+$ and a proper period
map to $\Gamma\backslash D$. Here finite-order monodromy is the
identity by neatness. One can check properness directly in this fixed quotient:
if periods converged along a sequence approaching an omitted
divisor, the horizontal Schwarz estimate would send its angular
loops to loops of length $O(1/\log(1/|t_j|))$. For large indices
they lie in an evenly covered neighborhood of the limiting period
point, so lift closed to $D$, forcing $T_j=1$, a contradiction.
The preceding proper case applies to $B^+$; restrict its algebraic
map to $B$. Density and properness identify the closed target with
$\overline{\Phi(B^{\mathrm{an}})}_{\mathrm{red}}$.

Finally let $B$ be any reduced separated finite-type algebraic
space. Chow's lemma \cite[Tag 088U]{Stacks} and projective resolution
give a proper surjection
$r:\widetilde B=\coprod_i\widetilde B_i\to B$ from smooth
quasi-projective schemes, isomorphic over a dense regular open in
each component. The pulled-back filtration satisfies Griffiths
transversality on that dense open, hence everywhere. Apply the smooth case to obtain
$g:\widetilde B\to Y_{\mathrm{alg}}$ and its specified closed
analytic target. Reducedness makes $\Phi$ factor through that
reduced target and makes $r$ schematically dominant. Lemma
\ref{aq:graph-descent} algebraizes the given analytic factor from
$B$, whose pullback by $r^{\mathrm{an}}$ is $g^{\mathrm{an}}$.
When $\Phi$ is proper, $g$ is proper and after every base change
$f_T(C)=g_T(r_T^{-1}C)$ for closed $C\subset B_T$. Thus $f$ is
universally closed, and hence proper. The smooth modification of
the proper extended source supplies the asserted $q:M\to Y$.

In the general case $f$ is dominant, and Chevalley's theorem makes
$f(B)$ a dense constructible subset of $Y$.
\end{proof}
\section{The Griffiths line bundle and quasi-projectivity}
\label{aq:griffiths-section}

Fix a polarized integral lattice, a pure weight, and Hodge numbers.
Let $D$ be a connected component of the corresponding period domain,
and let $\Gamma$ be a neat arithmetic subgroup preserving $D$.
Write $\mathcal H,F^\bullet,Q$ for the tautological filtered
bundle and polarization on $\Gamma\backslash D$. If $a$ and
$b$ are the smallest and largest indices with nonzero Hodge numbers, set
\begin{equation}\label{aq:griffiths-line}
 \Lambda=\bigotimes_{p=a+1}^{b}\det F^p.
\end{equation}
The indices $a$ and $b$ remain fixed. The neat hypothesis holds
until Theorem~\ref{aq:finite-level},
which treats arithmetic quotients with stabilizers.

\begin{theorem}\label{aq:quasiprojectivity}
Let $B$ be a reduced separated algebraic space of finite type over
$\C$, and let
$$
 \Phi:B^{\mathrm{an}}\longrightarrow \Gamma\backslash D
$$
be an everywhere defined period map, horizontal on the regular
locus. The reduced analytic closure
$$
 Y_{\mathrm{an}}
   =\overline{\Phi(B^{\mathrm{an}})}^{\,\Gamma\backslash D}_{\mathrm{red}}
$$
has a quasi-projective scheme structure $Y$ and a factorization
\begin{equation}\label{aq:projective-factor}
 B\xrightarrow{\,f\,}Y,\qquad
 i:Y^{\mathrm{an}}\hookrightarrow \Gamma\backslash D,\qquad
 \Phi=i f^{\mathrm{an}},
\end{equation}
where $f$ is algebraic and dominant, and $i$ is a closed analytic
immersion. The Griffiths line bundle $i^*\Lambda$ is the
analytification of an ample algebraic line bundle $L$.

Consequently there are $m>0$, a finite-dimensional
$V\subset H^0(Y,L^m)$, and a locally closed immersion
$$
 j:Y\hookrightarrow\mathbf P(V^\vee),\qquad
 j^*\cO(1)\simeq L^m.
$$
The scheme-theoretic closure $\bar Y\subset\mathbf P(V^\vee)$
has an ample line bundle $\cO_{\bar Y}(1)$ extending $L^m$.
If $\Phi$ is proper, then $f$ is proper and surjective and
$Y^{\mathrm{an}}=\Phi(B^{\mathrm{an}})_{\mathrm{red}}$.
Without properness, $f(B)$ is a dense constructible subset of $Y$.
\end{theorem}

The algebraization theorem in Section~\ref{aq:section} gives
$Y$, $f$, and a proper surjective algebraic map
\begin{equation}\label{aq:smooth-source}
 q:M\longrightarrow Y
\end{equation}
from a finite disjoint union of smooth quasi-projective varieties
carrying the pulled-back variation. We first algebraize the Hodge
bundles and filtration, compatibly with pullback to smooth bases.
Functoriality of the canonical extensions then allows us to apply
the general ampleness criterion, Theorem~\ref{aq:general-ampleness}.

\subsection{Hodge bundles on infinitesimal thickenings}

All algebraic spaces below are separated and of finite type over $\C$.

\begin{proposition}\label{aq:hodge-descent}
Let $T$ have a locally closed analytic
immersion $i_T:T^{\mathrm{an}}\hookrightarrow \Gamma\backslash D$. Suppose that
$$
 g:\Delta\longrightarrow T
 \quad\text{is proper and schematically dominant},\qquad
 h:\Delta\longrightarrow X,
 \qquad i_T g^{\mathrm{an}}=\Phi_Xh^{\mathrm{an}},
$$
where $\Delta$ is a scheme, $X$ is smooth quasi-projective, and
$\Phi_X$ is the period map of a polarized integral variation.
Let $\mathcal V$ be one of the tautological bundles
$\mathcal H,F^p,\operatorname{Gr}_F^p\mathcal H$ on
$\Gamma\backslash D$, or a fixed tensor construction from them.
Let $\mathcal V_X$ be its algebraic realization on $X$. Then
$N_T=i_T^*\mathcal V$ has an algebraic model $\mathcal V_T$
with the specified comparisons
\begin{equation}\label{aq:hodge-descent-comparison}
 \mathcal V_T^{\mathrm{an}}\simeq N_T,\qquad
 g^*\mathcal V_T\xrightarrow{\;\sim\;}h^*\mathcal V_X.
\end{equation}
The model is unique as a submodule of $g_*h^*\mathcal V_X$
with this analytification.
\end{proposition}

The space $T$ need not be reduced. The proposition applies in
particular to the Griffiths line bundle and every
$\bigotimes_p(\det F^p)^{n_p}$, $n_p\in\Z$.

\begin{proof}
Griffiths' algebraicity theorem
\cite[Theorem~(4.13)(a)]{Schmid1973} gives $\mathcal V_X$ with its natural
analytic identification. Put
$$
 M_\Delta=h^*\mathcal V_X,\qquad E=g_*M_\Delta.
$$
Proper comparison for coherent sheaves and the analytic projection
formula give the specified inclusion
\begin{equation}\label{aq:line-in-direct-image}
 N_T\hookrightarrow E^{\mathrm{an}},
 \qquad
 E^{\mathrm{an}}
 \simeq N_T\otimes g_*^{\mathrm{an}}\cO_{\Delta^{\mathrm{an}}}.
\end{equation}
The map is injective because $N_T$ is locally free and $g$ is
schematically dominant. We will algebraize this submodule.

First suppose that $C=T_{\mathrm{red}}$ is smooth. Work over an
affine \'{e}tale chart $U\to C$, with its unique lift $T_U\to T$
and the pullback $\Delta_U$ of $\Delta$.
The smooth-center construction of
Proposition~\ref{aq:smooth-center} gives a relative flag variety
$\mathcal Q_U$, a transverse slice $H$, and a map
$$
 j:\Delta_U\longrightarrow H_{\nu-1},
$$
where $\nu$ is a nilpotence bound for the nilradical of
$\Delta_U$, and $H_{\nu-1}$ is the corresponding finite
thickening of the central section.
The algebraic comparison in \eqref{aq:filtered-transport} preserves
the Hodge filtration:
$$
 \theta\bigl(h^*F_X^\bullet\bigr)
       =j^*F_{\mathcal Q_U}^\bullet.
$$
Applying the prescribed bundle construction gives
\begin{equation}\label{aq:line-taylor-transport}
 j^*\mathcal V_{\mathcal Q_U}
       \xrightarrow{\;\sim\;}M_{\Delta_U}.
\end{equation}
This comparison is the restriction of horizontal transport to the
infinitesimal neighborhood. The analytic evaluation map $e$ on
the locus satisfying the Hodge--Riemann bilinear relations satisfies
$$
 e^*\mathcal V\simeq\mathcal V_{\mathcal Q_U}.
$$
Thus \eqref{aq:line-taylor-transport} is the prescribed comparison.

Let $A=\operatorname{Im}_{\mathrm{sch}}(j)$. The construction
identifies $A^{\mathrm{an}}$ with $T_U^{\mathrm{an}}$, compatibly
with their maps from $\Delta_U$. Both are algebraic finite
thickenings of $U$. Lemma~\ref{aq:graph-descent}, applied to this
common proper source, makes the specified isomorphism algebraic.
Transporting $\mathcal V_{\mathcal Q_U}|_A$ gives a bundle on $T_U$
and the comparison \eqref{aq:hodge-descent-comparison}.
On overlaps, their adjunction maps analytify to the same inclusion
\eqref{aq:line-in-direct-image}. Faithful flatness of
analytification therefore descends the bundles and comparisons to $T$.

For general $C$, induct on $\dim C$, including all finite
nilpotent thickenings at each dimension. The zero-dimensional case
has smooth reduction.
For the induction step put
$$
 S=(C\setminus C_{\mathrm{reg}})_{\mathrm{red}},
 \quad I=\mathcal I_{S/T},\quad
 J=\sqrt{I\cO_\Delta},\quad
 \Delta_\ell=V(J^\ell),\quad
 T_\ell=\operatorname{Im}_{\mathrm{sch}}(\Delta_\ell\to T).
$$
The bundle is already algebraic on $T\setminus S$.
Each $T_\ell$ is an algebraic closed subspace of $T$, with
reduction $S$. Write
$$
 i_\ell:T_\ell\hookrightarrow T,\qquad
 g_\ell:\Delta_\ell\to T_\ell.
$$
The map $g_\ell$ remains proper and schematically dominant,
and $\Delta_\ell$ still maps to $X$. For every $\ell$, the
induction hypothesis gives the specified
bundle $\mathcal V_{T_\ell}$ and its comparison on $\Delta_\ell$.

Define algebraic coherent sheaves and maps
$$
\begin{aligned}
 E_\ell&=i_{\ell,*}g_{\ell,*}
               (M_\Delta|_{\Delta_\ell}),&
 r_\ell&:E\longrightarrow E_\ell,\\
 F_\ell&=\im(i_{\ell,*}\mathcal V_{T_\ell}\longrightarrow E_\ell),&
 K_\ell&=\ker r_\ell.
\end{aligned}
$$
Restriction of a locally free sheaf to a closed subspace is
surjective as a map of sheaves. Consequently,
\begin{equation}\label{aq:line-thickening-images}
 r_\ell^{\mathrm{an}}(N_T)=F_\ell^{\mathrm{an}},
 \qquad
 (r_\ell^{\mathrm{an}})^{-1}(F_\ell^{\mathrm{an}})
       =N_T+K_\ell^{\mathrm{an}}.
\end{equation}
The second equality follows by taking local lifts.

Lemma~\ref{aq:rees}, applied to the coherent sheaf $M_\Delta$,
gives
\begin{equation}\label{aq:line-cofinality}
 K_\ell=g_*(J^\ell M_\Delta),\qquad
 \text{for every }n\text{ there is }\ell
       \text{ with }K_\ell\subset I^nE.
\end{equation}
By \eqref{aq:line-thickening-images}, $N_T+K_\ell^{\mathrm{an}}$
is a specified algebraic submodule of $E^{\mathrm{an}}$.
Adding $(I^nE)^{\mathrm{an}}$ and using
\eqref{aq:line-cofinality} shows that
$$
 N_T+(I^nE)^{\mathrm{an}}
$$
is algebraic for every $n$. Lemma~\ref{aq:coherent-recovery}
algebraizes $N_T$ inside $E$. Faithful flatness of
analytification then gives a locally free model $\mathcal V_T$
of the prescribed rank.
Its adjunction map
$$
 g^*\mathcal V_T\longrightarrow M_\Delta
$$
analytifies to the specified isomorphism; its algebraic kernel and
cokernel therefore vanish. This proves the induction.
Uniqueness follows from the same faithful flatness.
\end{proof}

\subsection{Functoriality under pullback to smooth bases}

Apply Proposition~\ref{aq:hodge-descent} to
$$
 T=Y,\qquad \Delta=X=M,\qquad g=q,\qquad h=1_M.
$$
Since $Y$ is reduced, the proper surjection $q$ is
schematically dominant. We obtain
\begin{equation}\label{aq:actual-line-model}
 L^{\mathrm{an}}\simeq i^*\Lambda,\qquad
 q^*L\xrightarrow{\;\sim\;}L_M.
\end{equation}

\begin{corollary}\label{aq:filtered-algebraization}
The filtered bundle with polarization $i^*(\mathcal H,F^\bullet,Q)$
has an algebraic realization $(\mathcal H_Y,F_Y^\bullet,Q_Y)$
on $Y$.
Each $F_Y^p\subset\mathcal H_Y$ is a
subbundle with locally free graded quotients, $Q_Y$ is perfect,
and the natural comparison is algebraic:
\begin{equation}\label{aq:filtered-source-model}
 q^*(\mathcal H_Y,F_Y^\bullet,Q_Y)
 \xrightarrow{\;\sim\;}(\mathcal H_M,F_M^\bullet,Q_M).
\end{equation}
Moreover,
$$
 (F_Y^p)^{\perp_{Q_Y}}=F_Y^{w-p+1},\qquad
 L\simeq\bigotimes_{p=a+1}^{b}\det F_Y^p.
$$
\end{corollary}

This realization is defined also at the nonnormal points of $Y$.

\begin{proof}
Apply Proposition~\ref{aq:hodge-descent} to $\mathcal H$ and
each $F^p$, keeping their models inside the common coherent
sheaf $q_*\mathcal H_M$. The composite
$$
 F_Y^p\longrightarrow q_*F_M^p\longrightarrow q_*\mathcal H_M
$$
has zero class in
$q_*\mathcal H_M/\mathcal H_Y$ after analytification, hence
factors algebraically through $\mathcal H_Y$. The same argument
gives the specified inclusions $F_Y^{p+1}\subset F_Y^p$.
Both $\mathcal H_Y/F_Y^p$ and $F_Y^p/F_Y^{p+1}$ are locally
free: their analytifications are the tautological quotient bundles,
and analytic local rings are faithfully flat over the algebraic ones.

The polarization is an analytic section of
$\mathcal Hom(\mathcal H_Y\otimes\mathcal H_Y,\cO_Y)$.
Its pullback to the proper schematically dominant source $M$
is the algebraic $Q_M$. Apply Lemma~\ref{aq:graph-descent} to
the total space of this Hom bundle to algebraize that section.
Faithful analytification verifies its perfection, symmetry,
filtration orthogonality and compatibility with pullback to $M$.
Applying the same argument to the determinant tensor and
\eqref{aq:actual-line-model} identifies it with the specified $L$.
\end{proof}

\begin{lemma}\label{aq:smooth-pullback}
Let $u:V\to Y$ be an algebraic map from a smooth separated
finite-type algebraic space. The pullback from $\Gamma\backslash D$
is a polarized integral variation, and its canonical comparison with the
algebraic filtered bundle is algebraic:
\begin{equation}\label{aq:all-filtered-comparisons}
 u^*(\mathcal H_Y,F_Y^\bullet,Q_Y)
 \xrightarrow{\;\sim\;}(\mathcal H_V,F_V^\bullet,Q_V).
\end{equation}
In particular,
\begin{equation}\label{aq:all-smooth-comparisons}
 u^*L\xrightarrow{\;\sim\;}L_V
       =\bigotimes_{p=a+1}^{b}\det F_V^p.
\end{equation}
These comparisons commute with composition and base change.
\end{lemma}

\begin{proof}
Work first on a smooth affine \'{e}tale chart of $V$. Then
$$
 E_0=(M\times_YV)_{\mathrm{red}},\qquad
 v:E_0\to V,\qquad w:E_0\to M
$$
is a quasi-projective scheme: separatedness of $Y$ embeds it as
a closed subspace of $M\times V$. The map $v$ is proper and
surjective. Resolve $E_0$ and choose components dominating each
component of $V$. Generic smoothness gives a submersion onto a
dense open subset of $V$. The filtration pulled back from $\Gamma\backslash D$
is horizontal there because its pullback factors through the
variation on $M$. The morphism induced by the flat connection,
$$
 F_V^p\longrightarrow
  (\mathcal H_V/F_V^{p-1})\otimes\Omega_V^1
$$
is holomorphic and vanishes on that dense open set, hence vanishes
everywhere. The remaining polarized integral Hodge data pull back from
$\Gamma\backslash D$, so they define a variation of Hodge
structure on $V$.
Griffiths' algebraicity theorem gives its algebraic Hodge bundles.

Pull the natural horizontal comparison between
$w^*\mathcal H_M$ and $v^*\mathcal H_V$ to a proper smooth
quasi-projective resolution $r:\widetilde E_0\to E_0$.
Regular singularities are preserved by algebraic pullback, so
Deligne's full faithfulness
\cite[II, Proposition~4.6 and Theorem~5.9]{Deligne} makes the
comparison algebraic there. Since $E_0$ is reduced, $r$ is
schematically dominant. Apply Lemma~\ref{aq:graph-descent} to the
comparison as a map into the total Hom bundle over $E_0$.
The comparison descends algebraically, and
\eqref{aq:filtered-source-model} gives
$$
 v^*u^*\mathcal H_Y\simeq w^*\mathcal H_M\simeq v^*\mathcal H_V.
$$
A second application of the same lemma, now to the proper
schematically dominant map $v$ and the total bundle
$\mathcal Hom(u^*\mathcal H_Y,\mathcal H_V)$, algebraizes the
comparison on $V$. Faithful analytification verifies
its inverse and preservation of the filtration and polarization.
Taking determinants gives \eqref{aq:all-smooth-comparisons};
uniqueness gives descent between \'{e}tale charts and compatibility
with composition and base change.
\end{proof}

\begin{lemma}\label{aq:schmid-functoriality}
Let $\bar V$ and $\bar V'$ be smooth projective varieties with
simple normal crossings divisors $D_V$ and $D_{V'}$. Suppose
their complements carry variations obtained from algebraic maps
to $Y$, and that
$$
 \varphi:\bar V'\to\bar V,\qquad
 \varphi(V')\subset V,\qquad V=\bar V\setminus D_V,
$$
identifies the variations on the complements. Extending the Hodge
filtrations gives canonical extensions of the Griffiths line bundles
satisfying
\begin{equation}\label{aq:schmid-pullback}
 \varphi^*\bar L_V\xrightarrow{\;\sim\;}\bar L_{V'}.
\end{equation}
This is the extension of the specified comparison
\eqref{aq:all-smooth-comparisons}.
\end{lemma}

\begin{proof}
The integral monodromy theorem gives quasi-unipotent boundary
monodromy; neatness makes it unipotent. Locally write
$$
 t_i\circ\varphi=u_i\prod_j s_j^{a_{ij}},
 \qquad
 N'_j=\sum_i a_{ij}N_i,
$$
where the $u_i$ are units. The commuting $N_i$ are nilpotent,
so each $N'_j$ is nilpotent. In Deligne frames, the canonical
exponents are zero on both sides. The untwisted filtration
$$
 \exp\!\left(-\sum_i\frac{\log t_i}{2\pi\sqrt{-1}}N_i\right)
 F^\bullet
$$
pulls back to the corresponding filtration in the $s_j$
coordinates, up to the holomorphic invertible factor from the
$\log u_i$. Its holomorphic extension, supplied by the nilpotent
orbit theorem \cite[Theorem~(4.12)]{Schmid1973}, therefore gives the
same extended filtration. Taking the determinant tensor in
\eqref{aq:griffiths-line} proves
\eqref{aq:schmid-pullback}. Projective GAGA makes this comparison
algebraic.
\end{proof}

\subsection{A general ampleness criterion}

The following version of the argument in \cite[\S5]{BBT} includes
the pullback compatibility needed to compare auxiliary divisors with
their images.

\begin{theorem}[Ampleness criterion, after Bakker--Brunebarbe--Tsimerman]
\label{aq:general-ampleness}
Let $Y$ be a reduced separated algebraic space of finite type over $\C$,
and let $L$ be a line bundle on $Y$. Suppose that for every algebraic
map $u:V\to Y$ with $V$ smooth and quasi-projective, and every smooth
projective compactification $\bar V$ with simple normal crossings
boundary $B=\bar V\setminus V$, an extension $\bar L_V$ of $u^*L$
is specified. Assume compatibility with pullback by every morphism
of such pairs whose restriction to the interiors is over $Y$, using
the specified identifications on the interiors.
If $\bar L_V$ is nef and big on each component whenever
$V\to Z$ is a proper birational resolution of a reduced closed
subspace $Z\subset Y$, then $Y$ is a scheme and $L$ is ample.
Consequently, a finite-dimensional subspace of $H^0(Y,L^m)$ defines
a locally closed immersion for some $m>0$.
\end{theorem}

\begin{proof}
For a closed subspace $T\subset Y$, possibly nonreduced, call a
section of $L^m|_T$ \emph{boundary-vanishing} if its restriction to
$T_{\mathrm{red}}$, pulled back to a log resolution, extends to a
section vanishing on the reduced boundary. Common resolutions and
pullback compatibility make this independent of the resolution.
We prove by induction on dimension that, for every finite reduced
$P\subset T$, boundary-vanishing sections of some positive power
interpolate on $P$.

We first record three facts used in the induction.
First, interpolation implies ampleness. By Noetherianity of $T$ and
$T\times T\setminus\Delta_T$, finitely many homogeneous sections
have no common zero and separate points. Taking Veronese
monomials puts them in one degree without losing separation.
Their map to projective space is quasi-finite;
Zariski's main theorem factors it as an open immersion followed by
a finite map. Thus $T$ is a scheme and $L|_T$ is ample.

Second, interpolation passes from the reduction to any nilpotent
thickening. On $T_{\mathrm{red}}$, ampleness provides a section $a$
whose nonvanishing locus is affine and contains $P$. Choose a
boundary-vanishing section $h$ nonzero on $P$. For $N$ large,
$g=ah^N$ is boundary-vanishing: on a log resolution the zeros of
$h$ cancel the boundary poles of $a$. Its nonvanishing locus is
affine, since it is the principal open defined by
$h^{\deg a}/a^{\deg h}$ in $T_a$. Choose boundary-vanishing
sections with delta values on $P$, and multiply them by $g$;
each product has an affine nonvanishing locus. If $T_0\subset T_1$ is a
square-zero thickening with ideal $I$ and $s\in H^0(T_0,L^b)$ has
an affine nonvanishing locus, a sufficiently high power of $s$ lifts to
$T_1$. Indeed its obstruction is
$\delta(s^r)=r s^{r-1}\delta(s)\in H^1(T_0,I\otimes L^{br})$,
and multiplication by a high power of $s$ kills $\delta(s)$ because
cohomology localized at $s$ is the cohomology of the affine
$(T_0)_s$. Affineness persists under nilpotent thickenings, so this applies
along a finite square-zero filtration.
Common powers preserve the vanishing pattern on $P$ and boundary
vanishing; rescaling restores the prescribed nonzero values.

Third, powers extend across the boundary even on nonreduced schemes.
Let $A$ be a finite-type $\C$-algebra, let $f$ be a nonzerodivisor,
and let $s\in A[f^{-1}]$. Suppose that after a proper surjective map
from a normal scheme to $\operatorname{Spec}A_{\mathrm{red}}$, the pullback of
$s_{\mathrm{red}}$ is regular and vanishes over $f=0$.
Then $s$ is integral over $A$: coherence of proper direct images
gives integrality on the reduction, and a lifted monic relation
becomes zero after taking a power. Thus $B=A[s]$ is finite over $A$,
and $f^cB\subset A$ for some $c$. The induced proper map to $\operatorname{Spec}B_{\mathrm{red}}$
is surjective, so boundary vanishing gives $s^e\in fB$ for some $e$.
For every $k$,
$$
 s^{e(c+k)}\in f^{c+k}B\subset f^kA.
$$
On a finite affine cover, a common divisible power therefore
extends with any prescribed boundary order.

The zero-dimensional case is immediate. By the second fact, it
suffices in dimension $d>0$ to treat reduced $T$. Remove its
zero-dimensional connected components and take a log resolution
$$
 f:X=\bar X\setminus D\longrightarrow T,\qquad
 M=f_*\cO_X.
$$
Let $I$ define the reduced nonisomorphism locus; then
$I^a(M/\cO_T)=0$ for some $a$. The proper Rees argument in the
proof of Lemma~\ref{aq:rees} gives
$f_*(I^n\cO_X)\subset I^{n-c}M$ for one $c$ and all large $n$.
Fix $n$ with $f_*(I^n\cO_X)\subset\cO_T$ and put $S=V(I^n)$.
A section of $f^*L^m$ descends to $T$ if its restriction to
$f^{-1}S$ comes from $S$: locally lift that section from $S$ and
subtract; the difference lies in $f_*(I^n\cO_X)\otimes L^m$.
This also proves uniqueness. Put $Z=S\cup P$ and $W=f^{-1}Z$,
both scheme-theoretically.

Write $\bar L$ for the specified extension on $\bar X$, and choose
an ample divisor $A$. Bigness on every component gives an effective
divisor $E_0\in|\bar L^b(-A)|$ containing the scheme-theoretic closure of $W$.
For a fixed sufficiently large $r$, Fujita vanishing
\cite[Theorem~1.4.35]{Lazarsfeld2004} gives, for all $m\ge br$,
$$
 H^1\bigl(\bar X,\bar L^{m-br}(-D+rA)\bigr)=0.
$$
Hence sections of $\bar L^m(-D)|_E$, where $E=rE_0$, lift to
$\bar X$. Write $E=\bar E'+D'$, with $D'$ supported on $D$ and
no component of $\bar E'$ contained in $D$. For a fixed $k$ with
$kD\ge D+D'$, the exact sequence
$$
 0\longrightarrow\cO_{\bar E'}(-D')
 \longrightarrow\cO_E\longrightarrow\cO_{D'}\longrightarrow0
$$
embeds $\bar L^m(-kD)|_{\bar E'}$ into
$\bar L^m(-D)|_E$. Thus a section on $\bar E'$ with this boundary
order lifts and descends to a boundary-vanishing section on $T$
whenever its restriction to $W$ comes from $Z$.

Put $E'=\bar E'\cap X$ and
$F=\im_{\mathrm{sch}}(E'\to T)\cup Z$. This is a closed subspace
of dimension less than $d$, and $E'\to F$ is proper. By induction,
choose boundary-vanishing sections of one degree on $F$ with delta
values on $P$. Pull them back to $E'$. To compare extensions,
choose a log resolution $V_F\to F_{\mathrm{red}}$ and resolve
components of $E'_{\mathrm{red}}\times_{F_{\mathrm{red}}}V_F$
dominating $E'_{\mathrm{red}}$. The resulting smooth source is
proper over both spaces. Take a smooth projective compactification
with simple normal crossings boundary resolving both maps to
$\bar X$ and $\bar V_F$; the first factors through
$\bar E'_{\mathrm{red}}$. By properness, its boundary maps into
$D$ and the boundary of $\bar V_F$. Pullback compatibility
therefore makes each section regular, with zero restriction to the
boundary, on a proper surjective model of $\bar E'_{\mathrm{red}}$.
Since $\bar E'$ is a Cartier divisor in the smooth $\bar X$, it has
no embedded components; hence $D|_{\bar E'}$ is Cartier. The third
fact extends a common sufficiently large power to
$\bar L^m(-kD)|_{\bar E'}$, with $m\ge br$. Its restriction to $W$
comes from $Z\subset F$, so the preceding paragraph lifts and descends it.
After rescaling, these sections interpolate on $P$, completing
the induction. The first fact gives ampleness; a sufficiently high
power then has a finite system defining a locally closed immersion.
\end{proof}

\subsection{Positivity on all subspaces}

\begin{proof}[Proof of Theorem~\ref{aq:quasiprojectivity}]
Theorem~\ref{aq:algebraization} gives the factorization
\eqref{aq:projective-factor} and the proper surjection
\eqref{aq:smooth-source}. Equations
\eqref{aq:actual-line-model}--\eqref{aq:all-smooth-comparisons}
give the algebraic Griffiths line bundle and its compatibility
with pullback to smooth bases. It remains to prove ampleness.

Let $Z\subset Y$ be any reduced closed subspace, and take a
projective log resolution
$$
 r:V=\bar V\setminus D_V\longrightarrow Z,
$$
proper and birational on each irreducible component. Such a
resolution exists by Chow's lemma, compactification and resolution.
The preceding lemmas identify
$$
 r^*(L|_Z)\simeq L_V
$$
with the Griffiths line bundle and provide its functorial canonical
extension $\bar L_V$. Since $i$ is a closed analytic immersion
and $r$ is generically an isomorphism, the period map of each
positive-dimensional component of $V$ is generically immersive.
The classical positivity statement
\cite[Lemma~6.15]{BBT} therefore gives, component by component,
\begin{equation}\label{aq:nef-big}
 \bar L_V\ \text{nef},\qquad
 c_1(\bar L_V)^{d}
     =\int_V c_1(L_V,h)^{d}>0,\qquad d=\dim V>0.
\end{equation}
Thus $\bar L_V$ is big; zero-dimensional components are immediate.
Lemma~\ref{aq:schmid-functoriality} supplies the pullback compatibility
for all smooth sources required by Theorem~\ref{aq:general-ampleness}.
Together with \eqref{aq:nef-big}, that theorem shows that $Y$ is a
scheme and $L$ is ample. It also gives the asserted locally closed
immersion $j$ from a finite system in $H^0(Y,L^m)$.
Its scheme-theoretic projective closure carries the ample extension
$\cO_{\bar Y}(1)$.

Finally, if $\Phi$ is proper, the algebraization theorem
identifies its image with $Y^{\mathrm{an}}$ and
gives a proper surjective $f$. Otherwise Chevalley's theorem
applied to the algebraic dominant map $f$ gives a dense
constructible image.
\end{proof}

\subsection{Arithmetic quotients with stabilizers}

For an arithmetic group that is not neat, we require period maps
to be locally liftable to $D$. This holds for maps induced by
polarized integral variations.

\begin{theorem}\label{aq:finite-level}
The factorization and quasi-projectivity assertions of
Theorem~\ref{aq:quasiprojectivity} hold for any arithmetic
$\Gamma$, for an everywhere defined locally liftable period map
from a reduced separated finite-type algebraic space.
The Griffiths $\mathbb Q$-line bundle on the reduced closed image
is algebraic and ample.

More precisely, choose a normal neat subgroup
$\Gamma'\triangleleft\Gamma$ of finite index, and put
$$
 G=\Gamma/\Gamma'.
$$
There are a $G$-invariant reduced closed image
$Y'\subset\Gamma'\backslash D$,
an algebraic finite quotient $\pi:Y'\to Y=Y'/G$, and an ample
line bundle $L_0$ with its specified analytic identification such that
\begin{equation}\label{aq:coarse-line}
 \pi^*L_0\simeq(L')^{|G|},\qquad
 L=\frac{1}{|G|}[L_0]\in\operatorname{Pic}(Y)\otimes\mathbb Q.
\end{equation}
Here $L'$ is the Griffiths line bundle on $Y'$.
\end{theorem}

The quotient identification includes the nonnormal analytic local
rings. Stabilizers may act nontrivially on the individual Hodge bundles.

\begin{proof}
\emph{The level cover.}
First let $S$ be smooth quasi-projective and
$\phi:S^{\mathrm{an}}\to \Gamma\backslash D$ locally liftable. Form
\begin{equation}\label{aq:normalized-level-cover}
 E=(S^{\mathrm{an}}\times_{\Gamma\backslash D}(\Gamma'\backslash D))_{\mathrm{red}},\qquad
 E^\nu\longrightarrow S^{\mathrm{an}}.
\end{equation}
This map is finite. Near a point of $S$, choose a connected
normal neighborhood $U$ and a lift $F:U\to V\subset D$,
where the coarse quotient chart is $V/H$, $H$ finite.
Over every chart of $\Gamma'\backslash D$ above $V/H$, the reduced pullback
is a finite union of graphs
$$
 \{(u,hF(u)):u\in U\},\qquad h\in H,
$$
with repeated graphs counted once. Normalization separates the
distinct graphs, including where they meet. Thus
$$
 (E^\nu)|_U\simeq\coprod U.
$$
Consequently \eqref{aq:normalized-level-cover} is finite \'etale.
Riemann existence \cite[Expos\'e XII]{SGA1} algebraizes it as
$p:S'\to S$; its $G$-action is algebraic by full faithfulness
for finite \'etale covers. The map to $\Gamma'\backslash D$ is the period map of
a polarized integral variation on the smooth quasi-projective
$S'$. The $G$-action on this cover may have stabilizers.

\emph{The image and the group action.}
For the original reduced source choose a proper surjection
$r:S\to B$ from a finite disjoint union of smooth
quasi-projective schemes, isomorphic over a dense regular open
of each component. Horizontality pulls back by continuity from
that open. Apply the preceding construction to $\Phi r^{\mathrm{an}}$,
and Theorem~\ref{aq:quasiprojectivity} to $S'\to \Gamma'\backslash D$.
It gives an algebraic map $f':S'\to Y'$, the reduced closure
$Y'$, and a proper surjection $q:M\to Y'$ from smooth
quasi-projective varieties carrying the pulled-back variation.
The set $Y'$ is $G$-invariant.

Fix $g\in G$. The translate of the period map of $M$ is
again a neat period map, since $\Gamma'$ is normal. Apply
Lemma~\ref{aq:relation} to this translate and the original map:
$$
 R_g^{\mathrm{an}}
 =\{(m_1,m_2):gq(m_1)=q(m_2)\}_{\mathrm{red}},\qquad
 a=\operatorname{pr}_1:R_g\to M.
$$
The algebraic map $a$ is proper, surjective and schematically
dominant. Applying Lemma~\ref{aq:graph-descent} first to $a$ and
then to $q$ gives
$$
 gq:M\longrightarrow Y'\ \text{algebraic},\qquad
 g:Y'\longrightarrow Y'\ \text{algebraic}.
$$
These morphisms define an algebraic $G$-action on $Y'$.

\emph{The quotient and its structure sheaf.}
Every finite orbit in a quasi-projective scheme has an affine
neighborhood: take a projective closure and a section of a
sufficiently high power of its ample line bundle, vanishing on the boundary and
nonzero on the orbit.
Intersect its finitely many translates. Separatedness makes this
intersection affine, and it is $G$-invariant. The affine quotients
$$
 \operatorname{Spec}A\longmapsto\operatorname{Spec}A^G
$$
glue to a separated finite-type scheme $Y=Y'/G$ with finite
surjective $\pi$. Finite coherent comparison and exactness of
$G$-invariants in characteristic zero identify
\begin{equation}\label{aq:actual-invariant-rings}
 \cO_{Y^{\mathrm{an}}}
 =\bigl(\pi_*^{\mathrm{an}}\cO_{Y'^{\mathrm{an}}}\bigr)^G.
\end{equation}
The same exactness applied to the invariant ideal of
$Y'^{\mathrm{an}}\subset\Gamma'\backslash D$ makes
$Y^{\mathrm{an}}\to(\Gamma'\backslash D)/G=\Gamma\backslash D$
a closed analytic immersion: writing $\rho:\Gamma'\backslash D\to \Gamma\backslash D$,
$$
 0\longrightarrow(\rho_*\mathcal I_{Y'})^G
 \longrightarrow\cO_{\Gamma\backslash D}
 \longrightarrow(\rho_*\cO_{Y'^{\mathrm{an}}})^G
 \longrightarrow0.
$$
Finiteness of $\rho$ gives
$\rho(Y')=\overline{\Phi(B^{\mathrm{an}})}$. Equation
\eqref{aq:actual-invariant-rings} identifies its reduced analytic
structure, including its nonnormal local rings.

The analytic factor $S^{\mathrm{an}}\to Y^{\mathrm{an}}$
pulls back along $p$ to the algebraic $\pi f'$.
Lemma~\ref{aq:graph-descent}, first for $p$ and then for $r$,
therefore algebraizes the original map $f:B\to Y$.
Its dominance, properness when $\Phi$ is proper, and the
constructibility assertion follow as in
Theorem~\ref{aq:quasiprojectivity}.

\emph{Descent of the Griffiths line bundle.}
The tautological analytic action on the Griffiths line gives isomorphisms
$\alpha_g:g^*L'\to L'$ satisfying the cocycle law.
Pull them back to $M$. Lemma~\ref{aq:smooth-pullback} identifies the
two pullbacks with the Griffiths line bundles of the original and
translated variations. Their tautological identification comes
from a horizontal integral isometry; Deligne's full faithfulness
\cite[II, Proposition~4.6 and Theorem~5.9]{Deligne} algebraizes it,
and hence its determinant tensor. Apply
Lemma~\ref{aq:graph-descent} to the corresponding Hom bundle over
$Y'$. It algebraizes each $\alpha_g$. Faithful analytification
verifies the cocycle law, giving the specified algebraic
$G$-linearization of $L'$.

Every stabilizer acts trivially on $(L')^{|G|}$. Set
$$
 L_0=\bigl(\pi_*((L')^{|G|})\bigr)^G.
$$
To check descent, work on an invariant affine neighborhood of an
orbit. Choose a section of $(L')^{|G|}$ nonzero at one orbit point and
zero at the others; average over $G$. Triviality of the
stabilizer action makes the average nonzero on the entire orbit.
On its nonvanishing locus it is an invariant frame. Hence
$$
 L_0\ \text{is invertible},\qquad
 \pi^*L_0\xrightarrow{\;\sim\;}(L')^{|G|}.
$$
These analytic invariant frames identify $L_0^{\mathrm{an}}$
with the corresponding power of the Griffiths
$\mathbb Q$-line bundle.

To prove ampleness, fix a point $y\in Y$ and an open neighborhood
$W$. Using a projective closure of $Y'$, choose a section
$t\in H^0(Y',(L')^n)$ such that
$$
 \pi^{-1}(y)\subset Y'_t\subset\pi^{-1}(W),
 \qquad Y'_t\ \text{affine}.
$$
Its invariant product descends to a section
$$
 \prod_{g\in G}g\cdot t=\pi^*s,\qquad
 s\in H^0(Y,L_0^n),\qquad
 \pi^{-1}(Y_s)=\bigcap_{g\in G}gY'_t.
$$
The last intersection is affine and invariant, so its finite quotient
$Y_s$ is affine. These opens form a basis; the affine-open
criterion for ampleness proves that $L_0$ is ample.
The specified descents agree after common neat refinements,
giving \eqref{aq:coarse-line} and the claimed quasi-projectivity.
\end{proof}

Together with Corollary~\ref{aq:filtered-algebraization} and
Lemma~\ref{aq:smooth-pullback} in the neat case, this completes the
proof of Theorem~\ref{intro:algebraicity}.

\subsection{Jet interpolation and affine charts}

Fix an immersion supplied by
Theorem~\ref{aq:quasiprojectivity}:
$$
 j:Y\hookrightarrow\mathbf P(V^\vee),\qquad
 V\subset H^0(Y,L^{m_0}),\qquad
 \cO_{\bar Y}(1)|_Y\simeq L^{m_0}.
$$
Put
$$
 R_d=\im\bigl(\Sym^d V\longrightarrow H^0(Y,L^{m_0d})\bigr).
$$
This embedding gives the following interpolation bound.

\begin{proposition}\label{pol:finite-jets}
For every finite closed subscheme $T\subset Y$ of length
$\ell\ge1$,
\begin{equation}\label{pol:jet-evaluation}
 R_d\longrightarrow H^0(T,L^{m_0d}|_T)
 \quad\text{is surjective for }d\ge\ell-1.
\end{equation}
The same holds for finite analytic subschemes of $Y^{\mathrm{an}}$.
If $\mathcal T\subset Y\times S$ is finite flat of rank $\ell$
over an algebraic base $S$, then, for the same $d\ge\ell-1$,
\begin{equation}\label{pol:family-evaluation}
 R_d\otimes\cO_S\longrightarrow
 p_*\bigl(\operatorname{pr}_Y^*L^{m_0d}|_{\mathcal T}\bigr)
 \longrightarrow0
\end{equation}
is a surjection onto a vector bundle of rank $\ell$.
It has algebraic right inverses Zariski locally, and a global
algebraic right inverse when $S$ is affine.
\end{proposition}

The subscheme $T$ may be nonreduced and supported at singular points
of $Y$.

\begin{proof}
Choose a basis $x_0,\ldots,x_N$ of $V$ with $x_0$
nonzero at every support point of $T$. The immersion of the
proper scheme $T$ in projective space is closed. Thus, on
$x_0\ne0$,
$$
 B_T=H^0(T,\cO_T)=\C[z_1,\ldots,z_N],\qquad
 z_i=x_i/x_0,
$$
where the right side denotes the subalgebra generated by these
elements of $B_T$. Set
$$
 W_d=\operatorname{span}\{z^\alpha:|\alpha|\le d\}\subset B_T.
$$
Since $1\in W_0$,
$$
 W_d=W_{d+1}\ \Longrightarrow\ z_iW_d\subset W_d\ (1\le i\le N)
 \ \Longrightarrow\ W_d=B_T.
$$
Thus the dimensions strictly increase until they reach
$\dim B_T=\ell$, so $W_{\ell-1}=B_T$.
Homogenization by $x_0^d$ proves \eqref{pol:jet-evaluation}.

The defining ideal of a local analytic quotient of length $\ell$
contains $(\mathfrak m_y^{\mathrm{an}})^\ell$, and
$$
 \cO_{Y,y}/\mathfrak m_y^\ell
 \simeq
 \cO_{Y^{\mathrm{an}},y}/(\mathfrak m_y^{\mathrm{an}})^\ell.
$$
This ideal therefore algebraizes uniquely; glue over its finite
support and apply the first assertion.

For the family, the target of \eqref{pol:family-evaluation} is
locally free and commutes with base change. The first assertion
on geometric fibers and Nakayama's lemma give surjectivity.
A surjection onto a vector bundle splits locally; on an affine
base its module of sections is projective, hence it splits globally.
\end{proof}

For jets of orders $k_\nu$ at distinct points $y_\nu$, one may use
$$
 \ell=\sum_\nu\operatorname{length}
       (\cO_{Y,y_\nu}/\mathfrak m_{y_\nu}^{k_\nu+1})
 \le\sum_\nu\binom{e_\nu+k_\nu}{e_\nu},\qquad
 e_\nu=\dim_\C\mathfrak m_{y_\nu}/\mathfrak m_{y_\nu}^2.
$$
The degree bound is relative to the chosen $m_0$ and $V$.

\begin{corollary}\label{pol:affine-atlas}
There are $r>0$ and finitely many sections
$s_i\in H^0(Y,L^r)$, nonzero on every irreducible component,
such that
\begin{equation}\label{pol:affine-charts}
 Y=\bigcup_i U_i,\qquad U_i=Y_{s_i}\ \text{affine},\qquad
 L^r|_{U_i}\simeq\cO_{U_i},\qquad
 U_i\cap U_j=D_{U_i}(s_j/s_i).
\end{equation}
All finite intersections are affine and their analytifications
are Stein. Each zero divisor $(s_i=0)$ is Cartier.
\end{corollary}

\begin{proof}
Let $B_\infty=(\bar Y\setminus Y)_{\mathrm{red}}$.
For $k\gg0$, Serre's global generation theorem
\cite[II, Theorem~5.17]{Hartshorne1977} makes
$\mathcal I_{B_\infty}\otimes\cO_{\bar Y}(k)$ globally generated.
Choose a basis $\widetilde s_i$ of its section space avoiding
the finitely many proper subspaces of sections vanishing
identically on a component of $\bar Y$. Then
$$
 \bigcup_i\bar Y_{\widetilde s_i}=Y,\qquad
 \bar Y_{\widetilde s_i}\subset Y,\qquad
 s_i=\widetilde s_i|_Y,\quad r=m_0k.
$$
Since $\cO_{\bar Y}(k)$ is very ample, each complement
$\bar Y_{\widetilde s_i}$ is affine. Formula
\eqref{pol:affine-charts} gives the intersections.
Reducedness and nonvanishing on every component make
$s_i$ a non-zero-divisor in each local trivialization.
\end{proof}

These charts also carry potentials for the Hodge curvature.
For the Hodge metric $h$ on $L$, with
$\omega_L=i\Theta_h(L)$, put
\begin{equation}\label{pol:local-potential}
 \phi_i=-\frac1r\log\|s_i\|_{h^r}^2,\qquad
 i\partial\bar\partial\phi_i=\omega_L,\qquad
 \phi_j-\phi_i=-\frac1r\log|s_j/s_i|^2.
\end{equation}
Locally extend the nonvanishing $s_i$ to the Griffiths
line bundle in a quotient chart. Horizontal positivity and
Lemma~\ref{lem:singular-horizontal-psh} make $\phi_i$
strictly plurisubharmonic on $U_i^{\mathrm{an}}$, including its
singular points.

\begin{example}\label{aq:canonical-extension}
The canonical extension of the Griffiths line bundle need not be
ample on a prescribed compactification. Consider the product
variation of the two universal elliptic families over $Y(N)^2$,
$N\ge3$, with its
ordered splitting tensors and fixed product period quotient. Put
$$
 \bar X=X(N)^2,\qquad
 \bar L=\lambda\boxtimes\lambda.
$$
The modular discriminant is a section of $\lambda^{12}$, has
no zero in the open curve, and has vanishing order $w_c>0$
at a cusp of width $w_c$. Thus
$$
 \lambda^{12}\simeq
 \cO_{X(N)}\left(\sum_{\text{cusps }c}w_c[c]\right),
 \qquad \deg\lambda>0.
$$
Hence $\bar L$ is ample. Blow up an intersection of two boundary
components:
$$
 \pi:\widetilde X\longrightarrow\bar X,\qquad
 E\simeq\mathbf P^1.
$$
The open base is unchanged. By
Lemma~\ref{aq:schmid-functoriality} the canonical extension is
$\pi^*\bar L$, but
$$
 (\pi^*\bar L)|_E\simeq\cO_E,\qquad
 c_1(\pi^*\bar L)\cdot E=0.
$$
None of its positive powers is ample. On the other hand,
$\pi^*\bar L^m\otimes\cO_{\widetilde X}(-E)$ is ample for
$m\gg0$ and restricts to $L^m$ on the same open base.
\end{example}
\section{Splitting the Hodge filtration and recovering the period map}\label{sec:grading}

Let $H_{\mathbb C}$ be the underlying complex vector space of a
polarized Hodge structure of weight $w$, let $G$ be the connected
complex group preserving the polarization, and let
$\check D=G/P_o$ be the compact dual containing $D$. Put
$$
 S_o|_{H_o^{p,w-p}}=(2p-w)I,\qquad
 \mathcal A=G\cdot S_o=G/Z_G(S_o)\subset\mathfrak g.
$$
The polarization pairs complementary Hodge degrees, so
$S_o\in\mathfrak g$. The semisimple orbit $\mathcal A$ is closed
and affine. Its ordered eigenspaces define
$\varpi:\mathcal A\to\check D$.

\begin{theorem}\label{thm:grading-embedding}
The map $\varpi$ is a torsor under the cotangent bundle
$T^*\check D$. Every Stein complex subspace
$i:X\hookrightarrow\check D$ admits a holomorphic embedding
$$
 S:X\hookrightarrow\mathcal A\subset\mathfrak g,
 \qquad \varpi S=i.
$$
With $\mathcal P=\{p:h^{p,w-p}\ne0\}$,
the Hodge filtration is recovered by the polynomial identities
\begin{equation}\label{grad:spectral-recovery}
 \Pi_p(S)=\prod_{\substack{q\in\mathcal P\\q\ne p}}
       \frac{S-(2q-w)I}{2(p-q)},
 \qquad
 F^p=\bigoplus_{\substack{q\in\mathcal P\\q\ge p}}
       \im\Pi_q(S).
\end{equation}
\end{theorem}

Thus a holomorphic splitting gives matrix coordinates from which
the period map is recovered polynomially. The choice of splitting
is generally noncanonical.

\subsection{The affine bundle of compatible splittings}

Put $\mathfrak p=\operatorname{Lie}(P_o)$. Let $F\in\check D$.
A compatible grading of $F$ is a semisimple operator in
$\mathfrak g$ whose action on $\operatorname{Gr}_F^p$ is
multiplication by $2p-w$. Let $\mathfrak u_F$ be
the nilradical of the parabolic Lie algebra stabilizing $F$.

\begin{lemma}\label{lem:grading-torsor}
If $S_F$ is one compatible grading of $F$, then
$$
\varpi^{-1}(F)=S_F+\mathfrak u_F.
$$
These affine spaces form a torsor under $T^*\check D$.
\end{lemma}

\begin{proof}
Two compatible gradings induce the same scalar on each graded
piece. Their difference strictly raises the filtration and lies in
$\mathfrak g$, hence belongs to $\mathfrak u_F$.

Conversely, let $N\in\mathfrak u_F$. The diagonal blocks of
$S_F+N$ are distinct scalars, indexed by the degrees with nonzero
Hodge numbers. Conjugation by the unipotent radical successively
removes the off-diagonal blocks: at each positive filtration
degree, $\operatorname{ad}S_F$ is invertible, and higher
commutators affect only later degrees. Thus $S_F+N$ is conjugate
to $S_F$ by an element stabilizing $F$.

Writing $U$ for the unipotent radical of $P_o$, this gives
$$
\operatorname{Ad}(U)S_o=S_o+\mathfrak u,\qquad
\mathcal A\cong G\times_{P_o}(S_o+\mathfrak u).
$$
The linear part of the affine $P_o$-action is the adjoint action
on $\mathfrak u$. An invariant nondegenerate trace pairing identifies
$$
\mathfrak u=\mathfrak p^\perp\cong(\mathfrak g/\mathfrak p)^*.
$$
The associated direction bundle is therefore $T^*\check D$.
\end{proof}

\begin{proof}[Proof of Theorem \ref{thm:grading-embedding}]
Pull back the torsor of Lemma \ref{lem:grading-torsor} to $X$.
Differences of local sections define a \v{C}ech cocycle with values
in the locally free sheaf $i^*\Omega_{\check D}^1$. Since $X$ is Stein,
$$
H^1(X,i^*\Omega_{\check D}^1)=0
$$
by Cartan's Theorem B. Adjusting the local sections gives a global
holomorphic section, equivalently a map $S:X\to\mathcal A$ with
$\varpi S=i$.

Locally, $S$ is the graph of a section over the embedded subspace
$X$, so it is a holomorphic embedding. Lagrange interpolation
gives \eqref{grad:spectral-recovery}. The restriction of $\varpi$
is its inverse, preserving the original topology and complex-space
structure of $X$.
\end{proof}

The sections form an affine space over
$H^0(X,i^*\Omega_{\check D}^1)$. Their images are closed over
$X$, though not necessarily in $\mathcal A$.

\subsection{Algebraic gradings on affine charts of the period image}

Apply the same construction to the algebraic filtered bundle
$(\mathcal H_Y,F_Y^\bullet,Q_Y)$ on the quasi-projective image
from Corollary~\ref{aq:filtered-algebraization}. Put
$$
 \mathfrak g_{Q_Y}
 =\{A\in\mathcal End(\mathcal H_Y):A^{\dagger_{Q_Y}}=-A\},
 \qquad
 \mathfrak u_{F_Y}
 =\{A\in\mathfrak g_{Q_Y}:A(F_Y^p)\subset F_Y^{p+1}\ \forall p\}.
$$
These are algebraic vector bundles. The second is the nilradical
of the parabolic Lie algebra preserving the filtration, identified
with the pullback of the relative cotangent bundle of the flag
variety along the Hodge section. Here $\dagger_{Q_Y}$ denotes the
adjoint with respect to the bilinear polarization.

\begin{proposition}\label{pol:algebraic-grading}
Let $\mathcal P=\{p:\operatorname{Gr}_{F_Y}^p\mathcal H_Y\ne0\}$.
The compatible gradings
$$
 \mathscr G_{F_Y}
 =\{S\in\mathfrak g_{Q_Y}:
       (S-(2p-w)I)F_Y^p\subset F_Y^{p+1}\ \forall p\}
$$
form an algebraic torsor under $\mathfrak u_{F_Y}$.
It has an algebraic section $S_U$ on every affine open
$U\subset Y$. For each such section,
\begin{equation}\label{pol:algebraic-spectral-recovery}
 \begin{gathered}
 \prod_{p\in\mathcal P}(S_U-(2p-w)I)=0,\qquad
 \Pi_p(S_U)=
 \prod_{\substack{q\in\mathcal P\\q\ne p}}
       \frac{S_U-(2q-w)I}{2(p-q)},\\
 \mathcal H_Y|_U=\bigoplus_{p\in\mathcal P}\im\Pi_p(S_U),
 \qquad F_Y^p|_U=\bigoplus_{\substack{q\in\mathcal P\\q\ge p}}
                                  \im\Pi_q(S_U).
 \end{gathered}
\end{equation}
\end{proposition}

In particular, such sections exist on the finite affine cover
supplied by Corollary~\ref{pol:affine-atlas}, and their differences
on overlaps are sections of $\mathfrak u_{F_Y}$.

\begin{proof}
On an affine $U$, the locally free graded quotients are projective
modules. Successively split the filtration to obtain
$$
 \mathcal H_Y|_U=\bigoplus_{p\in\mathcal P}E^p,\qquad
 F_Y^p|_U=\bigoplus_{q\ge p}E^q,\qquad
 S_0|_{E^p}=(2p-w)I.
$$
Orthogonality $(F_Y^p)^\perp=F_Y^{w-p+1}$ implies that
$S_0^{\dagger_{Q_Y}}$ preserves the filtration. The perfect
pairing of the graded pieces of indices $p$ and $w-p$ gives
$$
 S_0^{\dagger_{Q_Y}}|_{\operatorname{Gr}_{F_Y}^p}
     =(w-2p)I.
$$
Thus
\begin{equation}\label{pol:polarization-correction}
 S_U=\frac12\bigl(S_0-S_0^{\dagger_{Q_Y}}\bigr)
\end{equation}
is $Q_Y$-skew-adjoint and induces the required scalar on each
graded piece. This proves existence on affine charts, including
singular and nonnormal ones. Two solutions differ by a section of
$\mathfrak u_{F_Y}$, and translation by every such section
preserves the defining conditions. In local frames adapted to the
filtration and polarization, the direction space is the nilradical.
This proves the torsor assertion.

Applying the commuting factors $S_U-(2p-w)I$ in increasing
order of $p\in\mathcal P$ annihilates the filtration. Since
the scalars $2p-w$ are distinct, the Lagrange polynomials in
\eqref{pol:algebraic-spectral-recovery} are consequently
orthogonal algebraic idempotents summing to $I$, and their
partial sums recover the filtration. These identities hold over
the coordinate ring.
\end{proof}

The torsor is defined in the endomorphism bundle
$\mathcal End(\mathcal H_Y)$, with direction bundle
$\mathfrak u_{F_Y}$. Its sections have regular matrix entries
in algebraic local frames. Under a change of frame, the matrix of
a section transforms by conjugation.

\subsection{The Griffiths line bundle and the grading orbit}

Suppose now that the Hodge structure is effective of weight
$w\geq1$, and put
$$
\Lambda=\bigotimes_{p=1}^{w}\det F^p.
$$
Let $\chi$ be the character by which $P_o$ acts on the reference
fiber. Since $\operatorname{tr}a=0$ for $a\in\mathfrak g$,
$$
d\chi(a)=\sum_p p\,\operatorname{tr}(a|_{\operatorname{Gr}_F^p})
=\frac12\operatorname{tr}(S_oa),\qquad a\in\mathfrak p.
$$
Use $b(A,B)=\frac12\operatorname{tr}(AB)$. The Atiyah sequence is
$$
0\longrightarrow\mathcal O_{\check D}
\longrightarrow\operatorname{At}(\Lambda)
\longrightarrow T_{\check D}\longrightarrow0;
$$
see \cite{Atiyah1957}. Choose the convention in which
$\operatorname{At}(\Lambda)_o=(\mathfrak g\oplus\mathbb C)/
\{(a,-d\chi(a)):a\in\mathfrak p\}$.
In the dual sequence, a functional taking the value $1$ on the
vertical direction descends precisely when its $\mathfrak g$-component
restricts to $d\chi$. Thus the affine fiber of the connection torsor is
$$
\{\mu\in\mathfrak g^*:\mu|_{\mathfrak p}=d\chi\}
=b(S_o+\mathfrak u,-).
$$
Its direction is $\mathfrak p^\perp$, so this convention identifies
the grading torsor with the connection torsor. The opposite convention
replaces the character by its negative and the line bundle by its dual.

On an open set $W\subset\check D$, a grading section therefore
specifies a holomorphic connection on $\Lambda|_W$. For a subspace
$i:X\hookrightarrow\check D$, the direction bundle pulls back to
$i^*\Omega_{\check D}^1$; the map
$i^*\Omega_{\check D}^1\to\Omega_X^1$ gives the intrinsic
connection data. These connections need not be flat.

For Hodge numbers $(1,g-1,g-1,1)$, the orbit has the explicit description
$$
\mathcal A=\left\{S\in\mathfrak{sp}(2g,\mathbb C):
(S^2-I)(S^2-9I)=0,\quad \operatorname{tr}(S^2)=2g+16\right\}.
$$
The square-free polynomial makes $S$ semisimple. The symplectic
pairing forces opposite eigenvalues to have equal multiplicities,
and the trace equation makes $3$ and $-3$ simple. In a fixed
symplectic basis,
$$
S=\begin{pmatrix}A&B\\ C&-A^{\mathsf T}\end{pmatrix},
\qquad B=B^{\mathsf T},\quad C=C^{\mathsf T}.
$$
Thus the entries give coordinates in $\mathbb C^{g(2g+1)}$, and
$$
F^3=\ker(S-3I),\qquad
F^2=\ker((S-3I)(S-I)),\qquad
F^1=(F^3)^{\perp_Q}.
$$
These formulas recover the Hodge filtration from the matrix entries.

\subsection{The horizontal differential}

Let $X\subset D$ be smooth and horizontal, and let $S$ be a
holomorphic grading. Write
$$
 H_x|_{H_x^{p,w-p}}=(2p-w)I,\qquad
 S=H+U,\qquad U=\sum_{k\ge1}U_k,
$$
where $U_k$ has Hodge degree $k$. Use the Hodge metric
$h_x$ and its Hilbert--Schmidt norm. With
$$
 \nabla=D_h+\theta+\theta^*,\qquad
 g_{\mathrm{GS}}(v,\bar z)=\tr(\theta_v\theta_z^*),
$$
the differential of $S$ has the following expression.

\begin{proposition}\label{grad:differential}
For $v\in T_x^{1,0}X$,
\begin{align}
 dS(v)&=2\theta_v+D'_{h,v}U+[\theta_v,U],\label{grad:dS}\\
 \|dS(v)\|_h^2
 &=4\|\theta_v\|_h^2+\|[\theta_v,U_1]\|_h^2\notag\\
 &\quad+\|D'_{h,v}U+[\theta_v,U_{\ge2}]\|_h^2.
 \label{grad:orthogonal}
\end{align}
At a fixed flag $F$, identify a horizontal tangent vector with
$\xi\in\mathfrak g_F^{-1}$. Then
\begin{equation}\label{grad:quotient-metric}
 \begin{aligned}
 \inf_{d\varpi_S(\eta)=\xi}\|\eta\|_{h_F}^2
 &=4\|\xi\|_{h_F}^2+\|[\xi,U_1]\|_{h_F}^2,\\
 \inf_{S\in\varpi^{-1}(F)}\inf_{d\varpi_S(\eta)=\xi}
       \|\eta\|_{h_F}^2&=4\|\xi\|_{h_F}^2.
 \end{aligned}
\end{equation}
\end{proposition}

\begin{proof}
Since $D_hH=0$ and $[\theta_v,H]=2\theta_v$, taking the
$(1,0)$ part of the flat derivative gives \eqref{grad:dS}.
The three summands
$$
 2\theta_v,\qquad [\theta_v,U_1],\qquad
 D'_{h,v}U+[\theta_v,U_{\ge2}]
$$
have degrees $-1,0,>0$, respectively, and are orthogonal.
For the quotient norm, $[\xi,S]$ is one lift of $\xi$, and
all lifts differ by
$\mathfrak u_F=\bigoplus_{k>0}\mathfrak g_F^k$. Minimization
cancels precisely its positive-degree component. The choice
$S=H_F$ gives the second minimum.
\end{proof}

For a fixed Hermitian form $h_0$, put
$B_x=h_0^{-1}h_x$ and
$\kappa(x)=\lambda_{\max}(B_x)/\lambda_{\min}(B_x)$. In an
$h_0$-orthonormal eigenbasis,
$$
 \|A\|_{h_x}^2
 =\sum_{i,j}\frac{\lambda_i}{\lambda_j}|A_{ij}|^2
 \le\kappa(x)\|A\|_0^2.
$$
Consequently the fixed Euclidean metric satisfies
\begin{equation}\label{grad:fixed-metric}
 g_E(v,\bar v):=\|dS(v)\|_0^2
 \ge\frac4{\kappa(x)}g_{\mathrm{GS}}(v,\bar v),
 \qquad
 \omega_E=i\partial\bar\partial\|S\|_0^2.
\end{equation}
The condition number $\kappa(x)$ need not be uniformly bounded.

Finally,
$$
 S=\sum_{p\in\mathcal P}(2p-w)\Pi_p(S),\qquad
 \Pi_p(S)\in\mathbb C[S]
$$
show that bounded grading coordinates are equivalent to uniformly
bounded spectral projections. Theorem~\ref{thm:grading-embedding}
gives a Euclidean embedding and polynomial recovery of the flag.
The bounded realization required by Griffiths'
Conjecture~\textup{(10.1)} \cite[p.~258]{Griffiths1970} requires
additional estimates.
\section{A quantitative estimate for the Hodge trace}\label{sec:trace}

Griffiths--Schmid constructed a proper exhaustion of $D$ with
positive Levi form in horizontal directions
\cite[Theorem~\textup{(8.1)}]{GriffithsSchmid1969}; potentials with
quantitative horizontal estimates were constructed in
\cite[Proposition~4.5.1 and Remark~4.5.3]{Eyssidieux1997}.
We give an explicit estimate for graded Hermitian vector spaces
and use it to construct bounded weights for the Bergman kernel
of square-integrable canonical forms.

Let $H_{\mathbb R}$ carry a real polarization $Q$, and put
$H_{\mathbb C}=H_{\mathbb R}\otimes_{\mathbb R}\C$. For the period
domain of a fixed Hodge type, write
$$
D=G_{\mathbb R}/V.
$$
Write $h_x$ for the positive Hodge metric at $x\in D$.
Here $G_{\mathbb R}$ preserves $Q$, $V$ fixes a reference
Hodge structure, and $K\supset V$ is the maximal compact subgroup
fixing its Hodge metric.

Fix $o\in D$ and a flat real $h_o$-orthonormal basis
$e_1,\ldots,e_N$, $N=\dim_{\mathbb C}H_{\mathbb C}$. Define
\begin{equation}\label{eq:trace}
E(x)=\sum_{\nu=1}^N h_x(e_\nu,e_\nu)
=\tr(h_o^{-1}h_x).
\end{equation}
The reference metric and frame remain fixed.

Let $r$ be the number of integers between the smallest and largest
indices $p$ with $h^{p,w-p}\ne0$, including both endpoints, and put
$$
c_r=\cos\!\left(\frac{\pi}{r+1}\right),\qquad
b_r=c_r^{-1}-1\quad(r\ge2).
$$

A reduced complex space is hyperconvex if it admits a bounded
negative plurisubharmonic exhaustion. Strict plurisubharmonicity
means that locally, in ambient coordinates, subtracting
$\varepsilon\|z\|^2$ preserves plurisubharmonicity for some
$\varepsilon>0$.

\begin{theorem}\label{thm:trace}\label{thm:trace-complete}
Let $X\subset D$ be a reduced, locally closed, positive-dimensional
analytic subspace with horizontal regular locus. Then $r\ge2$,
$E\ge N$, and
\begin{equation}\label{intro:trace}
 |\partial E(\xi)|^2\le c_r E\,
 (\partial\bar\partial E)(\xi,\bar\xi),
 \qquad \xi\in TX_{\mathrm{reg}}.
\end{equation}
For $0<a<b_r$, the function $-E^{-a}$ is bounded, negative,
continuous, and strictly plurisubharmonic on $X$. If $X$ is
closed in $D$, it is an exhaustion; hence $X$ is hyperconvex
and Stein. If $X$ is smooth,
$$
 \omega_{\mathrm{tr}}=i\partial\bar\partial\log E
$$
is K\"ahler with bounded primitive $d^c\log E$, where
$d^c=\frac i2(\bar\partial-\partial)$. It is complete when
$X$ is also closed in $D$.
\end{theorem}

\begin{lemma}\label{lem:graded}
Let $H=\bigoplus_{j=1}^rH_j$ be an orthogonal decomposition of a
Hermitian vector space, and let $A$ be an endomorphism satisfying
$A(H_j)\subset H_{j-1}$, where $H_0=0$. Then
$$
|\langle Av,v\rangle|^2
\leq \frac{c_r}{2}\|v\|^2
\left(\|Av\|^2+\|A^*v\|^2\right)
$$
for every $v\in H$.
\end{lemma}

\begin{proof}
Write $v=\sum v_j$, $a_j=\|v_j\|^2$, and
$b_j=\langle Av_j,v_{j-1}\rangle$. Orthogonality and
Cauchy--Schwarz give
$$
q:=\|Av\|^2+\|A^*v\|^2\geq
\sum_{j=2}^r |b_j|^2
\left(\frac{1}{a_j}+\frac{1}{a_{j-1}}\right),
$$
where terms with zero denominators are omitted, since their
$b_j$ vanish. Hence
\begin{align*}
|\langle Av,v\rangle|^2
&=\left|\sum_{j=2}^r b_j\right|^2
\leq q\sum_{j=2}^r\frac{a_ja_{j-1}}{a_j+a_{j-1}}\\
&\leq \frac{q}{2}\sum_{j=2}^r\sqrt{a_ja_{j-1}}.
\end{align*}
The path adjacency matrix on $r$ vertices has largest eigenvalue
$2c_r$. Applied to $u_j=\sqrt{a_j}$, this yields
$$
\sum_{j=2}^r u_ju_{j-1}
\leq c_r\sum_{j=1}^r u_j^2
=c_r\|v\|^2.
$$
\end{proof}

The same constant applies to a finite family by taking a direct sum.

On a polarized variation, the flat connection decomposes as
$$
\nabla=D_h+\theta+\theta^*.
$$
Here $D_h$ preserves the Hodge summands, $\theta$ is the
$(1,0)$ Higgs field and $\theta^*$ its adjoint. Put
$A=\theta(\partial_z)$; it lowers the
grading by one and satisfies $D_{h,\bar z}A=0$.

\begin{lemma}\label{lem:flat-derivatives}
For a flat section $v$, using a Hermitian product linear in the
first variable,
$$
\partial_z\|v\|_h^2=-2\langle Av,v\rangle_h,\qquad
\partial_z\partial_{\bar z}\|v\|_h^2
=2\left(\|Av\|_h^2+\|A^*v\|_h^2\right).
$$
\end{lemma}

\begin{proof}
Flatness gives
$$
D_{h,z}v=-Av,\qquad D_{h,\bar z}v=-A^*v.
$$
Metric compatibility gives the first identity, and the formula
$$
\partial_{\bar z}\langle Av,v\rangle_h
=-\langle AA^*v,v\rangle_h-\langle Av,Av\rangle_h
$$
gives the second.
\end{proof}

For the fixed frame, set
$$
\mu=\sum_{\nu=1}^N\langle Ae_\nu,e_\nu\rangle_h,\qquad
q=\sum_{\nu=1}^N
\left(\|Ae_\nu\|_h^2+\|A^*e_\nu\|_h^2\right).
$$
Then
\begin{equation}\label{eq:trace-derivatives}
E_z=-2\mu,\qquad E_{z\bar z}=2q.
\end{equation}
Lemma~\ref{lem:graded}, applied to the direct-sum vector
$(e_1,\ldots,e_N)$, gives
\begin{equation}\label{eq:trace-gradient}
|E_z|^2\leq c_r E E_{z\bar z}.
\end{equation}
In a nonzero horizontal tangent direction, $A\ne0$. Since the
frame spans each fiber, $q>0$ and $E_{z\bar z}>0$.

\begin{lemma}\label{lem:proper-trace}
The function $E:D\to\R$ is proper and satisfies $E\geq N$.
\end{lemma}

\begin{proof}
The homogeneous map defined by the Hodge metric,
$$
D=G_{\mathbb R}/V\longrightarrow G_{\mathbb R}/K,\qquad x\longmapsto h_x
$$
is proper, with compact fiber $K/V$. Its target is the space of
positive metrics compatible with $Q$.

For $B=h_o^{-1}h_x$, polarization invariance gives $\det B=1$.
Thus $\tr B\leq R$ implies, for every eigenvalue $\lambda_j$,
$$
\lambda_j\leq R,\qquad
\lambda_j\geq R^{-(N-1)}.
$$
Compatibility with $Q$ is a closed matrix condition. These
eigenvalue bounds therefore give compact trace sublevels on
$G_{\mathbb R}/K$, and hence on $D$. The arithmetic--geometric
mean inequality gives $E\geq N$.
\end{proof}

\begin{lemma}\label{lem:singular-horizontal-psh}
Let $u$ be a smooth real function on $D$ whose Levi form is
positive on every nonzero horizontal tangent vector. If
$X\subset D$ is a reduced locally closed analytic subspace with
horizontal regular locus, then $u|_X$ is continuous and strictly
plurisubharmonic in the sense defined above.
\end{lemma}

\begin{proof}
On a relatively compact neighborhood $U\subset D$ with coordinates $z$,
compactness of the horizontal unit sphere gives $c_U>0$ with
$$
i\partial\bar\partial u(\xi,\bar\xi)\geq
c_U\,i\partial\bar\partial\|z\|^2(\xi,\bar\xi)
$$
for horizontal $\xi$. For a proper resolution
$q:\widetilde X\to X\cap U$, put $f=i\circ q$, where
$i:X\cap U\hookrightarrow U$. The map
$$
T\widetilde X\xrightarrow{df}f^*TD
\longrightarrow f^*(TD/T_{\mathrm{hor}}D)
$$
vanishes on a dense open set, hence everywhere. For
$0<\varepsilon<c_U$,
$$
q^*\bigl(u|_X-\varepsilon\|z|_X\|^2\bigr)
$$
is plurisubharmonic.

For continuous functions, plurisubharmonicity descends under a proper
surjection. By the analytic-disc criterion
\cite[Theorem~5.3.1]{FornassNarasimhan1980}, it suffices to base change
by a disc and normalize a local curve dominating it. Subharmonicity
descends from this finite
branched cover off the branch values, and across them by continuity.
Thus $u|_X-\varepsilon\|z|_X\|^2$ is plurisubharmonic, which proves
strict plurisubharmonicity on $X$.
\end{proof}

\begin{proof}[Proof of Theorem \ref{thm:trace}]
Lemma~\ref{lem:proper-trace} gives $E\geq N$.
Equation~\eqref{eq:trace-gradient} holds on $X_{\mathrm{reg}}$.
Every horizontal vector of $D$ is tangent to a local horizontal
holomorphic curve, so the estimate holds in all horizontal directions.

Set $\Psi=\log E$. From \eqref{eq:trace-gradient},
\begin{equation}\label{eq:psi-strict}
\Psi_{z\bar z}
=\frac{E_{z\bar z}}{E}-\frac{|E_z|^2}{E^2}
\geq (1-c_r)\frac{E_{z\bar z}}{E}>0.
\end{equation}
The same estimate gives
\begin{equation}\label{eq:self-bounded}
\Psi_{z\bar z}\geq b_r|\Psi_z|^2.
\end{equation}
For $0<a<b_r$, put $\rho_a=-E^{-a}$. Then
\begin{align*}
(\rho_a)_{z\bar z}
&=aE^{-a-2}\left(EE_{z\bar z}-(a+1)|E_z|^2\right)\\
&\geq aE^{-a-1}\left(1-(a+1)c_r\right)E_{z\bar z}>0.
\end{align*}
By Lemma~\ref{lem:singular-horizontal-psh}, both potentials are
continuous and strictly plurisubharmonic on all of $X$.
The function $\rho_a$ is bounded and negative.
If $X$ is closed in $D$, properness of $E$ makes the sublevels
$\{\rho_a<c\}$ relatively compact for every $c<0$.
The function $\Psi=\log E$ is a continuous strictly
plurisubharmonic exhaustion; hence $X$ is Stein
\cite[Corollary~5.3.2]{FornassNarasimhan1980}.

If $X$ is smooth, $\omega=i\partial\bar\partial\Psi>0$ and
$$
\omega\geq b_r\,i\partial\Psi\wedge\bar\partial\Psi.
$$
Thus $\Psi$ is Lipschitz for $\omega$. If $X$ is also closed,
every $\omega$-Cauchy sequence lies in a compact $\Psi$-sublevel
set and converges, proving completeness. Finally,
$dd^c\Psi=\omega$, and the same inequality bounds
$\|d^c\Psi\|_\omega$.
\end{proof}

\begin{remark}\label{rem:trace-closedness}
Closedness is essential for exhaustion and completeness. On a
punctured coordinate disc in the weight-one period disc, the trace
extends smoothly across the puncture, so it is not proper and its
metric is incomplete. Nor is $\Delta^*$ hyperconvex: a bounded
negative subharmonic exhaustion would extend across $0$ with
value $0$, contradicting the maximum principle. The removable-singularity
theorem for $L^2$ canonical forms also identifies its Bergman
metric with the disc metric restricted to $\Delta^*$, which is incomplete.

For $X\subset D$ with its induced topology, $E|_X$ is proper
exactly when $X$ is closed: a sequence converging in $D$ has
bounded trace, and properness forces its limit into $X$.
For singular closed $X$, the exhaustion is defined on all of $X$,
whereas a path in $X_{\mathrm{reg}}$ may approach a singular point
at finite length.
\end{remark}

\section{Stein neighborhoods of the lifted period image}
\label{stein:section}

Griffiths' extension theorem and Siu's neighborhood theorem provide
a Stein neighborhood of the lifted period image, on which
Theorem~\ref{thm:grading-embedding} applies.

\begin{theorem}\label{stein:neighborhood}
Let $B$ be a smooth connected complex quasi-projective variety
carrying a polarized integral variation of pure Hodge structure.
Fix a marking and the lifted period map
$$
\widetilde\Phi:\widetilde B\longrightarrow D,\qquad
I=\widetilde\Phi(\widetilde B).
$$
The closure $Z=\overline I^{\,D}$ in $D$ is an analytic subspace,
which we equip with its reduced structure. It is Stein, and its regular
locus is horizontal. For every
open neighborhood $V$ of $Z$ in $D$, there is an
open Stein neighborhood $\Omega$ with
\begin{equation}\label{stein:containment}
I\subset Z\subset\Omega\subset V\subset D.
\end{equation}
The grading torsor admits a holomorphic section on $\Omega$, and
hence a holomorphic embedding
\begin{equation}\label{stein:grading}
S:\Omega\hookrightarrow\mathcal A\subset\mathfrak g,\qquad
\varpi\circ S=\operatorname{id}_{\Omega}.
\end{equation}
\end{theorem}

\begin{proof}
Passing to a connected finite \'etale cover with neat monodromy
does not change the universal cover or $I$. We may therefore assume
that
$$
\Gamma=\rho(\pi_1(B))
$$
is neat. It acts freely and properly on $D$, since point stabilizers
in $G_{\R}$ are compact. Thus
$$
D\longrightarrow\Gamma\backslash D
$$
is a holomorphic covering.

In a smooth projective compactification $\overline B$ with simple
normal crossings boundary, remove the divisors with infinite local
monodromy. The resulting open set $B^+$ has trivial local boundary
monodromy, since $\Gamma$ is neat. Griffiths' extension theorem
gives a proper horizontal holomorphic map
\begin{equation}\label{stein:proper-extension}
\Phi^+:B^+\longrightarrow\Gamma\backslash D
\end{equation}
\cite[Theorem~9.5, Propositions~9.10--9.11 and \S11(a)]{Griffiths1970III}.

Form the analytic fiber product
$$
P^+=B^+\times_{\Gamma\backslash D}D,\qquad
f^+:P^+\longrightarrow D.
$$
Then $f^+$ is proper and $P^+\to B^+$ is a covering. The restriction
$$
P=P^+|_B\longrightarrow B
$$
is the connected monodromy cover, since $\rho$ acts transitively
on the fiber $\Gamma$ by translation. Hence
$$
f^+(P)=I.
$$
By Remmert's theorem, $f^+(P^+)$ is a closed analytic subspace,
which we equip with its reduced structure. Density of $P$ in
$P^+$ gives
\begin{equation}\label{stein:closure}
f^+(P^+)=\overline{f^+(P)}^{\,D}=\overline I^{\,D}=Z.
\end{equation}
The horizontal differential of $f^+$ surjects onto the image
tangent space at generic regular points. Horizontality then extends
to all of $Z_{\mathrm{reg}}$ by continuity. The trace exhaustion
of Theorem~\ref{thm:trace-complete} shows that $Z$ is Stein;
the constant-image case is immediate.

Siu's Stein neighborhood theorem
\cite[Main Theorem]{Siu1976}, applied to the closed Stein subspace
$Z$ of $V$, gives \eqref{stein:containment}.
Theorem~\ref{thm:grading-embedding} then gives \eqref{stein:grading}.
\end{proof}

\begin{corollary}\label{stein:restriction}
Under the hypotheses of Theorem~\ref{stein:neighborhood}, the matrix
coordinates of $S$ are holomorphic on an open neighborhood of the
entire lifted image. Their restrictions give a topological embedding
of $I$ with its subspace topology. The spectral projections
$$
\Pi_p(S)=\prod_{q\ne p}
\frac{S-(2q-w)I}{2(p-q)},
$$
where $p,q$ range over the degrees with nonzero Hodge numbers,
recover the Hodge filtration by
$$
F^p=\bigoplus_{q\ge p}\im\Pi_q(S).
$$
If $I$ is endowed with its induced locally closed analytic-subspace
structure, the restriction is a holomorphic embedding.
\end{corollary}

\begin{proof}
Restrict \eqref{stein:grading}. The inverse is the restriction of
the holomorphic map $\varpi$; the polynomial formulas recover
the Hodge filtration. When $I$ has the induced analytic-subspace
structure, the same argument applies locally in $\Omega$.
\end{proof}

\begin{example}\label{stein:nonstein-image}
The lifted image in Theorem~\ref{stein:neighborhood} need not be
Stein. Let $Y(\ell)$, for sufficiently large full level $\ell$,
be the fine moduli curve of elliptic curves, and retain the ordered
decomposition in the product family over $Y(\ell)^2$. Its lifted
weight-one period map identifies its image with the diagonal product
domain
$\mathfrak H\times\mathfrak H\subset\mathfrak H_2$.
Remove one point $b$ from $Y(\ell)^2$, and restrict this
smooth projective family to
$$
B=Y(\ell)^2\setminus\{b\}.
$$
Writing $p:\mathfrak H^2\to Y(\ell)^2$, its lifted image is
$$
I=\mathfrak H^2\setminus p^{-1}(b).
$$
The removed set is nonempty, closed, and discrete. Removing these
points does not change the fundamental group in complex dimension
two, so this description agrees with the lifted period image of
the indicated source.

Hartogs extension gives
$\mathcal O(I)=\mathcal O(\mathfrak H^2)$.
Choose a small closed Euclidean ball around one removed point,
containing no other removed points, and let $K\subset I$ be its
boundary sphere. Every function in $\mathcal O(I)$ extends across
that ball. The maximum principle therefore puts the punctured
ball in the holomorphic hull of $K$ relative to $I$. This hull
contains a sequence approaching the missing point and is not
compact in $I$. Thus $I$ is not Stein. Its closure
$Z=\mathfrak H^2$, however, is Stein, and the grading construction
on a Stein neighborhood applies.
\end{example}

The neighborhood $\Omega$ need not be horizontal or
$\Gamma$-invariant, and the section $S$ need not be
$\Gamma$-equivariant.
\section{Exterior-power trace potentials and boundary asymptotics}\label{sec:boundary-hodge}

We construct a trace potential on exterior powers whose Levi form is
comparable to the Griffiths--Schmid metric, and compute its normal and
tangential asymptotics at the boundary. Reference metrics are fixed
in a flat trivialization throughout.

\subsection{The exterior-power trace potential}

Let $n\ge2$ be the rank of the variation and $r\ge2$ the length
of its Hodge filtration, counting every index between the first and
last nonzero graded pieces.
Fix any positive Hermitian metric $h_0$ on the underlying
vector space $V$. Put
\begin{equation}\label{bdy:exterior-potential}
 E_k=\tr_{\Lambda^kh_0}(\Lambda^kh),\qquad
 \Phi=\sum_{k=1}^{n-1}\log E_k.
\end{equation}
All metrics on exterior powers are induced by the corresponding metrics
on $V$. Set
$$
 c_k=\cos\frac{\pi}{k(r-1)+2},\quad
 \delta_{n,r}=\min_{1\le k<n}\frac{1-c_k}{\binom nk},\quad
 b_*=\min_{1\le k<n}\frac{1-c_k}{c_k}.
$$

\begin{theorem}\label{bdy:exterior-comparison}
Normalize the Griffiths--Schmid metric by
$g_{\rm GS}(v,\bar v)=\|\theta_v\|_h^2$. Then
\begin{equation}\label{bdy:exterior-bounds}
 \delta_{n,r}g_{\rm GS}\le
 g_\Phi:=\partial\bar\partial\Phi\le2^n g_{\rm GS},
 \qquad
 g_\Phi\ge\frac{b_*}{n-1}\,
 \partial\Phi\otimes\bar\partial\Phi .
\end{equation}
On the base of the variation, these inequalities hold as semipositive
forms, with $\ker g_\Phi=\{v:\theta_v=0\}$. The constants are independent
of the reference metric $h_0$.
\end{theorem}

\begin{proof}
Fix a point and put $A=\theta_v$. On $\Lambda^kV$, let $A^{[k]}$
be the induced derivation. The trace estimate of
Section~\ref{sec:trace}, applied to this exterior-power variation, gives
\begin{equation}\label{bdy:exterior-one}
 |\partial_vE_k|^2\le c_kE_k(E_k)_{v\bar v},\quad
 (\log E_k)_{v\bar v}\ge(1-c_k)\frac{(E_k)_{v\bar v}}{E_k},
 \quad
 \partial\bar\partial\log E_k
 \ge\frac{1-c_k}{c_k}|\partial\log E_k|^2.
\end{equation}
The Hodge length of $\Lambda^kV$ is at most $k(r-1)+1$, so the stated
constant is valid even when some of its Hodge summands vanish.

Choose an $h$-orthonormal basis diagonalizing $h_0^{-1}h$, with
eigenvalues $\lambda_1\ge\cdots\ge\lambda_n>0$. This basis need not
respect the Hodge grading. For $I=\{i_1<\cdots<i_k\}$, put
$\lambda_I=\prod_{i\in I}\lambda_i$. The flat-section derivative
identity gives
$$
 E_k=\sum_{|I|=k}\lambda_I,\qquad
 (E_k)_{v\bar v}
 =2\sum_{|I|=k}\lambda_I
 \bigl(\|A^{[k]}e_I\|^2+\|(A^{[k]})^*e_I\|^2\bigr).
$$
The index set $I_k=\{1,\ldots,k\}$ has maximal weight, hence
$\lambda_{I_k}/E_k\ge\binom nk^{-1}$. If
$s_k=\sum_{j\le k}A_{jj}$, direct expansion of the wedge gives
$$
 \|A^{[k]}e_{I_k}\|^2+\|(A^{[k]})^*e_{I_k}\|^2
 =2|s_k|^2+
 \sum_{i\le k<j}(|A_{ij}|^2+|A_{ji}|^2)=:Q_k.
$$
Each off-diagonal entry occurs for at least one $k$. Also
$\tr A=0$, since $A$ lowers Hodge degree, and, with $s_0=s_n=0$,
$$
 \sum_j|A_{jj}|^2
 =\sum_j|s_j-s_{j-1}|^2\le4\sum_{k=1}^{n-1}|s_k|^2.
$$
Consequently $\sum_kQ_k\ge\frac12\|A\|_{\rm HS}^2$, and
\eqref{bdy:exterior-one} proves the lower bound.

The flat-section identity also gives
$(\log E_k)_{v\bar v}\le4\|A^{[k]}\|_{\rm HS}^2$. Counting diagonal
and single-index replacement terms proves
$$
 \|A^{[k]}\|_{\rm HS}^2
 =\binom{n-2}{k-1}\|A\|_{\rm HS}^2
  +\binom{n-2}{k-2}|\tr A|^2.
$$
Summing and using $\tr A=0$ proves the upper bound. The last
inequality in \eqref{bdy:exterior-one}, followed by Cauchy--Schwarz
over $k$, proves the gradient estimate.
\end{proof}

\begin{remark}\label{bdy:exterior-position}
Plurisubharmonic potentials whose Levi forms control the Hodge metric
and their gradients were constructed by Eyssidieux
\cite[Proposition 4.5.1 and Remark 4.5.3]{Eyssidieux1997}.
The explicit formula \eqref{bdy:exterior-potential} gives constants
independent of the reference metric and admits the boundary expansion
below.
The tensor-product trace identity identifies $\exp\Phi$ with a single
trace on $\bigotimes_{k=1}^{n-1}\Lambda^kV$.

The determinant of $h_0^{-1}h$ is constant. Hence every $E_k$ has a
positive lower bound, and $\Phi$ is proper on $D$, since $E_1$ is
proper. For $0<a<b_* /(n-1)$, the function $-e^{-a\Phi}$ is bounded,
negative and strictly plurisubharmonic in horizontal directions.
On a closed horizontal image it is an exhaustion. On a closed smooth
horizontal image, $g_\Phi$ is complete. These conclusions follow as in
Section~\ref{sec:trace}. Neither $\Phi$ nor $g_\Phi$ is generally
monodromy invariant.
\end{remark}

\subsection{Asymptotics for one-parameter degenerations}

Consider an admissible polarized real variation of pure weight $w$
on $\Delta^*\times S$, where $S$ is a small ball. After a finite
base change making the monodromy unipotent, write
\begin{equation}\label{bdy:normal-form}
 F(z,s)=e^{zN}e^{\Gamma(q,s)}F_\infty(s),\qquad
 q=e^{2\pi iz},\quad z=x+iy,\quad \Gamma(0,s)=0 .
\end{equation}
Here $N$ is the logarithm of monodromy, and its monodromy weight
filtration $W$ is centered at $w$. Assume that the limiting mixed
Hodge structures, polarized by $N$, have constant Hodge numbers and
constant ranks of the intersections of the Hodge and weight filtrations
along $S$.

Write $G_d=\operatorname{Gr}_{w+d}^WV$, $m_d=\dim G_d$, and list the
centered weights with multiplicity as
$$
 d_1\ge\cdots\ge d_n,\qquad \sum_i d_i=0.
$$
Define
\begin{equation}\label{bdy:kappa}
 a_k=\sum_{i=1}^kd_i,\qquad
 \kappa(W)=\sum_{k=1}^{n-1}a_k
 =\frac12\sum_{i<j}|d_i-d_j|.
\end{equation}
On $G_d$, use the fixed reference metric
$$
 h_{0,d}=(h_0|_{W_{w+d}})/W_{w+d-1}.
$$
This quotient metric depends only on $h_0$ and the fixed filtration
$W$.

Equip $G_d$ with the pure Hodge metric $h_d^{\rm sl_2}(s)$
normalized as follows. Split the limiting mixed Hodge structure over
$\R$, restrict the Hodge metric of $e^{iN}\widehat F_\infty(s)$
to the centered weight-$d$ subspace, and transport it to $G_d$.
This is the Hodge metric of a polarized
variation of pure Hodge structure:
the flat primitive decomposition
$$
 G_d=\bigoplus_{\substack{\ell\ge|d|\\\ell\equiv d\ (2)}}
 N^{(\ell-d)/2}P_\ell,\qquad
 P_\ell=\ker(N^{\ell+1}:G_\ell\to G_{-\ell-2}),
$$
identifies it with a direct sum of the primitive polarizations
$Q(u,N^\ell v)$, with fixed positive factors from the
$\operatorname{Sym}^\ell\R^2$ metric at $i$.
The factors depend only on $\ell$ and the position in the
$\mathfrak{sl}_2$-string.
In particular, $\det h_d^{\rm sl_2}$ is constant in a flat volume
frame. Put
$$
 J_{d,j}=\tr_{\Lambda^jh_{0,d}}(\Lambda^jh_d^{\rm sl_2}),\qquad
 \Phi_d=\sum_{j=1}^{m_d-1}\log J_{d,j},\qquad
 g_{\rm gr}=\sum_d\partial_S\bar\partial_S\Phi_d.
$$
Empty sums are zero.

\begin{theorem}\label{bdy:finite-part}
Uniformly for $s$ in a transverse compact set and $x$ in a bounded
interval,
\begin{equation}\label{bdy:finite-part-formula}
 \Phi-\kappa(W)\log y
 =C_0+\sum_d\Phi_d(s)+O_{C_s^2}(y^{-1/2}).
\end{equation}
More precisely, if $k=\sum_{e>d}m_e+j$, where
$1\le j\le m_d$, then
$$
 y^{-a_k}E_k\longrightarrow
 P_k(s):=\left(\prod_{e>d}J_{e,m_e}\right)J_{d,j}
 \quad\hbox{in }C_s^2.
$$
Moreover,
\begin{equation}\label{bdy:full-levi}
 g_\Phi\big|_{(y\partial_z,\partial_s)}
 =
 \begin{pmatrix}\kappa(W)/4&0\\0&g_{\rm gr}\end{pmatrix}
 +O(y^{-1/2}),\qquad
 \ker g_{\rm gr}=\bigcap_d\{v:(\theta_{G_d})_v=0\}.
\end{equation}
Here $\theta_{G_d}$ is the Higgs field of the graded variation
$G_d$.
\end{theorem}

The matrix is expressed in the rescaled smooth tangent frame
$(y\partial_z,\partial_s)$.

\begin{proof}
We use Schmid's nilpotent orbit and $\mathrm{SL}_2$-orbit theorems
\cite[Theorems 4.12, 5.13 and 6.6]{Schmid1973} and the splitting results
of Cattani--Kaplan \cite[Theorems 3.12--3.14]{CattaniKaplan1989}.
For general asymptotic estimates of Hodge norms, see also
\cite[Section 12]{KatoNakayamaUsui2008}.

Choose a fixed $h_0$-orthogonal splitting $V=\bigoplus_dV_d^0$
of $W$, and put $D_y=\operatorname{diag}(y^{d_i/2})$.
Write $F_\infty=e^{i\delta_s}\widehat F_s$, with real centered
grading $Y_s$. Then
$[Y_s,N]=-2N$, $[\delta_s,N]=0$, and every component of
$\delta_s$ has weight at most $-2$. The comparison $U_s$ between
the fixed splitting and the $Y_s$-splitting induces the identity on
each $G_d$; thus $U_s-I$ strictly lowers weight. The change of
splitting $U_s$ may be complex; the transformations
$y^{-Y_s/2}$ and $e^{xN}$ belong to $G_{\R}$ and preserve the
Hodge metric.

Put $\varepsilon=y^{-1/2}$ and $\tau=x/y$. In the fixed splitting,
$$
 N_0=U_s^{-1}NU_s
$$
is independent of $s$: its only blocks have weight $-2$, and are
the fixed maps on $\operatorname{Gr}^W$.
The normalized pure flag
$$
 e^{iN}\exp\!\bigl(i\operatorname{Ad}(y^{Y_s/2})\delta_s\bigr)
 \widehat F_s
$$
stays in a compact subset of $D$. Its pulled-back metric is a smooth
matrix $K(\varepsilon,s)$ with
$$
 K(0,s)=H_{\rm gr}(s)=\bigoplus_dH_d(s),\qquad
 K-H_{\rm gr}=\varepsilon^2K_2(\varepsilon,s).
$$
Derivatives of $K_2$ of any fixed order are bounded on smaller
compact parameter sets. Also
$$
 B(\varepsilon,s)=D_yU_s^{-1}D_y^{-1}=I+O(\varepsilon)
$$
is smooth at $\varepsilon=0$, because each strictly weight-lowering
entry gains a positive integral power of $\varepsilon$. The
nilpotent-orbit Hodge matrix has the exact expression
\begin{equation}\label{bdy:normalized-matrix}
 D_y^{-1}HD_y^{-1}
 =B^*e^{-\tau N_0^*}K(\varepsilon,s)e^{-\tau N_0}B.
\end{equation}
The stars here denote matrix conjugate transpose.

For the finite part, set $\tau=0$ in the limit and use the principal
minor formula
$$
 E_k=\sum_{|I|=k}y^{d_I}
       \det(D_y^{-1}HD_y^{-1})_I,\qquad d_I=\sum_{i\in I}d_i.
$$
The maximum of $d_I$ is $a_k$, attained by taking every index in
the blocks above $d$ and $j$ indices in the $d$-block.
The limits of the corresponding minors sum to $P_k(s)$. All other exponents
are lower by at least one.
The positive leading coefficient has a uniform lower bound on compact
sets, so logarithms may be taken in $C_s^2$.
In $\sum_k\log P_k$, each nontrivial $J_{d,j}$, $j<m_d$, occurs
once, and the determinant terms are constant. This proves
\eqref{bdy:finite-part-formula}. The identity for $\kappa$ follows
by summing the partial sums and using $\sum_i d_i=0$.

To compute the normal and mixed Levi terms, retain $\tau$. The highest
$Y$-weight in $\Lambda^kV$ is $a_k$. Its weight space
consists of highest-weight vectors for copies of the irreducible
$\mathfrak{sl}_2$-module of highest weight $a_k$. If $v$ belongs
to this space, the split metric satisfies
$$
 \|N_0^jv\|^2=j!\frac{a_k!}{(a_k-j)!}\|v\|^2,\qquad
 \|e^{-\tau N_0}v\|^2=(1+\tau^2)^{a_k}\|v\|^2.
$$
The first equality follows from the raising--lowering relations and
adjointness; different weights are orthogonal, giving the second.
Thus \eqref{bdy:normalized-matrix} and the minor formula yield
\begin{equation}\label{bdy:log-expansion}
 \log E_k
 =a_k\log y+\log P_k(s)+a_k\log(1+x^2/y^2)
   +\varepsilon r_k(\varepsilon,\tau,s),
\end{equation}
where $r_k$ is smooth with bounded fixed-order derivatives for
bounded $\tau$ and $s$ in a compact set.

The vector fields
$$
 L_x=y\partial_x=\partial_\tau,\qquad
 L_y=y\partial_y=-\tfrac12\varepsilon\partial_\varepsilon
                      -\tau\partial_\tau
$$
preserve the $O(\varepsilon)$ size of the remainder and its
transverse derivatives. Since
$y^2\partial_y^2=L_y^2-L_y$, it follows that the rescaled normal,
mixed and transverse second derivatives of the remainder are all
$O(\varepsilon)$. The leading normal term is
$$
 a_k\log y+a_k\log(1+x^2/y^2)
 =a_k\log|z|^2-a_k\log y,
$$
whose Levi coefficient is $a_k/(4y^2)$. Its mixed derivatives
vanish. This proves the matrix limit for the nilpotent orbit.

For the period map, $\Gamma(q,s)$ and its derivatives of any fixed order
are exponentially small. In $(\varepsilon,\tau)$ coordinates,
$$
 q=e^{-2\pi/\varepsilon^2}e^{2\pi i\tau/\varepsilon^2};
$$
all derivatives of any fixed order vanish to infinite order at
$\varepsilon=0$, even after multiplication by the finite powers
introduced by conjugation and scaling.
Hence \eqref{bdy:log-expansion} retains the same derivative estimates.
Finally Theorem~\ref{bdy:exterior-comparison} on each pure $G_d$
shows that $\ker\partial\bar\partial\Phi_d
=\{v:(\theta_{G_d})_v=0\}$.
Summing semipositive forms proves the assertion about the kernel.
\end{proof}

\begin{corollary}\label{bdy:noncollapse}
Assume $N\ne0$ and that the joint differential of the graded period
maps is injective on a neighborhood of a transverse compact set.
Then, for sufficiently large $y$,
\begin{equation}\label{bdy:cusp-comparison}
 g_\Phi\asymp g_{\rm GS}\asymp
 \frac{|dz|^2}{y^2}+g_{\rm gr}
 \asymp
 \frac{|dq|^2}{|q|^2(\log|q|)^2}+g_S ,
\end{equation}
where $g_S$ is any fixed positive metric on that compact set.
The comparisons hold for all tangent vectors, including mixed
normal--transverse vectors. Their constants depend on the local
degeneration and the compact set.
\end{corollary}

\begin{proof}
If $N\ne0$, symmetry of the centered weights gives $\kappa(W)>0$.
The kernel assertion in Theorem~\ref{bdy:finite-part} makes the limiting
matrix positive definite, uniformly on compact sets. The result follows
by perturbation and \eqref{bdy:exterior-bounds}. The Griffiths--Schmid
metric is invariant under monodromy, so its estimate extends from a
bounded real strip to every $x\in\R$. Since the constants in
\eqref{bdy:exterior-bounds} are independent of the flat reference
metric, the same bounds hold for $g_\Phi$ on every strip.
\end{proof}

\begin{proposition}\label{bdy:charts}
Under the assumptions of Corollary~\ref{bdy:noncollapse}, the
normalized period maps around $p=(x+iy,s_0)$,
$$
 v\longmapsto F\bigl(x+y(i+v_0),s_0+v'\bigr),
$$
have a uniform holomorphic embedding radius. After applying an element
of $G_{\R}$ fixed on each chart, their images lie in a fixed
relatively compact coordinate neighborhood in $D$, their derivatives
of every prescribed order are uniformly bounded, and their
differentials have a uniform positive lower bound on a smaller ball.
\end{proposition}

These are charts of the selected parametrized branch; no separation
from other branches of the image is required.

\begin{proof}
On a fixed small complex ball in $v$, the imaginary part of
$x+y(i+v_0)$ is comparable to $y$, and $s_0+v'$ stays in a
larger compact set. Thus \eqref{bdy:cusp-comparison} holds throughout
the ball, with respect to the Euclidean metric in $v$.
Translate its center to a fixed point $o\in D$ by a single element
of $G_{\R}$. The upper metric
bound controls the ambient length of every radial segment.
Shrinking the radius uniformly keeps all images in a fixed
relatively compact holomorphic coordinate neighborhood of $o$, by a
first-exit argument. Cauchy estimates give uniform second and
higher derivative bounds. The lower metric bound at the center
gives a positive lower bound for the least singular value of the
coordinate differential $A_p$.

Choose a complex-linear left inverse $L_p$ of $A_p$, with
uniformly bounded norm. For the centered coordinate map $f_p$,
the map $G_p=L_pf_p$ satisfies $dG_p(0)=I$. The second derivative
bound gives $\|dG_p-I\|\le1/2$ on a common smaller ball. Integrating
along line segments yields
$$
 |G_p(v)-G_p(w)|\ge\tfrac12|v-w|.
$$
Hence $G_p$, and therefore the period map, is injective there,
with a uniform inverse Lipschitz bound. The quantitative inverse
function theorem also supplies a fixed image ball and a holomorphic graph
description.
\end{proof}

\section{Boundary asymptotics for nodal quintic threefolds}\label{bdy:nodal-section}

Over the locus of smooth quintic threefolds, a local slice to the action
of projective coordinate changes carries a weight-three variation on
$H^3$, of rank $204$ and Hodge numbers $(1,101,101,1)$.
We verify the hypotheses
of Corollary~\ref{bdy:noncollapse} in all deformation directions near
a general one-nodal quintic and prove the corresponding two-parameter
estimate near a general two-nodal quintic.

\subsection{The period map on the equisingular stratum}

\begin{proposition}\label{bdy:nodal-residue}
For a general quintic threefold with $\nu=1$ or $2$ distinct
ordinary double points and no other singularities, a local deformation
slice is
$\Delta^\nu\times S$, with $\dim S=101-\nu$.
The vanishing cycles are independent, and the weight-three graded
piece $M$ of the limiting mixed Hodge structure has Hodge numbers
$$
 h(M)=(1,101-\nu,101-\nu,1).
$$
The differential of its period map is injective in all $101-\nu$ equisingular
directions. In the one-node case, the logarithmic Higgs map has rank
$101$, while its component acting on $F^3$ is locally equivalent
to $\operatorname{diag}(I_{100},q)$.
\end{proposition}

\begin{proof}
Write $S_{\rm pol}=\C[x_0,\ldots,x_4]$, and let $J(f)$ be
the ideal generated by the partial derivatives of the defining
quintic polynomial $f$.
Quintic polynomials independently prescribe first and second jets
at one or two fixed points. A general polynomial with zero first
jets and nondegenerate quadratic parts has precisely these nodes,
by Bertini away from them. The smoothing differential is
$$
 (S_{\rm pol}/J(f))_5\longrightarrow\C^\nu,\qquad
 [g]\longmapsto (g(p_1),\ldots,g(p_\nu)).
$$
It is surjective. After taking a local slice to the action of projective
coordinate changes, the dimension is $101$, with independent
smoothing coordinates $q_j$.

The local Morse residue calculation for the global holomorphic
three-form gives
$$
 a_j(q,s)=\int_{\delta_j}\Omega=q_jU_j(q,s),\qquad U_j(0,s)\ne0.
$$
The factor $q_j$ follows by the rescaling $x=\sqrt{q_j}\,u$
in a four-variable quadratic equation. A constant linear relation
between the $\delta_j$ would give
$\sum c_jq_jU_j=0$; differentiating in the independent $q_j$
at the origin gives $c_j=0$. Their Milnor balls are disjoint,
so their intersections vanish. Picard--Lefschetz therefore gives
independent rank-one $N_j$, with $N_iN_j=0$.

The vanishing-cycle exact sequence identifies
$H^3(X_{0,s})$ with $\bigcap_j\ker N_j$.
The extreme weight pieces are $\R(-1)^\nu$ and $\R(-2)^\nu$;
the limiting cohomology has rank $204$.
Consequently
$$
 \dim F^3H^3(X_{0,s})=1,\qquad
 \dim F^2H^3(X_{0,s})=102-\nu,\qquad
 \dim\operatorname{Gr}_F^2H^3(X_{0,s})=101-\nu.
$$
These ranks are constant on the equisingular stratum.
Its tangent space is
$I_{\{p_1,\ldots,p_\nu\}}(5)/J(f)_5$.
Differentiating the residue class spanning $F^3$ in a direction
$g$ gives
$$
 \operatorname{Res}\!\left(-g\Omega_{\mathbb P}/f^2\right)
          \pmod {F^3}.
$$
Put $U=\mathbb P^4\setminus X_{0,s}$. Restriction and residue give
the natural injection
\cite[Section~3.1.2 and Proposition~3.3(iv)]{WangNodal2017}
$$
 \operatorname{Gr}_F^2H^3(X_{0,s})
 \hookrightarrow\operatorname{Gr}_F^3H^4(U)
 \simeq(S_{\rm pol}/J(f))_5.
$$
The last isomorphism follows from \cite[Theorem~2.2]{DimcaSaitoWotzlaw2009}
with $n=4$, $d=5$, and
$q=1<m=\lfloor n/2\rfloor=2$, so its unconditional
Jacobian-ring formula applies. For the equisingular family,
\cite[Corollary~4.5]{DimcaSaitoWotzlaw2009} identifies the
Gauss--Manin map
$\operatorname{Gr}_F^4H^4(U)\to\operatorname{Gr}_F^3H^4(U)$
with multiplication by $-g$, from $(S_{\rm pol}/J(f))_0$
to $(S_{\rm pol}/J(f))_5$. Its additional hypotheses are vacuous
for $q=1<m=2$. Naturality of restriction and residue therefore
sends the derivative above to $-[g]$ under this injection.
Thus a zero derivative forces $g\in J(f)_5$, proving injectivity.
Since $W_2$ has only type $(1,1)$, the differential of $F^3$
remains injective after passage to $M$.

For one node, the logarithmic Higgs field in the normal direction
at $q=0$ is induced by $N/(2\pi i)$. Its component
$F_\infty^2/F_\infty^3\to F_\infty^1/F_\infty^2$ is nonzero,
while its component on $F_\infty^3$ vanishes. The $100$
transverse columns are linearly independent on $F_\infty^3$, so
together with the normal column they give rank $101$. In a frame of
the canonical extension, the normal column acting on $F^3$ is
represented by
$$
 q\partial_q\widehat\Omega+\frac{N}{2\pi i}\widehat\Omega=qB(q,s).
$$
Pairing with the vanishing cycle gives
$Q(\delta,B(0,s))=U(0,s)\ne0$, whereas this functional kills all
$100$ transverse columns and $F_\infty^3$. Thus the maximal
minor of this component vanishes exactly to first order.
An invertible $100$-minor permits row and column elimination,
giving $\operatorname{diag}(I_{100},q)$.
\end{proof}

\begin{corollary}\label{bdy:quintic-one-node}
At a general one-nodal quintic, the centered weights are
$1,0^{202},-1$, and
$$
 \kappa(W)=203,\qquad h(M)=(1,100,100,1),\qquad
 g_\Phi\big|_{(y\partial_z,\partial_s)}
 =\begin{pmatrix}203/4&0\\0&g_{\Phi_M}\end{pmatrix}
       +O(y^{-1/2}).
$$
Here $g_{\Phi_M}$ is the Levi form of the exterior-power trace
potential of $M$. The transverse block is positive definite in all
$100$ equisingular directions. Hence the $101$-dimensional period map satisfies
\eqref{bdy:cusp-comparison} and has the uniform embedded charts of
Proposition~\ref{bdy:charts}.
\end{corollary}

\begin{proof}
The two extreme graded pieces have rank one, so their exterior-power
potentials vanish. Apply Proposition~\ref{bdy:nodal-residue} and
Theorem~\ref{bdy:finite-part}.
\end{proof}

\subsection{Two nodes with arbitrary relative degeneration rates}

\begin{theorem}\label{bdy:quintic-corner}
For the general two-nodal quintic of
Proposition~\ref{bdy:nodal-residue}, write
$q_j=e^{2\pi iz_j}$, $y_j=\operatorname{Im}z_j$.
Let $g_M$ be the Griffiths--Schmid metric of the graded variation
$M$. On transverse compact sets and for all $y_1,y_2\ge Y$,
\begin{equation}\label{bdy:double-cusp}
 g_{\rm GS}\asymp g_\Phi\asymp
 \frac{|dz_1|^2}{y_1^2}+\frac{|dz_2|^2}{y_2^2}+g_M.
\end{equation}
The estimate holds in all $101=2+99$ deformation directions, with
constants independent of $y_1/y_2$. There are also uniform embedded
holomorphic period charts in the coordinates
$$
 z_j=x_j^0+y_j^0u_j,\quad |u_j-i|<r,\qquad s=s_0+v,
$$
with uniform bounds on derivatives of every fixed order after applying
an element of $G_{\R}$ fixed on each chart.
\end{theorem}

\begin{proof}
Choose a real symplectic basis with $N_jf_j=e_j$ for $j=1,2$.
At a reference point of $S$, on the complex Lagrangian
$F_M^2$, diagonalize the Hermitian form $iQ(v,\bar w)$ in a basis
$v_a$, with signed diagonal entries. Then
$$
 E_M=\operatorname{span}_{\R}\{v_a+\bar v_a\},\qquad
 Q(E_M,E_M)=0,\qquad
 F_M^2\cap E_{M,\C}=0.
$$
Lift $E_M$ to $W_3$ and adjoin $W_2$ to obtain a real
Lagrangian $E$. Since $F_\infty^2\cap W_{2,\C}=0$,
$$
 F_\infty^2\cap E_\C=0.
$$
Shrinking the parameter polydisc preserves this transversality. The
untwisted Hodge filtration therefore lies in a fixed Lagrangian graph
chart. Its symmetric matrix $T$ is holomorphic on the polydisc,
including the boundary, and the period map has the form
$$
 F^2=\operatorname{span}\{f_\alpha+\sum_\beta
                    \tau_{\beta\alpha}e_\beta\},\qquad
 \tau=\operatorname{diag}(z_1,z_2,0,\ldots,0)+T(q,s).
$$
The Hodge line $F^3$ is generated by
$\sum_\alpha a_\alpha(q,s)
(f_\alpha+\sum_\beta\tau_{\beta\alpha}e_\beta)$.
The first two coefficients pair with the vanishing cycles, so
$a_j=q_jU_j$ for $j=1,2$.

At a center $(x_j^0+iy_j^0,s_0)$, fix the real symplectic transformation
$$
 g_y e^{-x_1^0N_1-x_2^0N_2},\qquad
 g_ye_j=(y_j^0)^{-1/2}e_j,\quad
 g_yf_j=(y_j^0)^{1/2}f_j ,
$$
acting trivially on the central complement. Put
$D_y=\operatorname{diag}(\sqrt{y_1^0},\sqrt{y_2^0},1,\ldots,1)$.
In these holomorphic coordinates the graph is given by
$$
 \tau_p=\operatorname{diag}(u_1,u_2,0,\ldots,0)
                     +D_y^{-1}T(q,s)D_y^{-1},\qquad
 a_p=D_ya(q,s).
$$
On a fixed larger polydisc $|u_j-i|<2r<1/2$, these converge to
$$
 \tau_\infty=\operatorname{diag}(u_1,u_2,T_{MM}(0,s)),\qquad
 a_\infty=(0,0,a_M(0,s)).
$$
This limit belongs to $D$: it is the orthogonal sum of two
standard nodal blocks and the polarized Hodge structure $M(s)$.

Uniformly on that larger polydisc, the error is bounded by a constant
times
\begin{equation}\label{bdy:corner-error}
 (\min y_j^0)^{-1/2}
 +\sum_{j=1}^2\sqrt{y_j^0}
                  e^{-2\pi(1-2r)y_j^0}.
\end{equation}
The mixed graph entries gain $1/\sqrt{y_j^0}$; node--node entries
gain $1/\sqrt{y_1^0y_2^0}$; and the entries of $F^3$ in the
nodal directions are $\sqrt{y_j^0}q_jU_j$. Each exponential factor
is thus paired with the corresponding $y_j^0$, giving a bound
independent of $y_1^0/y_2^0$. Cauchy estimates on a smaller polydisc
give the same bound for derivatives of every fixed order. Choose a
nonzero coordinate of $F_M^3$ to trivialize the projective
coordinates; a finite cover of the transverse compact set bounds the
denominators away from zero. The estimates therefore hold in local
coordinates on the flag variety.

The two normal derivatives of the limit change distinct nodal
blocks. The differential of the period map of $M$ is injective in
the remaining $99$ directions. The combined differential has a
uniform positive lower bound on compact sets. Smoothness and
$G_{\R}$-invariance of the metric give \eqref{bdy:double-cusp} for all
tangent vectors, including mixed vectors. The quantitative left-inverse argument in
Proposition~\ref{bdy:charts} gives a uniform embedding radius.
Finally Theorem~\ref{bdy:exterior-comparison} supplies the same
comparison for $g_\Phi$.
\end{proof}

\section{Canonical Bergman metrics and monodromy covers}\label{berg:section}

Let $Y$ be a complex manifold of dimension $d>0$. The Hilbert
space of square-integrable holomorphic top forms is
$$
 \mathcal H_Y=H^0_{(2)}(Y,K_Y),\qquad
 \|s\|^2=\int_Y i^{d^2}s\wedge\overline s.
$$
This norm depends only on the complex structure. If
$(s_j)$ is an orthonormal basis and $s_j=f_j e$ in a local
holomorphic frame $e$ of $K_Y$, put $b_Y=\sum_j|f_j|^2$.
When evaluation is everywhere nonzero, these functions define a
Hermitian metric and its curvature form by
\begin{equation}\label{berg:kernel}
 h_Y(e,e)=b_Y^{-1},\qquad
 \omega_{B,Y}=i\Theta(K_Y,h_Y)
             =i\partial\bar\partial\log b_Y.
\end{equation}
The $L^2$ estimate for $\bar\partial$ gives interpolation of
finite jets of holomorphic top forms.

\begin{proposition}\label{berg:jets}
Suppose that $Y$ admits a complete K\"ahler metric and a smooth
bounded strictly plurisubharmonic function $\varrho$. For distinct
points $y_1,\ldots,y_s$ and integers $k_1,\ldots,k_s\geq0$, the map
$$
 \mathcal H_Y\longrightarrow
 \bigoplus_{i=1}^s K_{Y,y_i}/\mathfrak m_{y_i}^{k_i+1}K_{Y,y_i}
$$
is surjective, where $K_{Y,y_i}$ denotes the stalk.
Consequently $\omega_{B,Y}$ is a K\"ahler form,
and projective evaluation defines a holomorphic injective immersion
$$
 \kappa_Y:Y\longrightarrow\mathbb P(\mathcal H_Y^*),\qquad
 \kappa_Y^*\omega_{\rm FS}=\omega_{B,Y}.
$$
If $\varrho<0$ is also an exhaustion, then $\omega_{B,Y}$ is complete.
\end{proposition}

\begin{proof}
Fix disjoint coordinate balls at the prescribed points, local canonical
frames $e_i$, polynomial representatives $P_i$, and cut-offs
$\chi_i=1$ near their centers. Extend by zero and set
$$
 \eta=\sum_i\chi_iP_ie_i,\qquad v=\bar\partial\eta,\qquad
 \phi_\varepsilon=A\varrho+
       \sum_i(d+k_i)\chi_i\log(|z_i|^2+\varepsilon).
$$
The negative Levi terms of the cut-off logarithms lie in fixed compact
annuli. Thus an $A$ independent of $\varepsilon$ gives
\begin{equation}\label{berg:variable-bound}
 i\partial\bar\partial\phi_\varepsilon
 \geq\tfrac A2i\partial\bar\partial\varrho>0.
\end{equation}
Apply the $L^2$ estimate
\cite[Theorem~4.1 and Remark~4.5]{Demailly1982}, using the complete
K\"ahler metric as a background and
$i\partial\bar\partial\varrho$ in the estimate.
The form $v$ is supported in the fixed annuli, where both the
inverse curvature and the weight are uniformly bounded. Hence
$$
 \bar\partial u_\varepsilon=v,\qquad
 \int_Y i^{d^2}u_\varepsilon\wedge\overline{u_\varepsilon}
                 e^{-\phi_\varepsilon}\leq C.
$$
The weights are uniformly bounded above. Weak compactness and
monotone convergence give $\bar\partial u=v$ with the same estimate
for $\phi_0$. If a holomorphic coefficient of $u$ at $y_i$
had lowest nonzero Taylor degree $\ell\leq k_i$, its norm would contain
$$
 \int_0^\epsilon r^{2\ell-2k_i-1}\,dr=\infty.
$$
Thus $\eta-u\in\mathcal H_Y$ has the prescribed jets.

The mean-value and Cauchy estimates imply continuity of evaluation
and its derivatives. Interpolation at one point gives a positive
kernel; interpolation at two points makes projective evaluation
injective; and interpolation of first jets makes the pullback of the
Fubini--Study form positive definite. If $\varrho$ is a negative
exhaustion, $Y$ is hyperconvex, so completeness follows from
\cite[Theorem~1]{Chen2004}.
\end{proof}

\begin{corollary}\label{berg:horizontal}
Let $Y$ admit a complete K\"ahler metric, and let
$\Phi:Y\to D$ be a horizontal holomorphic immersion. Then the
jet interpolation and positivity conclusions of
Proposition~\ref{berg:jets} hold. If $\Phi$ is proper, the Bergman
metric is complete.
In particular, every closed smooth horizontal complex submanifold
$X\subset D$ has a complete canonical Bergman metric.
\end{corollary}

\begin{proof}
For $0<a<b_r$, Theorem~\ref{thm:trace} gives the smooth bounded
strictly plurisubharmonic function
$\varrho=-(E\circ\Phi)^{-a}$, whose Levi form is positive
definite because $d\Phi$ is injective.
If $\Phi$ is proper, properness of $E:D\to\R$ makes $\varrho$
an exhaustion. Proposition~\ref{berg:jets} applies.
For a closed smooth horizontal submanifold, the complete background
metric is supplied by Theorem~\ref{thm:trace-complete}.
\end{proof}

\begin{proposition}\label{berg:descent}
Let $M$ be a connected smooth quasi-projective variety of positive
dimension carrying a polarized integral variation of Hodge structure
whose period differential is everywhere injective. Let
$$
 \rho:\pi_1(M)\longrightarrow\Gamma\subset G_{\R},\qquad
 \pi:Y\longrightarrow M
$$
be its monodromy representation, with $\Gamma=\im\rho$, and the
connected covering corresponding to $\ker\rho$. Then
$\mathcal H_Y$ interpolates arbitrary finite jets of holomorphic top forms.
Its Bergman form and Hermitian metric on $K_Y$ descend uniquely to
$\omega_{M,\rho}$ and $h_{M,\rho}$, with
\begin{equation}\label{berg:descended}
 \pi^*\omega_{M,\rho}=\omega_{B,Y},\qquad
 \pi^*h_{M,\rho}=h_Y,\qquad
 i\Theta(K_M,h_{M,\rho})=\omega_{M,\rho}>0.
\end{equation}
\end{proposition}

These metrics depend only on the complex structure of the covering
$Y\to M$.

\begin{proof}
A smooth quasi-projective variety has a complete K\"ahler metric,
constructed using a Poincar\'e metric near the simple normal crossings
boundary of a smooth projective compactification. Its pullback to a
covering is complete. The lifted period map $\widetilde\Phi:Y\to D$ is a horizontal
holomorphic immersion, so Corollary~\ref{berg:horizontal} applies.

Every deck transformation $\gamma$ acts unitarily on
$\mathcal H_Y$ by $U_\gamma s=(\gamma^{-1})^*s$, because
integration of canonical volume densities is invariant under change
of variables. Hence the kernel, its metric on $K_Y$, and its
curvature are invariant. Since $K_Y=\pi^*K_M$, they descend and
satisfy \eqref{berg:descended}. Both metrics are defined by the canonical
$L^2$ norm on $Y$; the choices made in the interpolation argument
do not affect them.
\end{proof}

Projective evaluation is equivariant for the dual unitary action
$U_\gamma^\vee\ell=\ell\circ U_\gamma^{-1}$. It gives a holomorphic
section of $Y\times_\Gamma\mathbb P(\mathcal H_Y^*)\to M$.
Its injectivity implies that this projective representation of
$\Gamma$ is faithful. If $q:M'\to M$ is a connected finite
\'{e}tale covering defined by a subgroup of $\pi_1(M)$ containing
$\ker\rho$, the pulled-back variation has the same monodromy cover
$Y$. Consequently
$$
 \omega_{M',\rho'}=q^*\omega_{M,\rho},\qquad
 h_{M',\rho'}=q^*h_{M,\rho}.
$$

\begin{remark}\label{berg:scope}
Although the Bergman metric descends to $M$, the individual
canonical forms on $Y$ need not descend to $M$ or to the period
image. Completeness of $\omega_{M,\rho}$ and the topological embedding
property of $\kappa_Y$ also require further hypotheses. The positive
curvature in
\eqref{berg:descended} is that of $h_{M,\rho}$ on $K_M$;
the Ricci curvature of $\omega_{M,\rho}$ instead depends on the
Hermitian metric that $\omega_{M,\rho}$ induces on $K_M$.
\end{remark}
\section{Bergman estimates on closed period images}
\label{bberg:section}

For a reduced pure $n$-dimensional closed horizontal analytic
subspace $Z\subset D$, consider the Hilbert space of
square-integrable holomorphic $n$-forms on its regular locus:
\begin{equation}\label{bberg:space}
 \mathcal H_Z=H^0_{(2)}(Z_{\rm reg},K_{Z_{\rm reg}}),\qquad
 \|\sigma\|^2=\int_{Z_{\rm reg}}i^{n^2}\sigma\wedge\overline{\sigma}.
\end{equation}
The diagonal kernel defines the canonical density
$$
 dV_{B,Z}=i^{n^2}\sum_\nu\sigma_\nu\wedge\overline{\sigma_\nu}
$$
for an orthonormal basis. This density is not, in general, the
Riemannian volume density of the Bergman metric
$g_{B,Z}= \partial\bar\partial\log b_Z$.
We work on $Z_{\rm reg}$, using the global Hilbert space
$\mathcal H_Z$ throughout.

We use Guan--Zhou's point-extension theorem and Zimmer's method
of constructing jets and estimating the Bergman metric from uniform
charts and bounded plurisubharmonic weights \cite{Zimmer2020}.
Here the charts may meet singularities of the horizontal space $Z$.

\subsection{A lower bound for the canonical density}

\begin{theorem}\label{bberg:density}
For fixed $D$, $n\geq1$, and a normalization of the
Griffiths--Schmid metric, there is $c(D,n)>0$ such that every reduced pure
$n$-dimensional closed horizontal analytic space $Z\subset D$
satisfies
$$
 dV_{B,Z}\geq c(D,n)\,dV_{\rm GS}
 \quad\hbox{on } Z_{\rm reg}.
$$
\end{theorem}

\begin{proof}
Fix $o\in D$, and choose ambient holomorphic coordinates $w$
centered at $o$, whose differential is an isometry for the invariant
Hermitian metric $g_D$. Fix a coordinate ball of radius $r<1$,
and a smooth cut-off $\chi$, equal to one near zero and supported
inside this ball. For an exponent allowed by the trace theorem,
write
$$
 \rho_o=-E_o^{-a},\qquad E_o(F)=\tr(h_o^{-1}h_F).
$$
This function is bounded and negative on $D$, with positive Levi
form on nonzero horizontal vectors. The negative Levi terms of
$n\chi\log|w|^2$, extended by zero, occur
in one fixed compact annulus. Positivity of $dd^c\rho_o$ on the
horizontal unit sphere over this annulus gives a fixed $A>0$ for
which
$$
 u_o=n\chi\log|w|^2+A\rho_o
$$
is negative and has nonnegative Levi form on horizontal vectors.
Its only logarithmic pole is $o$. For $p\in Z_{\rm reg}$,
choose $g\in G_{\R}$
with $g(o)=p$, and put $u_p=u_o\circ g^{-1}$. All constants
remain unchanged.

The trace exhaustion makes $Z$ Stein. On the irreducible component
containing $p$, Cartan's Theorem A supplies
$f\in H^0(Z,\mathcal I_{Z_{\rm sing}})$ with $f(p)\ne0$.
The principal open set
$$
 W=Z\setminus\{f=0\}
$$
is a smooth Stein manifold: it is the closed subspace $tf=1$ of
$Z\times\C$. The restriction $\Psi=u_p|_W$ is negative
plurisubharmonic and smooth off $p$.
Choose local coordinates $z$ on $Z_{\rm reg}$ orthonormal for
$g_{\rm GS}(p)$. Then
$$
 w(g^{-1}F(z))=Tz+O(|z|^2),\quad T^*T=I,\qquad
 \Psi=n\log|z|^2-A(\dim H_\C)^{-a}+o(1).
$$
In particular, the constant term of the pole is independent of
the point, tangent plane, and second fundamental form.

Apply the point-extension case of
\cite[Definition~1.1 and Theorem~2.2]{GuanZhou2015} to the smooth
Stein pair $(W,\{p\})$, with the trivial line bundle and its
constant metric,
and the function $c_0(t)=1$. Its curvature condition is just
plurisubharmonicity of $\Psi$. For an auxiliary volume
$a(z)d\lambda_z$ and
a pole $\Psi=n\log|z|^2+b+o(1)$, the residue measure at $p$ is
$a(0)e^{-b}$, up to the fixed normalization. Multiplication by
the squared norm of the prescribed value in $K_{W,p}$ cancels
$a(0)$. A value with coefficient one therefore
extends to $\sigma_p\in H^0(W,K_W)$ with
$$
 \|\sigma_p\|^2_{L^2(W)}\leq C_n e^{A(\dim H_\C)^{-a}}.
$$
Square-integrable holomorphic top forms extend across a proper analytic
subset of a smooth manifold, with unchanged norm. Thus $\sigma_p$
extends across $\{f=0\}\cap Z_{\rm reg}$ to
$\mathcal H_Z$; on other components take zero on their regular
loci (their intersections lie in $Z_{\rm sing}$). The extremal
characterization of the Bergman kernel gives the stated lower bound.
\end{proof}

\begin{remark}\label{bberg:density-scope}
The extension across $\{f=0\}\cap Z_{\rm reg}$ preserves the
$L^2$ norm. Thus the kernel is defined on the regular locus,
without choosing an extension of the canonical bundle to
$Z_{\rm sing}$. Comparison of the Bergman and Griffiths--Schmid
metrics also requires estimates for first derivatives, which we
establish below.
\end{remark}

\subsection{Uniform interpolation and the Bergman metric}

Fix an invariant Hermitian metric $g_D$ on $D$ restricting to
$g_{\rm GS}$ on horizontal tangent spaces.

\begin{theorem}\label{bberg:uniform-jets}
Let $P\subset Z_{\rm reg}$. Suppose there are constants
$R>0$, $0<c_0\leq C_0$, and holomorphic embeddings
$$
 \kappa_p:B_R^n\hookrightarrow D,\qquad
 \kappa_p(0)=p,\qquad \kappa_p(B_R^n)\subset Z,\qquad p\in P,
$$
such that on the whole ball
\begin{equation}\label{bberg:chart-bound}
 c_0g_{\rm eucl}\leq\kappa_p^*g_D\leq C_0g_{\rm eucl}.
\end{equation}
Each chart has dimension $\dim Z$. For every fixed $k\geq0$
and every polynomial $Q$ of degree at most $k$, there is
$\sigma_{p,Q}\in\mathcal H_Z$ such that
$$
 j_0^k\left(\frac{\kappa_p^*\sigma_{p,Q}}
 {dz_1\wedge\cdots\wedge dz_n}\right)=j_0^k Q,\qquad
 \|\sigma_{p,Q}\|^2\leq C_k\|Q\|_{\rm coeff}^2 .
$$
The constants are uniform in $p$.
Moreover,
\begin{equation}\label{bberg:metric-comparison}
 c\,g_{\rm GS}\leq g_{B,Z}\leq C\,g_{\rm GS}
 \quad\hbox{at all points of }P.
\end{equation}
For every integer $m\geq0$, the norm of
$\nabla_{g_{B,Z}}^m R(g_{B,Z})$ is uniformly bounded on $P$.
\end{theorem}

Each chart parametrizes a smooth branch that may meet other branches
of $Z$; only its center is required to lie in $Z_{\rm reg}$.
The pullback of a form in $\mathcal H_Z$ extends across the inverse
image of the singular set by $L^2$ removability.

\begin{proof}
Choose $g_p\in G_\R$ with $g_p(o)=p$.
The upper bound in \eqref{bberg:chart-bound} controls the lengths
of radial paths in the parameter ball. After decreasing $R$
uniformly, a first-exit argument shows that
$g_p^{-1}\kappa_p(B_R)$ lies in a fixed relatively compact ambient
coordinate neighborhood $U$ of $o$.
The functions
$$
 \rho_p=\rho_o\circ g_p^{-1}
$$
are uniformly bounded and negative on all of $Z$, strictly
plurisubharmonic on $Z_{\rm reg}$, and satisfy
$$
 b\,\omega_{\rm eucl}\leq
 dd^c(\rho_p\circ\kappa_p)\leq B\,\omega_{\rm eucl}
$$
on the entire ball, by compactness of the horizontal unit sphere
over $\overline U$ and
\eqref{bberg:chart-bound}. The chart is horizontal even where it
meets $Z_{\rm sing}$: its differential modulo the holomorphic
horizontal distribution vanishes on a dense open subset and hence
everywhere.

As in the preceding proof, choose
$f_p\in H^0(Z,\mathcal I_{Z_{\rm sing}})$ with $f_p(p)\ne0$
and let $W_p=Z\setminus\{f_p=0\}$.
Put $A_p=\kappa_p^{-1}\{f_p=0\}$. This is a proper analytic set
not containing zero. Since the chart has dimension $\dim Z$,
$\kappa_p(B_R\setminus A_p)$ is open in the smooth space $W_p$.
Cut-off functions supported in a fixed compact subset of $B_R$
extend smoothly by zero to $W_p$. Indeed, a boundary accumulation
point of their support lying in $W_p$ is still in the image of
that compact subset; the equal-dimensional inverse function
theorem then puts it inside the chart image. This argument also
applies when the deleted divisor meets the cut-off annulus.

Take fixed cut-offs $\chi,\chi_0$ equal to one near zero, and set
$$
 e_z=dz_1\wedge\cdots\wedge dz_n,\qquad
 \eta=\chi Q e_z,\qquad v=\bar\partial\eta ,
$$
and extend these forms by zero as above. For small $\delta>0$, put
$$
 \varphi_{p,\delta}
 =A_k\rho_p+(n+k)\chi_0\log(|z|^2+\delta)\leq0.
$$
The negative Levi part of the cut-off logarithm is supported in
a fixed annulus and uniformly bounded. Thus an $A_k$ independent of
$p$ and $\delta$ gives
$$
 \beta_{p,\delta}:=dd^c\varphi_{p,\delta}
 \geq \tfrac12 A_k dd^c\rho_p>0
 \quad\hbox{on }W_p.
$$
On the support of $v$, both the inverse curvature and
$e^{-\varphi_{p,\delta}}$ have uniform bounds in chart coordinates.

The Stein manifold $W_p$ admits a complete K\"ahler metric.
The $L^2$ estimate of
\cite[Theorem~4.1 and Remark~4.5]{Demailly1982} permits a different,
possibly incomplete, K\"ahler metric for the estimate; use
$\omega_p=dd^c\rho_p$.
For an $(n,1)$-form
$v=\sum_jv_{\bar j}e_z\wedge d\bar z_j$, the curvature-inverse
density is
\begin{equation}\label{bberg:curvature-cancel}
 \big\langle[\beta,\Lambda_{\omega_p}]^{-1}v,v\big\rangle_{\omega_p}
 dV_{\omega_p}
 =c_n\sum_{j,l}\beta^{j\bar l}v_{\bar j}
       \overline{v_{\bar l}}\,d\lambda_z .
\end{equation}
The auxiliary metric cancels from this density. The constant is
therefore independent of its behavior near the deleted divisor.
The form $v$ has finite weighted
$L^2$ norm, and the right side of the estimate is bounded by
$C_k\|Q\|_{\rm coeff}^2$. We obtain
$$
 \bar\partial u_{p,\delta}=v,\qquad
 \int_{W_p}i^{n^2}u_{p,\delta}\wedge\overline{u_{p,\delta}}\,
 e^{-\varphi_{p,\delta}}
 \leq C_k\|Q\|_{\rm coeff}^2.
$$
Since $\varphi_{p,\delta}\leq0$, this controls the unweighted
canonical norm. Weak compactness, lower semicontinuity for a
fixed regularized weight, and then monotone convergence give a
solution $u_p$ for the limiting singular weight.
Near zero the solution is holomorphic and integrable against
$|z|^{-2(n+k)}$. Its Taylor coefficients through order $k$
therefore vanish. The section $\eta-u_p$ has the prescribed jet
and a uniform unweighted norm. Removability across
$\{f_p=0\}\cap Z_{\rm reg}$ gives a section in
$\mathcal H_Z$.

Conversely, every $\sigma\in\mathcal H_Z$ pulls back to
$B_R\setminus A_p$ with canonical norm no greater than its global
norm, because $\kappa_p$ is injective. It extends to the whole
ball by $L^2$ removability. The mean-value and Cauchy estimates give
$$
 |f_\sigma(0)|^2+|df_\sigma(0)|^2\leq C\|\sigma\|^2,\qquad
 \kappa_p^*\sigma=f_\sigma e_z .
$$
Let
$$
 b_p=\sup_{\|\sigma\|\leq1}|f_\sigma(0)|^2,\qquad
 J_p(\xi)=\sup_{\substack{\|\sigma\|\leq1\\ f_\sigma(0)=0}}
             |df_\sigma(0)\xi|^2 .
$$
Constant and linear interpolation and the preceding upper bounds give
$$
 0<c_b\leq b_p\leq C_b,\qquad
 c_J|\xi|^2\leq J_p(\xi)\leq C_J|\xi|^2.
$$
The formula
$\kappa_p^*g_{B,Z}(0;\xi)=J_p(\xi)/b_p$, together with
\eqref{bberg:chart-bound}, proves
\eqref{bberg:metric-comparison}.

To bound curvature derivatives, pull back the forms to the fixed
ball and extend them as above. On a smaller ball the off-diagonal
Bergman kernel is uniformly bounded by Cauchy--Schwarz and the mean-value
inequality. Cauchy estimates bound all of its fixed-order
derivatives. At the center, the kernel and the metric have the
positive lower bounds just proved. Differentiating $\log b_p$
and the curvature formula therefore bounds every fixed covariant
derivative of curvature. At a regular center these are exactly
the derivatives of the original metric, even when other points of the
chart lie in $Z_{\rm sing}$.
\end{proof}

\begin{remark}\label{bberg:method}
Passing to a smooth Stein open subset and extending by $L^2$
removability allows Zimmer's interpolation argument
\cite{Zimmer2020} to apply to charts meeting $Z_{\rm sing}$, while
still using $\mathcal H_Z$.
\end{remark}

\subsection{Boundary comparison under nondegeneracy assumptions}

\begin{theorem}\label{bberg:boundary}
Let an admissible one-parameter polarized variation of pure Hodge
structure of arbitrary weight have unipotent monodromy with
logarithm $N\ne0$, fixed limiting mixed Hodge numbers, and
jointly injective differentials of its graded period maps on a
neighborhood of a transverse compact set $K$.
Suppose its lifted period branch is full dimensional in a reduced
closed horizontal space $Z$.
At all centers of this branch lying in $Z_{\rm reg}$, for
$y=\operatorname{Im}z$ sufficiently large,
$$
 g_{B,Z}\asymp g_{\rm GS}\asymp
 \frac{|dq|^2}{|q|^2(\log|q|)^2}+g_K,\qquad q=e^{2\pi iz}.
$$
The constants are uniform on $K$ and in $\operatorname{Re}z$.

For the full effective deformation spaces of general one-node and
two-node quintic threefolds, the conclusions hold in all $101$
directions. In the two-node case they read
\begin{equation}\label{bberg:double-node}
 g_{B,Z}\asymp g_{\rm GS}\asymp
 \sum_{j=1}^2
 \frac{|dq_j|^2}{|q_j|^2(\log|q_j|)^2}+g_K ,
\end{equation}
uniformly for $y_1,y_2\geq Y$, with constants independent of $y_1/y_2$.
\end{theorem}

The comparisons are between Hermitian forms. The statement allows
arbitrary nilpotency order and nonsplit limiting mixed Hodge structures.

\begin{proof}
Corollary~\ref{bdy:noncollapse} and Proposition~\ref{bdy:charts} give
embedded period charts of fixed radius satisfying
\eqref{bberg:chart-bound}, in the coordinates
$$
 z=x+y(i+v_0),\qquad s=s_0+v'.
$$
Theorem~\ref{bberg:uniform-jets} applies, with the one-node quintic
case supplied by Corollary~\ref{bdy:quintic-one-node}.

In the two-node quintic case, Theorem~\ref{bdy:quintic-corner}
provides the corresponding charts in
$z_j=x_j+y_j(i+v_j)$ and $99$ transverse coordinates.
The error estimate \eqref{bdy:corner-error} is uniform in both rates,
so Theorem~\ref{bberg:uniform-jets} gives
\eqref{bberg:double-node}.
\end{proof}

\begin{proposition}\label{bberg:source-version}
Let $M$ be a smooth quasi-projective variety carrying a polarized
integral variation of pure Hodge structure whose period differential
is everywhere injective. Let $M_\rho$ be its monodromy cover.
If a full-dimensional boundary family in
Theorem~\ref{bberg:boundary} has a locally biholomorphic classifying
map to $M$, the same comparisons hold for
$$
 H^0_{(2)}(M_\rho,K_{M_\rho})
$$
at all the corresponding points of $M_\rho$. The resulting metric
comparisons descend to $M$.
\end{proposition}

\begin{proof}
The cover has a complete K\"ahler background pulled back from
$M$, and the period immersion pulls back the bounded strictly
plurisubharmonic function obtained from the trace. A normalized
parameter ball is simply connected and its classifying map lifts to $M_\rho$. If two
points of that lift were equal, their periods would be equal;
the embedded period chart makes the lift injective. It is
locally biholomorphic, hence an open coordinate chart.
The proof of Theorem~\ref{bberg:uniform-jets} applies without the
deletion of a divisor. Deck transformations act unitarily on
this Hilbert space, as in Proposition~\ref{berg:descent}.
\end{proof}

Similarly, when $Z$ is the closure of the entire lifted image,
the monodromy group preserves $Z$ and acts unitarily on
$\mathcal H_Z$. On a neat level the comparison descends to the
quotient of the corresponding boundary region in $Z_{\rm reg}$.
This uses the Hilbert space $\mathcal H_Z$, whereas
Proposition~\ref{bberg:source-version} uses canonical forms on
$M_\rho$.
The curvature derivative bounds in
Theorem~\ref{bberg:uniform-jets} hold on these boundary regions.
Parameter paths remaining in such a region and approaching a
cusp have infinite length, and the volume on a local monodromy
fundamental domain with transverse compact parameters is finite.
Both statements follow by integrating the cusp model in that region.

\subsection{A complete period surface with no Bergman upper bound}

\begin{lemma}\label{bberg:strictification}
If a smooth Stein manifold has a bounded negative continuous
plurisubharmonic exhaustion, it has a smooth bounded negative
strictly plurisubharmonic exhaustion.
\end{lemma}

\begin{proof}
Normalize the given exhaustion by $-1\leq\rho<0$, and choose a
smooth strictly plurisubharmonic function $\psi\geq0$.
Choose $a_j<c_j<b_j<0$ with $a_j\nearrow0$ and
$$
 0<\epsilon_j<
 \frac{b_j-c_j}{1+\sup_{\{\rho\leq b_j\}}\psi}.
$$
On $\{\rho<b_j\}$, let
$u_j=\max\{\rho,c_j+\epsilon_j\psi\}$, and put $u_j=\rho$
elsewhere. Near $\rho=b_j$, on the side $\rho<b_j$, the maximum
equals $\rho$, so the two expressions define a plurisubharmonic
function. We have $\rho\leq u_j<0$, and $u_j$
is strictly plurisubharmonic on $\{\rho<a_j\}$.
The uniformly convergent sum $u=\sum_j2^{-j}u_j$ is continuous,
strictly plurisubharmonic, satisfies $\rho\leq u<0$, and is an
exhaustion. Richberg smoothing with allowed error $-u/2$
\cite[Chapter~I, Theorem~5.21]{DemaillyCADG} gives
$\rho\leq u\leq r\leq u/2<0$, where $r$ is smooth and
strictly plurisubharmonic. Thus $r$ is a bounded negative exhaustion.
\end{proof}

\begin{proposition}\label{bberg:upper-counterexample}
There is a smooth quasi-projective surface carrying a polarized
integral variation of Hodge structure of weight one whose period
map is an immersion. Its lifted image $X$ is closed and smooth
in the period domain, with complete Griffiths--Schmid metric, but
$$
 \sup_{p\in X,\,v\ne0}\frac{g_{B,X}(v,\bar v)}
                                  {g_{\rm GS}(v,\bar v)}=\infty.
$$
\end{proposition}

\begin{proof}
Put $Y_0=\C\setminus\{0,1\}$, $a=1/2$, and
$$
 B=\{(t,s)\in\C^2:t,a+ts,a+s^2\in Y_0\},\qquad
 F(t,s)=(t,a+ts,a+s^2).
$$
This is a closed embedding into $Y_0^3$, with image and inverse
$$
 C=\{(\lambda_2-a)^2=\lambda_1^2(\lambda_3-a)\},\qquad
 t=\lambda_1,\quad s=(\lambda_2-a)/\lambda_1.
$$
The derivative of its defining equation in $\lambda_3$ is
$-\lambda_1^2\ne0$. Thus $B$ and $C$ are smooth.
The product of Legendre families
$$
 E_t\times E_{a+ts}\times E_{a+s^2},\qquad
 E_\lambda:\ y^2=x(x-1)(x-\lambda),
$$
is a smooth projective polarized abelian family over $B$.
Its period differential is injective: the first coordinate detects
$\partial_t$, and $t\ne0$ makes the second detect
$\partial_s$.

Let $\lambda:\mathfrak H\to Y_0$ be the modular lambda covering,
and let $X$ be a connected component of $(\lambda^3)^{-1}(C)$.
Path lifting for the covering identifies $X$ with a lifted
period image of this family. It is closed and smooth in $\mathfrak H^3$.
The block-diagonal inclusion $\mathfrak H^3\hookrightarrow
\mathfrak H_3$ is closed and totally geodesic, so $X$ is a
closed smooth horizontal surface in the weight-one period domain.
The complete product hyperbolic metric restricts to a complete
metric on the closed embedded submanifold $X$.
Writing the metric on $Y_0$ as $h(\lambda)|d\lambda|^2$,
the induced Griffiths--Schmid metric, up to a fixed normalization, is
$$
 g_{\rm GS}=h(t)|dt|^2+h(a+ts)|s\,dt+t\,ds|^2
                         +h(a+s^2)|2s\,ds|^2.
$$
Consequently
\begin{equation}\label{bberg:legendre-collapse}
 g_{\rm GS}(\partial_s,\overline{\partial_s})\big|_{s=0}
                         =h(a)|t|^2.
\end{equation}

Consider the smooth curve
$$
 R=\{(s,\tau_3)\in\C\times\mathfrak H:
                            s^2=\lambda(\tau_3)-a\}.
$$
Its projection to $\mathfrak H$ is finite and proper; hence it
is Stein and the pullback of a bounded negative subharmonic
exhaustion is a continuous plurisubharmonic exhaustion.
Lemma~\ref{bberg:strictification} gives a smooth strictly subharmonic
function $-1\leq r_R<0$. The global map
$$
 X\longrightarrow\mathfrak H\times R,\qquad
 (\tau_1,\tau_2,\tau_3)\longmapsto
 \left(\tau_1,\frac{\lambda(\tau_2)-a}{\lambda(\tau_1)},\tau_3\right)
$$
is locally biholomorphic, since $\lambda'\ne0$.

Choose $\lambda(\tau_{20})=\lambda(\tau_{30})=a$ and the component
containing the connected path $p_y=(iy,\tau_{20},\tau_{30})$.
On $X$, put
$$
 \rho_y=
 \left|\frac{\tau_1/y-i}{\tau_1/y+i}\right|^2-1+r_R(s,\tau_3),
 \qquad -2\leq\rho_y<0.
$$
These functions are smooth and strictly plurisubharmonic. On a
fixed small bidisc in $(u,s)$ about $(i,0)$, the inverse branches
$$
 \tau_1=yu,\qquad
 \lambda(\tau_2)=a+\lambda(yu)s,\qquad
 \lambda(\tau_3)=a+s^2
$$
give embedded charts in $X$ for all large $y$:
$\lambda(yu)\to0$ uniformly on a larger bidisc.
In these coordinates $\rho_y$ is independent of $y$, and its
Levi form has a fixed positive lower bound.

The cut-off and logarithmic-pole argument in the proof of
Theorem~\ref{bberg:uniform-jets}, now with these global weights and
charts, gives $\sigma_y=f_y\,du\wedge ds\in\mathcal H_X$ such that
$$
 \|\sigma_y\|^2\leq C,\qquad
 f_y(i,0)=0,\qquad \partial_sf_y(i,0)=1.
$$
Here $X$ is smooth Stein. The uniform constants depend only on
the bounded weights and their Levi bounds on the fixed cut-off
annulus.
Interpolation of values makes the kernel positive. Restricting global
forms to the fixed bidisc and applying the mean-value inequality gives
$b_y(i,0)\leq C'$. The derivative extremal is at least $C^{-1}$,
so
$$
 g_{B,X}(\partial_s,\overline{\partial_s})(p_y)
              =\frac{J_y(\partial_s)}{b_y(i,0)}\geq (CC')^{-1}.
$$
Since $t_y=\lambda(iy)\to0$, equation
\eqref{bberg:legendre-collapse} proves the assertion.
\end{proof}

\clearpage
\raggedbottom

\end{document}